\documentclass[reqno,11pt]{amsart}
\usepackage{amsfonts}
\usepackage{graphicx}
\usepackage{amsfonts,amsmath, amssymb,mathrsfs}
\usepackage{marginnote}
\usepackage{color}
\usepackage[numbers,sort&compress]{natbib} 

\numberwithin{equation}{section}

\newcommand{\dis}{\displaystyle}

\newcommand{\R}{\mathbb{R}}

\newcommand{\fb}{\mathfrak{b}}
\newtheorem{theorem}{Theorem}[section]
\newtheorem{corollary}[theorem]{Corollary}
\newtheorem{lemma}[theorem]{Lemma}

\newtheorem{remark}[theorem]{Remark}
\newtheorem{definition}[theorem]{Definition}

\makeatletter
\newcommand{\rmnum}[1]{\romannumeral #1}
\newcommand{\Rmnum}[1]{\expandafter\@slowromancap\romannumeral #1@}
\makeatother

\def\v{\varepsilon}

\def\t{\theta}
\def\k{\kappa}

\def\g{\gamma}

\def\r{\rho}
\def\s{\sigma}

\def\f{\frac}
\def\dd{{\rm d}}

\usepackage{tikz}
\newcommand{\lr}[1]{\left(#1\right)}

\begin{document}

\title[Compressible Euler Equations with nonlocal Interactions]{Global solutions of the one-dimensional compressible Euler equations with nonlocal interactions \\via the inviscid limit}

\author[J.A. Carrillo]{Jos$\acute{\mathrm{E}}$ A. Carrillo}
\address{J.A. Carrillo:\,Mathematical Institute, University of Oxford, Oxford OX2 6GG, UK}
\email{carrillo@maths.ox.ac.uk}

\author[G.-Q. Chen]{Gui-Qiang G. Chen}
\address{Gui-Qiang G.Chen:\, Mathematical Institute, University of Oxford, Oxford OX2 6GG, UK}
\email{chengq@maths.ox.ac.uk}

\author[D.~F. Yuan]{Difan Yuan}
\address{D.~F. Yuan:\, School of Mathematical Sciences, Beijing Normal University and Laboratory of
Mathematics and Complex Systems; Mathematical Institute, University of Oxford, Oxford OX2 6GG, UK}
\email{yuandf@amss.ac.cn; difan.yuan@maths.ox.ac.uk}

\author[E. Zatorska]{Ewelina Zatorska}
\address{E. Zatorska:\, Mathematics Institute, University of Warwick, Coventry CV4 7AL, UK }
\email{ewelina.zatorska@warwick.ac.uk}

\begin{abstract}
We are concerned with the global existence of finite-energy entropy solutions of the one-dimensional
compressible Euler equations with (possibly) damping, alignment forces, and nonlocal interactions: Newtonian repulsion and quadratic confinement.
Both the polytropic gas law and the general gas law are analyzed.
This is achieved by constructing a sequence of solutions of the one-dimensional compressible Navier-Stokes-type
equations with density-dependent viscosity under the stress-free boundary condition
and then taking the vanishing viscosity limit.
The main difficulties in this paper arise from the appearance of the nonlocal terms.
In particular, some uniform higher moment estimates for the compressible Navier-Stokes equations on expanding intervals with stress-free boundary conditions are obtained by careful design of the approximate initial data. 
\end{abstract}

\keywords{Euler equations, nonlocal interactions, compressible flows,
large data, finite energy, concentration,
Navier-Stokes equations, {\it a priori} estimate,  higher integrability, vanishing viscosity,
compactness framework, approximate solutions, free boundary.}
\subjclass[2010]{\, 35Q35, 35Q31, 35B25, 35B44, 35L65, 35L67, 76N10, 35R09, 35R35, 35D30, 76X05, 76N17}
\date{\today}
\maketitle

\section{Introduction}
Hydrodynamic models of collective behavior provide a comprehensive framework for characterizing
the behavior of vast assemblies of interacting individuals.
In most of the interesting cases, these models can only be formally derived from the particle-type
systems capturing precise interactions between the individuals.
In cases where full mathematical rigor is attainable,
it is usually established through
the mean-field limit techniques or the BBGKY hierarchies;
see \cite{Carrillo,Ha,Carrillo-Choi21}
and the classical references therein on the matter.
In this paper, we are interested in a specific
one-dimensional (1-D) example of such models, which captures
local repulsion, and nonlocal attraction and repulsion forces, as well as nonlocal alignment. 
This model is described by the compressible Euler equations (CEEs) with the corresponding interaction forces{\rm :}
\begin{align}\label{1.1}
\begin{cases}
\rho_t+m_x=0,\\[1mm]
m_t+(\frac{m^2}{\rho}+P(\rho))_x =\lambda m-\rho \partial_xW\ast\rho+\displaystyle\int_{\R}\varpi(x-y)\big(\rho(x) m(y)-\rho(y) m(x)\big)\,\dd y,
\end{cases}
\end{align}
for $t>0$ and $x\in \R$,
where $\rho\geq0$ is the density, $m$ is the momentum, $P$ is the pressure.

The nonlocal attraction-repulsion interaction forces are described by  potential $W$ of the form:
\begin{equation*}
W(x)=-|x|+\frac{x^2}{2}.
\end{equation*}
The long-range attraction between the individuals is captured by the quadratic confinement part $\frac{1}{2}x^2,$ while the short-range repulsion is described by the Newtonian part of the potential $-|x|.$

The consensus in velocities (the alignment) among individuals is described by the weight $\varpi$ that satisfies:
 $$\varpi\in L^{\infty}(\R),\quad  \varpi(x)=\varpi(-x), \quad  \varpi\geq0.$$

Finally, the linear term $\lambda m$ in our system  stands for damping, if the  coefficient $\lambda$ is non-positive.

\subsection*{The equation of state}
We consider a general pressure law
$P(\rho)$ satisfying  hypotheses ($\mathcal{H}$) formulated below:
\begin{itemize}
\item [($\mathcal{H}.1$)] The pressure function $P(\rho)\in C^1([0,\infty))\cap C^4((0,\infty))$ satisfies the hyperbolic and genuinely nonlinear conditions{\rm :}
\begin{equation}\label{pressure2}
		P'(\rho)>0,\quad 2P'(\rho)+\rho P''(\rho)>0\qquad\,\, \text{for }\rho>0.
	\end{equation}
	
\smallskip	
\item [($\mathcal{H}.2$)] There exist constants $\rho_{*}>0$, $\gamma_1\in (1,3)$, and $\kappa_1>0$ such that
	\begin{equation}\label{pressure3}
		P(\rho)=\kappa_{1} \rho^{\gamma_1}\big(1+\mathcal{P}_1(\rho)\big) \qquad \text { for } \rho \in[0, \rho_{*}),
\end{equation}
for some function $\mathcal{P}_1(\rho)\in C^4((0,\infty))$
satisfying 
\[|\mathcal{P}_1^{(j)}(\rho)|\leq C_{*}\rho^{\gamma_1-1-j} \qquad \text { for }\rho\in (0,\rho_{*}),\ \text { and } j=0,\cdots,4,\]
 where $C_{*}>0$ depends only on $\rho_{*}$. 

\smallskip
\item [($\mathcal{H}.3$)] There exist constants $\rho^{*}> \rho_{*}>0$, $\gamma_2\in (1,\gamma_{1}]$, and $\kappa_{2}>0$, such that
	\begin{equation}\label{pressure4}
		P(\rho)=\kappa_2\rho^{\gamma_2}\big(1+\mathcal{P}_2(\rho)\big) \qquad \text { for } \rho \in[\rho^{*},\infty),
	\end{equation}
for some  function $\mathcal{P}_2(\rho)\in C^{4}((0,\infty))$ satisfying
\[|\mathcal{P}_{2}^{(j)}(\rho)|\leq C^{*}\rho^{-\epsilon-j}\qquad \text { for }\rho\in [\rho^{*},\infty),\ \epsilon>0,\  \text { and } j=0,\cdots,4,\]
where $C^{*}>0$ depends only on $\rho^{*}$.
\end{itemize}
Without loss of generality, by scaling,  $\kappa_1$ may be chosen to be equal to $\kappa_1=\frac{(\gamma_1-1)^2}{4\gamma_1}$.

A special example of such an equation of state is the pressure-density relation for the {\emph{polytropic gases}}:
\begin{equation}\label{pressure1}
P(\r)=\kappa\r^{\gamma}\qquad \text{ for $\gamma\in(1,\infty)$},
\end{equation}
where $\gamma$ is the adiabatic exponent. In this case, we can choose $\kappa=\frac{(\gamma-1)^2}{4\gamma}$.

\subsection*{Initial data} We consider the Cauchy problem for \eqref{1.1} with the initial data{\rm :}
\begin{equation}\label{1.6}
(\rho, m)|_{t=0}=(\rho_0, m_0)(x)
\end{equation}
such that $\rho_0\in L^1_+(\mathbb{R})$ and $\frac{m_0}{\sqrt{\rho_0}}\in  L^2(\mathbb{R})$
with finite initial total mass and initial second moment{\rm :}
\begin{equation}\label{1.7}
M:=\int_{\mathbb{R}} \rho_0(x)\,\dd x<\infty, \qquad  M_2:=\int_{\R}x^2\rho_0(x)\,\dd x<\infty.
\end{equation}
We further assume that the initial data are of finite energy:
\begin{align}\label{2.1}
\mathcal{E}_0:=
 \int_{\mathbb{R}}\Big(\frac{1}{2} \Big|\frac{m_0}{\sqrt{\rho_0}}\Big|^2
+ \rho_0e(\rho_0)+\frac{1}{2}\rho_0 W\ast\r_0\Big)(x) \,\dd x<\infty,
\end{align}
where $e(\rho)$ is the internal energy:
\begin{equation*}
	e(\rho)=\int_0^\rho
 \frac{P(\s)}{\s^2} \,\dd \s.
\end{equation*}
In particular,  for the polytropic case \eqref{pressure1},
\begin{equation*}
e(\rho):=\frac{\kappa}{\gamma-1}\rho^{\gamma-1}.
\end{equation*}

\subsection*{Main objectives of the paper}
The goal of the present paper is to establish the global-in-time existence of finite-energy entropy solutions
of system \eqref{1.1} without restriction on the size of initial data. Our method of choice is
the vanishing viscosity limit for the strong solutions of the compressible
Navier-Stokes equations (CNSEs). More precisely, we construct the solution by means of a sequence of approximate problems:
\begin{align}\label{1.8}
\begin{cases}
 \dis \r^{\varepsilon}_t+ (\rho^{\varepsilon}u^{\varepsilon})_x=0,\\[1mm] \dis (\r^{\varepsilon}u^{\varepsilon})_t+\big(\rho^{\varepsilon}(u^{\varepsilon})^2+P(\rho^{\varepsilon})\big)_x= \varepsilon(\mu(\r^{\varepsilon})u^{\varepsilon}_x)_x+\lambda\r^{\varepsilon}u^{\varepsilon}-\r^{\varepsilon}\partial_xW\ast\r^{\varepsilon}\\[1mm]
\quad\quad\quad\quad\quad\quad\quad\quad\quad\quad\quad\quad\quad\quad\displaystyle+\rho^{\v}\int_{\R}\varpi(x-y)(u^{\v}(y)-u^{\v}(x))\rho^{\v}(y)\,\dd y,
\end{cases}
\end{align}
on bounded intervals: 
\begin{equation}\label{3.1}
\Omega_T^\varepsilon=\big\{(t,x)\ :\,  b^-_\varepsilon(t)\leq r\leq b^+_\varepsilon(t),\,0\leq t\leq T \big\},
\end{equation}
with the free boundaries $b^\pm_\varepsilon(t)$,  expanding to infinity as $\varepsilon\rightarrow 0^+$.

In the above system, the viscosity coefficient is equal to
$\varepsilon\mu(\rho)=\varepsilon\rho^{\alpha},$ $\alpha\geq0$,
and the parameter $\v\in (0,1]$ denotes the inverse of the Reynolds number.
A derivation of CNSEs \eqref{1.8} from the Boltzmann equations (without the nonlocal term)
may be found in
Liu-Xin-Yang \cite{Liu}, in which
the viscosity is not constant but depends on the temperature. This dependence 
can be translated into the dependence of the viscosity on the density for the isentropic flow.

\subsection*{The state of the art} There is a huge literature devoted to the study of the existence of solutions to CEEs
either via the analysis of the vanishing artificial viscosity limit (see \cite{Bianchini,Dafermos,D. Hoff-Liu} for instance),
or by constructing the finite difference scheme.
The global existence of solutions with large initial data in $L^{\infty}$
was first established by DiPerna \cite{R.J.DiPerna2}
for $\gamma=1+\frac{2}{2n+1}$ with $n\geq2$ integer
by using the artificial viscosity method for the density equation.
For the general interval $1<\gamma\leq\frac{5}{3},$
the global existence problem was solved by Ding-Chen-Luo \cite{Ding1985,Ding1989} and Chen \cite{Chen2} by approximation via the Lax-Friedrichs scheme. The case of adiabatic exponent $\gamma>\frac{5}{3}$
was covered by Lions-Perthame-Tadmor \cite{Lions P.-L.2} and
Lions-Perthame-Souganidis \cite{Lions P.-L.1} via the vanishing viscosity method.
We also refer to \cite{Li-Wang-2006} for blow-up results on CEEs and
to Chen-LeFloch \cite{Chen-LeFloch-2000,Chen-LeFloch-2003}
for relevant results on CEEs with general pressure law.

Concerning approximation by the CNSEs, the $L^{\infty}$ estimates for the approximate solutions
are not expected when the initial data only are of finite energy.
Therefore, we work with the finite-energy framework, which was first considered
by LeFloch-Westdickenberg \cite{Ph. LeFloch} for $ 1<\gamma<\frac{5}{3}$, and was late 
developed
by Chen-Perepelitsa \cite{Perepelitsa}
to the whole physical range of adiabatic exponents $\gamma>1.$
In particular, they established
a compensated compactness framework (based on the earlier work by Tartar \cite{L. Tartar} and Murat \cite{F. Murat})
and proved the vanishing viscosity limit of the solutions of the 1-D CNSEs 
to the corresponding finite-energy solutions
of CEEs for large initial data for all $\gamma>1$ in \cite{Perepelitsa}.
In this result, the initial density at the far-field is allowed to be positive, and the viscosity
is independent of the density, {\it i.e.,} $\alpha=0.$
This framework has been subsequently extended by Chen-Schrecker \cite{Chen-Schrecker-2018} and
Chen-Wang \cite{Chen-Wang-2018} to study spherically symmetric Euler equations.
More recently, Chen-He-Wang-Yuan \cite{Chen2021} established the global existence of finite-energy solutions
of the multi-dimensional Euler-Poisson equations for both compressible gaseous stars and plasmas
with large spherically symmetric initial data.
Adapting the approach developed in \cite{Chen2021}, He-Wang \cite{He} proved the vanishing viscosity limit
for the 1-D CNSEs.
For CEEs with general pressure law, Shrecker-Schulz \cite{Schrecker-Schulz-2019,Schrecker-Schulz-2020}
proved the vanishing viscosity limit for the 1-D CNSEs under asymptotically isothermal assumptions.
More recently, Chen-Huang-Li-Wang-Wang \cite{Chen2023} proved the vanishing physical viscosity
limit for CNSPEs with general pressure law with large spherically symmetric initial data,
in which an $L^p$ compensated compactness framework for general pressure law was also established.
Finally, the existence of solutions for a general class of density-degenarate viscous models has been 
studied in \cite{CDS20,CDNP20}.

The main challenges tackled in this paper in comparison to the previous results
are due to the presence of the nonlocal terms. The existence of global-in-time solutions to the nonlocal Euler system \eqref{1.1} on the whole line has never been proven before. In contrast, it is known that the classical solutions of \eqref{1.1} without pressure and alignment, {\it i.e.,} $P(\rho)=\varpi=0$, may blow up in finite time \cite{Carrillo2016} and that they can be approximated by a degenerate vanishing viscosity method \cite{Carrillo2020}. More recently, Carrillo-Galtung \cite{CG23} showed that for such pressureless Euler systems in 1-D, the Lagrangian and entropy solutions are equivalent, which was then used to explore the long-time-asymptotics of the solutions. Moreover, Chaudhuri {\emph{et al.}} \cite{CCTZ23} studied the two-velocity reformulation of system \eqref{1.8} and derived an energy-type inequality, in the spirit of the Bresch-Desjardins estimate \cite{BD-2003-CRMASP,BD-2003-CPDE}. It was then used to construct the weak solutions and to study their long-time behaviour leading to the same density profiles as in \cite{Carrillo2016} or \cite{CG23}.
The  weak solutions to CNSEs with nonlocal attraction-repulsion forces in three dimensions were recently constructed in \cite{MSZ24}, see also \cite{Carrillo2019} for the constant viscosity case, where the long-time behavior of the solutions was considered.
Finally, global-in-time well-posedness theory for the Euler equations with  Riesz interactions and linear damping was recently established in \cite{Choi-Jung-Lee23} for the initial data near the equilibrium state and on the torus. We construct global entropy solutions to \eqref{1.1} via global weak solutions to the nonlocal Navier-Stokes system \eqref{1.8} with stress-free boundary conditions.
The proof that these solutions are uniformly bounded, in terms of $\varepsilon$,
despite the presence of nonlocal terms is another novelty. A new procedure for approximating the initial data, with the bounded second moment of the density, is implemented to obtain the propagation in time of the second moment of the density so that the confinement is controlled; see Lemma \ref{lemma4.2}.

\medskip

The paper has the following structure. In \S \ref{Def}, we first introduce the definitions of finite-energy entropy solutions for CEEs with nonlocal interactions
and present the main result -- Theorem \ref{thm:merged}. We then describe the construction of the approximate solutions for CNSEs and state the inviscid limit theorem, Theorem \ref{thm3.3}. The rest of the paper is then devoted to the proof of this result for two types of equations of state: the polytropic equation of state \eqref{pressure1}, and the general pressure satisfying hypotheses ($\mathcal{H}$).
In \S \ref{section3}, we derive  the basic estimates for both the density and velocity, and prove the $H^{-1}_{\rm loc}$--compactness of entropy dissipation measures  for the approximate solutions for the polytropic case.
In \S  \ref{section4}, we prove the global existence of finite-energy entropy weak solutions for this case, $\it{i.e.},$  we prove Theorem \ref{thm3.3}. 
In \S \ref{section5}, we develop the arguments from \S \ref{section3} and \S \ref{section4} to prove the global existence of finite-energy entropy solutions for CEEs for the general pressure law case. In this case, we prove the inviscid limit by proving the $W^{-1,p}_{\rm loc}$ compactness of entropy dissipation measures for the approximate solutions for $1\leq p<2$. Finally, in Appendix \ref{SecA}, we explain how the approximate initial data sequences $(\rho^\v_0(x),\rho^\v_0u^{\v}_0(x))$ can be constructed with desired estimates, regularity, and boundary compatibility for the polytropic equation of state.

\section{Definitions and Main Theorems}\label{Def}

In this section, we define the notion of entropy solutions of system \eqref{1.1}
and then formulate our main results of this paper.

\subsection{Entropy and entropy flux pairs for CEEs}
Recall that a pair 
$(\eta(\rho,m),q(\rho,m))$
is called an entropy-entropy flux pair (entropy pair, for short) of the Euler system \eqref{1.1}
if they satisfy
\begin{equation}\label{entropypair}
\nabla q(\rho,m)=\nabla \eta(\rho,m)\nabla \Big(\begin{matrix}m\\
	\frac{m^2}{\rho}+P(\rho)\end{matrix}\Big);
\end{equation}
see Lax \cite{P. D. Lax}.
Furthermore, $\eta(\rho,m)$ is called a weak entropy if
\begin{align}\nonumber
\eta|_{\rho=0}=0\qquad\,\, \text{for any fixed $u=\frac{m}{\rho}$}.
\end{align}
From \cite{Lions P.-L.2},
it is well known ({\it cf}. \cite{Perepelitsa,Chen6,Lions P.-L.2})
that, for the polytropic case, any weak entropy $(\eta,q)$ can be represented by
\begin{align}\label{weakentropy}
\begin{cases}
\displaystyle \eta^\psi(\rho,m)=\eta(\rho,\rho u)=\int_{\R}\chi(\rho;s-u)\psi(s)\,\dd s,\\[3mm]
\displaystyle q^\psi(\rho,m)=q(\rho,\rho u)=\int_{\R}(\theta s+(1-\theta)u)\chi(\rho;s-u)\psi(s)\, \dd s,
\end{cases}
\end{align}
where $\chi(\rho;s-u)=[\rho^{2\theta}-(s-u)^2]_{+}^{\fb}$
with $\fb=\frac{3-\gamma}{2(\gamma-1)}>-\frac{1}{2}$
and $\theta=\frac{\gamma-1}{2}$ is the weak entropy kernel.
In particular, when $\psi(s)=\frac{1}{2}s^2,$
the entropy pair is the pair of the mechanical energy and the associated flux{\rm :}
\begin{align*}\label{m-entropy}
\eta^{*}(\rho,m)=\frac{m^2}{2\rho}+\rho e(\rho),\quad q^{*}(\rho,m)=\frac{m^3}{2\rho^2}+m (\rho e(\rho))',
\end{align*}
From \eqref{entropypair}, any entropy satisfies
\begin{equation}\label{2.7}
	\eta_{\rho\rho}-\frac{P'(\rho)}{\rho^{2}}\eta_{uu}=0
\end{equation}
with $u=\frac{m}{\rho}$. It has been proved in \cite{Chen-LeFloch-2000, Chen-LeFloch-2003, Lions P.-L.1,Lions P.-L.2}
that any regular weak entropy can be generated by the convolution of a smooth function $\psi(x)$
with a fundamental solution $\chi(\rho,u,s)$ of the entropy equation \eqref{2.7}, {\it i.e.},
\begin{equation}\label{2.8}
	\eta^{\psi}(\rho,u)=\int_{\R}\chi(\rho,u,s)\psi(s)\,\mathrm{d}s.
\end{equation}
The corresponding entropy flux is generated from the flux kernel $\sigma(\rho,u,s)$ (see \eqref{2.10}), {\it i.e.},
\begin{align}\label{2.9}
	q^{\psi}(\rho,u)=\int_{\R}\sigma(\rho,u,s)\psi(s)\,\mathrm{d}s.
\end{align}
The entropy kernel $\chi=\chi(\r,u,s)$ is a fundamental solution of the entropy equation \eqref{2.7}{\rm :}
\begin{equation}\label{6.1}
	\left\{\begin{aligned}
		\dis&\chi_{\r\r}-\frac{P'(\r)}{\r^2}\chi_{uu}=0,\\
		\dis&\chi\vert_{\r=0}=0,\quad \chi_{\r}\vert_{\r=0}=\delta_{u=s}.
	\end{aligned}
	\right.
\end{equation}
As pointed out in \cite{Chen-LeFloch-2000} that equation \eqref{6.1}
is invariant under the Galilean transformation, which implies that
$\chi(\rho,u,s)=\chi(\r,u-s,0)=\chi(\r,0,s-u)$.
For simplicity, we write it as  $\chi(\r,u,s)=\chi(\r,u-s)$ below
when no confusion arises.
\smallskip
The corresponding entropy flux kernel $\sigma(\r,u,s)$ satisfies the Cauchy problem
for $\sigma-u\chi${\rm :}
\begin{equation}\label{2.10}
	\left\{\begin{aligned}
		\dis&(\sigma-u\chi)_{\r\r}-\frac{P'(\r)}{\r^2}(\sigma-u\chi)_{uu}=\frac{P''(\r)}{\r}\chi_{u},\\
		\dis&(\sigma-u\chi)\vert_{\r=0}=0,\quad  (\sigma-u\chi)_{\r}\vert_{\r=0}=0.
	\end{aligned}
	\right.
\end{equation}
It can be checked that $P(\rho)$, described by \eqref{pressure2}--\eqref{pressure4},
satisfies all the conditions in \cite{Chen-LeFloch-2000,Chen-LeFloch-2003}. In particular, $\sigma-u\chi$
is  Galilean invariant; see  \cite{Chen-LeFloch-2000}.

\subsection{Existence of solutions to  CEEs}
We are now ready to define the notion of solutions to system \eqref{1.1} and to formulate our main result.
\begin{definition}\label{definition-Euler}
A pair of functions $(\rho, m)(t,x)$ with $\rho \in L^\infty(\mathbb{R}_+; L^1_+(\mathbb{R}))$
and $\frac{m}{\sqrt{\rho}}
\in L^\infty(\mathbb{R}_+; L^2(\mathbb{R}))$ for $\mathbb{R}_+:=(0,\infty)$
is a finite-energy entropy solution of the Cauchy problem \eqref{1.1} and \eqref{1.6}, with $P(\rho)$ satisfying hypotheses {\rm (}$\mathcal{H}${\rm )}, if the following conditions hold{\rm :}
\begin{enumerate}
\item[\rm (i)]
The total mass is conserved{\rm :}
$$
\int_{\mathbb{R}} \rho(t,x)\,\dd x=\int_{\mathbb{R}} \rho_0(x)\,\dd x=:M \qquad\mbox{ a.e. $t\geq 0$},
$$
and $(\frac{m}{\sqrt{\rho}})(t,x)=0$ {\it a.e.} on
the vacuum states $\{(t,x)\,:\, \rho(t,x)=0\}$.

\smallskip
	
\item[\rm (ii)] For a.e. $t>0$, the total energy is not increasing{\rm :}
\begin{align*} 
\int_{\mathbb{R}}\Big(\frac{1}{2} \Big|\frac{m}{\sqrt{\rho}}\Big|^2
+\rho e(\rho)+\frac{1}{2} \rho W\ast\r\Big)(t,x)\, \dd x
\leq \mathcal{E}_0.
\end{align*}
		
\smallskip
\item[\rm (iii)] For any $\Psi(t,x)\in C^1_0([0,\infty)\times\mathbb{R})$,
\begin{align}
&\int_{\mathbb{R}^{2}_+} \big(\rho \Psi_t + m\Psi_x\big)\,\dd x\dd t
+\int_{\mathbb{R}} \rho_0(x) \Psi(0,x)\,\dd x=0,\label{2.5}\\[1mm]
&\int_{\mathbb{R}^{2}_+}\Big(m\Psi_t +\big(\frac{m^2}{\r}
+P(\rho)\big)\Psi_x\Big)\,\dd x\dd t
+\int_{\mathbb{R}} m_0(x)\Psi(0,x)\,\dd x\nonumber\\[-1mm]
&\quad=\int_{{\mathbb{R}^2_+}}\Big(-\lambda m+\r \partial_xW\ast\r-\int_{\R}\varpi(x-y)\big(\rho(x) m(y)-\rho(y) m(x)\big)\,\dd y\Big) \Psi \,\dd x\dd t.\label{2.6}
\end{align}

\item[\rm (iv)]
For any convex function $\psi(s)$ with subquadratic growth at infinity and any entropy pair $(\eta^{\psi},q^{\psi})$ defined in \eqref{2.8}--\eqref{2.9},
\begin{align*}
&\eta^{\psi}(\rho,m)_t+q^{\psi}(\rho,m)_x\nonumber\\
&-\eta^{\psi}_m\Big(\lambda m+\rho \int_{\R}\varpi(x-y)\big(\rho(x) m(y)-\rho(y) m(x)\big)\,\dd y-\rho \partial_x W\ast\rho \Big)\leq 0
\end{align*}
is satisfied in the sense of distributions.
\end{enumerate}
\end{definition}

\smallskip
We now state the main result of this paper.
\begin{theorem}[Existence of solutions of CEEs with nonlocal interactions]\label{thm:merged}
Consider problem \eqref{1.1}  with initial data \eqref{1.6} satisfying \eqref{1.7}--\eqref{2.1}. Let $P(\rho)$ satisfy hypotheses {\rm (}$\mathcal{H}${\rm )}.

Then there exists a global-in-time finite-energy entropy solution $(\rho, m )(t,x)$ of problem \eqref{1.1} and \eqref{1.6},
in the sense of   {\rm Definition \ref{definition-Euler}}.
In particular, there exists a global-in-time finite-energy entropy solution for
problem \eqref{1.1} and \eqref{1.6} with polytropic gases  equation of state \eqref{pressure1}.
\end{theorem}

\begin{remark}
The interaction potential can also be extended to the more general form:
$W(x)=-|x|+\frac{x^2}{2}+\tilde W (x),$ with
$|\tilde W (x)|\leq C x^2$ at infinity and smooth.
Even more, we can consider
$W(x)=-|x|+\frac{x^2}{2}+\frac{|x|^{\nu}}{\nu}$ for $1<\nu<2.$
\end{remark}

\begin{remark}
We can allow the communication kernel $\varpi(x)\in C(\R\backslash\{0\})$, and
$$\varpi(x)\mathbf{I}_{B(0,R)}\in L^{\frac{\gamma}{\gamma-1}}(\R), \quad \varpi(x)\mathbf{I}_{\R\backslash B(0,R)}\in L^{\infty}(\R) \qquad \text{for any } R>0.$$
For example, we can choose the singular communication weight $\varpi(x)$ as $\varpi(x)=\frac{1}{|x|^{a}}$ for $a<\frac{\gamma-1}{\gamma}.$
\end{remark}

\begin{remark} Our results also hold for an asymptotically isothermal pressure law $(\gamma_2=1)$, that is, $P(\rho)/\rho= O(1)$ in the limit $\rho\rightarrow \infty$; see {\rm \cite{Schrecker-Schulz-2019,Schrecker-Schulz-2020}}.
\end{remark}

\subsection{Existence of solutions to CNSEs and their inviscid limit}
In this paper, we do not prove Theorem \ref{thm:merged} directly. Instead, the solution from Definition \ref{definition-Euler} is obtained as a limit of the regular solutions of CNSEs \eqref{1.8}
 on a truncated domain $\Omega_T^\varepsilon$ with the moving boundary. We define
\begin{equation}\label{3.1}
\Omega_T^\varepsilon=\big\{(t,x)\ :\,  b^-_\varepsilon(t)\leq r\leq b^+_\varepsilon(t),\,0\leq t\leq T \big\},
\end{equation}
where $\{x=b^{\pm}_\varepsilon(t)\,:\,0<t\leq T\}$ are the free boundaries determined by
\begin{equation}\label{3.2}
\displaystyle \frac{\dd }{\dd t}b^\pm_\varepsilon(t)=u^\varepsilon(t,b^\pm_\varepsilon(t)), \quad b^\pm_\varepsilon(0)=\pm b_\varepsilon  \qquad\quad \mbox{for $t>0$}.
\end{equation}
On the free boundaries $x=b^\pm_\varepsilon(t)$,
the stress-free boundary condition is imposed{\rm :}
\begin{equation}\label{3.3}
\big(P(\rho^{\v})-\v\mu(\rho^{\varepsilon}) u^{\varepsilon}_x\big)(t,b^\pm_\varepsilon(t))
=0  \qquad \mbox{for $t>0$}.
\end{equation}
The initial intervals $[-b_\varepsilon, b_\varepsilon]$ approximate the whole line when $\v\to0^+$. More specifically, we  assume explicitly that
 $$
 b_\varepsilon=\varepsilon^{-p}\quad  \text{with}\quad  p>\left\{\begin{array}{ll}\frac{\gamma}{\gamma-\alpha} & \quad\text{for polytropic gases,}\\
 \frac{\gamma_1}{\gamma_1-\alpha}& \quad\text{for general pressure law.}
 \end{array}\right.
$$ 

Let $(\rho_0, m_0)(x)$ satisfy \eqref{1.7}--\eqref{2.1}.
The initial data for the approximate system are given by
\begin{equation}\label{3.4}
(\rho^{\varepsilon},\rho^{\varepsilon} u^{\varepsilon})(0,x)
=(\r_0^{\v},\rho_0^{\v}u_0^{\v})(x)\qquad\,\mbox{for $x\in [-b_\varepsilon,b_\varepsilon]$},
\end{equation}
where $(\rho_0^{\v}, u_0^{\v})(x)$ are obtained by the smooth approximation of
the original data $(\rho_0,m_0)$, constructed in Appendix A. In particular,  they satisfy Lemmas \ref{lemA.2}--\ref{propA.1}.

We further define
\begin{align}\label{3.5}
\mathcal{E}_{0}^{\v}:&=
\int_{-b_\varepsilon}^{b_\varepsilon}\Big(\frac{1}{2}\rho_0^{\v}|u_0^{\v}|^2+ \rho_0^{\v}e(\rho_0^{\v})+\frac{1}{2}\rho_0^{\v}W\ast\r^{\varepsilon}_0\Big)\,
\dd x,\,\,\, &\\[1.5mm]
\mathcal{E}_1^{\v}
 :&= \v^2 \int_{-b_\varepsilon}^{b_\varepsilon} (\rho^{\varepsilon}_0)^{2\alpha-3}|\rho^{\varepsilon}_{0x}|^2\,\dd x,\label{3.6}
\end{align}
where the convolution in \eqref{3.5} is understood as $\int^{b_\varepsilon}_{-b_\varepsilon}W(x-y)\rho_0^{\varepsilon}(y)\,\dd y$.
\smallskip

We  are now ready to state the main theorem of our paper.

\begin{theorem}[Inviscid limit for CNSEs with nonlocal interactions]\label{thm3.3}
Let the hypotheses of  \mbox{ {\rm Theorem \ref{thm:merged}}} be satisfied.
Let $\gamma_1\geq\gamma_2>1$ and $\frac{2}{3}<\alpha\leq1$,  and let
$\{(\rho^{\v},\rho^{\v}u^{\v})\}_{\v\in (0,1]}$ be a sequence of the strong solutions to \eqref{1.8} with initial and boundary data  specified in \eqref{3.1}--\eqref{3.4} and satisfying:
\begin{align}
& 0< C^{-1}_{\varepsilon}\leq \r^{\varepsilon}_0(x)\leq C_{\varepsilon},\quad
(\r^{\varepsilon}_0(b))^{\gamma_1} b_\v= (\r^{\varepsilon}_0(-b_\varepsilon))^{\gamma_1} b_\varepsilon \leq C_0
\qquad\mbox{ for $\varepsilon \in (0,1]$},\nonumber\\[2mm]
&(\rho^\v_0,\rho_0^\v u^{\v}_0)(x) \rightarrow (\rho_0,m_0)(x)
 \qquad \mbox{as $\v\to 0^+$ in $L^{q}_{\rm loc}(\R)\times L^1_{\rm loc}(\R)$} \mbox{ for $q\in\{1,\gamma_2\}$},\nonumber\\[2mm]
&(\mathcal{E}_0^\v,\mathcal{ E}_1^\v)\rightarrow (\mathcal{E}_0,0) \qquad \mbox{as $\v\to 0^+$}, \nonumber\\[2mm]
&\int_{-\infty}^{\infty}\rho^\v_0(x)\,\dd x= M,
\qquad \mathcal{E}^{\v}_0+\v^{-1}\mathcal{E}_1^\v\le C_0,
\quad \nonumber \\[2mm]
&\int^{\infty}_{-\infty}x^2\rho^{\v}_0(x)\,\dd x\,\rightarrow \,\int^{\infty}_{-\infty}x^2\rho_0 (x)\,\dd x=M_2\qquad\mbox{as $\v\to 0^+$},\nonumber
\end{align}
where $C_0>0$ is a constant independent of $\varepsilon\in (0,1]$ which may depend on $(\mathcal{E}_0, M, M_2, \gamma,\alpha)$, while
$C_\v>0$ is a constant depending on $\v>0$.

Then there exist both a subsequence {\rm (}still denoted{\rm )} $(\rho^{\v},m^{\v} )(t,x)$ with $m^{\v}=\rho^\v u^\v$, and a vector-valued function $(\rho, m)(t,x)$
such that, as $\v\rightarrow0^+$,
\begin{align}\nonumber
\begin{split}
&(\rho^{\v},m^{\v})(t,x)\rightarrow (\rho,m)(t,x)\qquad
  \mbox{a.e. $(t,x)\in \R_+\times\R$},\\
&(\rho^{\v},m^{\v})(t,x)\rightarrow (\rho,m)(t,x)\qquad
  \mbox{in $L^{q_1}_{\rm loc}(\R^2_+)\times L^{q_2}_{\rm loc}(\R^2_+)$ with $q_1\in[1,\gamma_2+1)$ and $q_2\in [1,\frac{3(\gamma_2+1)}{\gamma_2+3})$}.
\end{split}
\end{align}
Moreover, $(\rho,m)(t,x)$
is a global finite-energy entropy
solution of the Cauchy problem \eqref{1.1} and \eqref{1.6}
in the sense of  {\rm Definition \ref{definition-Euler}}.
In particular, the result is true for the polytropic case \eqref{pressure1}.
\end{theorem}

\begin{remark}
In   {\rm Theorem \ref{thm3.3}}, we need $\frac{2}{3}<\alpha\leq 1$; see \eqref{A.4222}. However, note that  the condition $\frac{2}{3}<\alpha<\gamma$ is needed only for obtaining the uniform estimates
for the approximate solutions{\rm ;} see {\rm Lemmas \ref{lem4.1}--\ref{lem4.2}} and
 {\rm Lemmas \ref{lem5.2}--\ref{lem5.9}}.
In particular,  our results cover the Saint-Venant model for shallow water ($\alpha=1$ and $\gamma=2$),
even with the nonlocal terms.
\end{remark}

\begin{remark}   The regular solutions for the compressible Navier-Stokes problem \eqref{1.8} with the new nonlocal terms  and with stress-free boundary conditions \eqref{3.2}--\eqref{3.3} can be obtained by following the arguments similar to {\rm \cite{Jiang,Qin}}; see also {\rm\cite{CCTZ23}} for the Dirichlet boundary conditions. 
\end{remark}

\section{Uniform Estimates of Approximate Solutions}\label{section3}
This section is dedicated to several uniform estimates with respect to $\v\in (0, 1]$
for the regular solutions of \eqref{1.8} on the bounded time-dependent domain \eqref{3.1} with corresponding boundary conditions. Although some of the estimates hold also for general $P(\rho)$ specified through hypotheses $(\mathcal{H})$, for the purposes of this section, we restrict ourselves to the polytropic equation of state \eqref{pressure1}. The general pressure case is discussed in \S\ref{section5}.

We drop the $\v$-superscript of approximate solutions when no confusion arises.
In all of the estimates from now on,
$C>0$ is a universal constant independent of $\v$, which may depend on $\mathcal{E}_0, M, M_2, C_0,T,\gamma$, and $\alpha$, and can be different at each occurrence. For abbreviation, we also denote
\begin{equation}\label{DefV}
V=V(t,x):=\int_{\R}\varpi(x-y)\Big(m(y)-\frac{m(x)}{\rho(x)}\rho(y)\Big)\,\dd y
=\int^{b^+(t)}_{b^-(t)}\varpi(x-y)\Big(m(y)-\frac{m(x)}{\rho(x)}\rho(y)\Big)\,\dd y
\end{equation}
for $x\in (b^-(t),b^+(t))$.

\subsection{Conservation of mass and Lagrangian coordinates}
We present the reformulation of IBVP \eqref{1.8} and \eqref{3.1}--\eqref{3.4}
in the so-called Lagrangian coordinates,
which is equivalent to the Eulerian formulation for sufficiently regular solutions.
It follows from $\eqref{1.8}_1$ and \eqref{3.2} that
\begin{equation*}
\frac{\dd}{\dd t}\int_{b^-(t)}^{b^+(t)}\rho(t,x)\,\dd x= (\rho u)(t,b^+(t))-(\rho u)(t,b^-(t))-\int_{b^-(t)}^{b^+(t)}(\rho u )_x(t,x)\,\dd x=0,
\end{equation*}
which, due to \eqref{A.23-1}, yields that
\begin{equation}\label{3.7}
\int_{b^-(t)}^{b^+(t)}\rho(t,x)\,\dd x=\int_{b^-(t)}^{b^+(t)}\rho^{\v}_0(x)\,\dd x= M \qquad\,\, \mbox{for any $t\geq0$}.
\end{equation}
For $x\in[b^-(t),b^+(t)]$ and $t\in[0,T]$, we define the Lagrangian
coordinates $(\tau,\xi)$ as
\begin{equation*}
\xi=\int_{b^-(t)}^x \rho(t,y)\,\dd y,\qquad  \tau=t,
\end{equation*}
which transform the domain with the free boundary: $[0,T]\times[b^-(t),b^+(t)]$
into the fixed domain: $[0,T]\times[0,M]$.
A direct calculation shows that
\begin{align*}\label{2.9}
\begin{cases}
\displaystyle \frac{\partial \xi}{\partial t}=-\rho u,\quad \frac{\partial \xi}{\partial x}=\rho,\quad\frac{\partial\tau}{\partial t}=1, \quad\frac{\partial\tau}{\partial x}=0,\\[3mm]
\displaystyle \frac{\partial x}{\partial \tau}=u,\quad \frac{\partial x}{\partial\xi}=\frac{1}{\rho},\quad \frac{\partial t}{\partial\tau}=1,\quad \frac{\partial t}{\partial \xi}=0.
\end{cases}
\end{align*}

Applying the Euler-Lagrange transformation, IBVP \eqref{1.8} and \eqref{3.1}--\eqref{3.4}
becomes
\begin{equation}\label{3.8}
\begin{cases}
\displaystyle\rho_\tau+\rho^2u_{\xi}=0,\\[2mm]
\displaystyle u_\tau+ P(\rho)_\xi-\varepsilon (\mu(\rho)\rho u_{\xi})_{\xi}-\lambda u-V+\rho(W\ast\rho)_{\xi}=0,
\end{cases}
\end{equation}
for $(\tau,\xi)\in[0,T]\times[0,M]$, and
\begin{equation}\label{3.9}
 \big(P(\rho)-\v \mu(\rho)\rho u_{\xi}\big)(\tau,0)=0, \quad\big(P(\rho)-\v \mu(\rho)\rho u_{\xi}\big)(\tau,M)=0
 \qquad\,\,  \mbox{ for $\tau\in[0,T]$}.
\end{equation}
With a slight abuse of notation, we denote by $V$ a function defined in \eqref{DefV}, taken as a function of $\xi$.
The fixed boundaries $\xi=0,M$ correspond to the free boundaries $x=b^{\pm}(t)$ in the Eulerian coordinates.

\subsection{Basic energy estimate}
To begin with, we obtain the basic energy estimate for problem \eqref{1.8}
and \eqref{3.1}--\eqref{3.4}.

\begin{lemma}\label{lem4.1}
For a smooth solution of problem \eqref{1.8} and \eqref{3.1}--\eqref{3.4}, the following estimate holds{\rm :}
\begin{align}
&\int^{b^+(t)}_{b^-(t)}\Big(\frac{1}{2}\rho u^2+\rho e(\rho)
+\frac{1}{2}\rho W\ast\rho\Big)(t,x)\,\dd x
+\int_0^t\int_{b^-(s)}^{b^+(s)}
  \lr{\v\rho^{\alpha}|u_x|^2-\lambda \rho u^2}(s,x)\,\dd x\dd  s\nonumber\\
&+\frac{1}{2}\int^t_0\int^{b^+(s)}_{b^-(s)}\int^{b^+(s)}_{b^-(s)}\varpi(x-y)\rho(x)\rho(y)|u(y)-u(x)|^2\,\dd y\dd x\dd s
= \mathcal{E}^{\varepsilon}_0\le C_0.\label{4.1}
\end{align}
Here and hereafter, $C_0>0$ is the constant from the statement of {\rm Theorem \ref{thm3.3}}, independent of $\varepsilon\in (0,1]$.
\end{lemma}

\noindent{\bf Proof.}
It follows from $\eqref{3.8}_1$ that
\begin{align}\label{4.2}
-P(\rho)u_{\xi}=\frac{P(\rho)}{\rho^2}\rho_{\tau}=\kappa\rho^{\gamma-2}\rho_{\tau}=e(\rho)_{\tau}.
\end{align}
Multiplying $\eqref{3.8}_2$ by $u,$ we have
\begin{align}\label{4.3}
\Big(\frac{u^2}{2}+e(\rho)\Big)_{\tau}+\varepsilon\rho^{1+\alpha}(u_{\xi})^2
+ \big((P(\rho)-\varepsilon\mu(\rho)\rho u_{\xi})u\big)_{\xi}
-\lambda u^2-u V+\rho u (W\ast\rho)_{\xi}=0.
\end{align}
Notice that
\begin{align}
-\int^M_0u V\,\dd \xi
&=-\int^{b^+(t)}_{b^-(t)}\rho(x)u(x)\Big(\int^{b^+(t)}_{b^-(t)}\varpi(x-y)
  \big(u(y)-u(x)\big)\rho(y)\,\dd y\Big)\dd x\nonumber\\
&=\int^{b^+(t)}_{b^-(t)}\int^{b^+(t)}_{b^-(t)}\varpi(x-y)\rho(x)
 \rho(y)\big(|u(x)|^2-u(x)u(y)\big)\,\dd y\dd x\nonumber\\
&=\frac{1}{2}\int^{b^+(t)}_{b^-(t)}\int^{b^+(t)}_{b^-(t)}\varpi(x-y)\rho(x)\rho(y)|u(y)-u(x)|^2\,\dd y\dd x,\nonumber
\end{align}
where we have used the symmetry of $\varpi(x)$.

For the other nonlocal term, our aim is to prove that
\begin{equation}\label{result}
\int^M_0\rho u(W \ast \rho)_{\xi}\,\dd \xi=\frac{\dd}{\dd t}\int^{b^+(t)}_{b^-(t)}\frac{1}{2}\rho W
\ast \rho \,\dd x.
\end{equation}
First we rewrite the right-hand-side (RHS) to obtain:
\begin{align}\label{4.8a}
&\frac{\dd}{\dd t}\int^{b^+(t)}_{b^-(t)}\frac{1}{2}\rho W\ast\rho \,\dd x\nonumber\\
&\,\,=\int^{b^+(t)}_{b^-(t)}\Big(\frac{1}{2}\rho W\ast\rho\Big)_t\,\dd x
+\frac{1}{2}(\rho W\ast\rho)(t,b^+(t))\frac{\dd b^+(t)}{\dd t}
-\frac{1}{2}(\rho W\ast\rho)(t,b^-(t))\frac{\dd b^-(t)}{\dd t}\nonumber\\
&\,\,=\frac{1}{2}\int^{b^+(t)}_{b^-(t)}\lr{\rho_t W\ast\rho+\rho W\ast\rho_t}\,\dd x
+\frac{1}{2}(u\rho W\ast\rho)(t,b^+(t))
-\frac{1}{2}(u\rho W\ast\rho)(t,b^-(t)).
\end{align}
Using the change of the order of integration
and the explicit form of $W(x)=-|x|+\frac{|x|^2}{2}$  (in particular $W(x-y)=W(y-x)$), we obtain
\begin{align*}
&\int^{b^+(t)}_{b^-(t)}\rho W\ast\rho_t\,\dd x\nonumber\\
&=\int^{b^+(t)}_{b^-(t)}\rho(t,x)\big(\int^{b^+(t)}_{b^-(t)}W(x-y)\rho(t,y)\,\dd y\big)_t\,\dd x\nonumber\\
&=\int^{b^+(t)}_{b^-(t)}\rho_t(t,y)\big(\int^{b^+(t)}_{b^-(t)}W(x-y)\rho(t,x)\,\dd x\big)\,\dd y
  +(\rho u W\ast \rho)(t,b^+(t))-(\rho u W\ast \rho)(t,b^-(t))\nonumber\\
&=\int^{b^+(t)}_{b^-(t)}\rho_t(t,y)\big(\int^{b^+(t)}_{b^-(t)}W(y-x)\rho(t,x)\,\dd x\big)\,\dd y
  +(\rho u W\ast \rho)(t,b^+(t))-(\rho u W\ast \rho)(t,b^-(t))\nonumber\\
&=\int^{b^+(t)}_{b^-(t)}\rho_t(t,x)\big(\int^{b^+(t)}_{b^-(t)}W(x-y)\rho(t,y)\,\dd y\big)\,\dd x
  +(\rho u W\ast \rho)(t,b^+(t))-(\rho u W\ast \rho)(t,b^-(t))\nonumber\\
&=\int^{b^+(t)}_{b^-(t)}\rho_t W\ast\rho \,\dd x+(\rho u W\ast \rho)(t,b^+(t))-(\rho u W\ast \rho)(t,b^-(t)).
\end{align*}
Thus, it follows that the first term on the RHS of \eqref{4.8a} is equal to
$$
\frac{1}{2}\int^{b^+(t)}_{b^-(t)}\big(\rho_t W\ast\rho+\rho W\ast\rho_t\big)\,\dd x
=\int^{b^+(t)}_{b^-(t)}\rho_t W\ast\rho \,\dd x+\frac{1}{2}(\rho u W\ast \rho)(t,b^+(t))-\frac{1}{2}(\rho u W\ast \rho)(t,b^-(t)).
$$
Therefore, we have proven that
\begin{align}\label{new3.10}
&\frac{\dd}{\dd t}\int^{b^+(t)}_{b^-(t)}\frac{1}{2}\rho W\ast\rho \,\dd x=\int^{b^+(t)}_{b^-(t)}\rho_t W\ast\rho\,\dd x
+(u\rho W\ast\rho)(t,b^+(t))
-(u\rho W\ast\rho)(t,b^-(t)).
\end{align}
On the other hand, it follows  directly from $\eqref{1.8}_1$ that
\begin{align}
\int^M_0\rho u (W\ast\rho)_{\xi}\,\dd \xi
&=\int^{b^+(t)}_{b^-(t)}\rho u (W\ast\rho)_x\,\dd x\nonumber\\
&=(\rho u W\ast\rho)(t,b^+(t))-(\rho u W\ast\rho)(t,b^-(t))-\int^{b^+(t)}_{b^-(t)}(\rho u)_x W\ast\rho\,\dd x\nonumber\\
&=(\rho u W\ast\rho)(t,b^+(t))-(\rho u W\ast\rho)(t,b^-(t))+\int^{b^+(t)}_{b^-(t)}\rho_t W\ast\rho\,\dd x.\label{4.4}
\end{align}
Comparing \eqref{new3.10} with \eqref{4.4}, we deduce \eqref{result}.

Since $W=-|x|+\frac{|x|^2}{2},$ 
we have that $W+\frac{1}{2}\geq0.$ and so, using \eqref{3.7}, we obtain
$$
\frac{\dd}{\dd t}\int^{b^+(t)}_{b^-(t)}\frac{1}{2}\rho W\ast\rho \,\dd x=\frac{\dd}{\dd t}\int^{b^+(t)}_{b^-(t)}\frac{1}{2}\rho\Big(W+\frac{1}{2}\Big)\ast\rho \,\dd x.
$$
Integrating \eqref{4.3} over $[0,\tau]\times[0,M]$, using the Gr\"{o}nwall inequality
and the stress-free boundary conditions \eqref{3.9}, pulling the resultant equation
back to the Eulerian coordinates, we obtain \eqref{4.1}.
This completes the proof.
$\hfill\Box$

\smallskip
\subsection{Second moment estimate}
We now derive the second moment estimate for the density.
\begin{lemma}\label{lemma4.2}
There exists $C=C>0$, independent of $\v$, such that
\begin{align}\label{4.10}
\int^{b^+(t)}_{b^-(t)}x^2\rho(t,x) \,\dd x\leq C.
\end{align}
\end{lemma}
\noindent{\bf Proof.}
Using $\eqref{1.8}_1,$ we have
\begin{align*}
\frac{\dd}{\dd t}\int^{b^+(t)}_{b^-(t)}x^2\rho  \,\dd x
&=(b^+(t))^2(\rho u)(t,b^+(t))
-(b^- (t))^2(\rho u)(t,b^-(t))+\int^{b^+(t)}_{b^-(t)} x^2\rho_t\,\dd x\\
&=2\int^{b^+(t)}_{b^-(t)} x\rho u\,\dd x,
\end{align*}
such that 
$$
\frac{\dd}{\dd t}\int^{b^+(t)}_{b^-(t)}x^2\rho  \,\dd x
\leq\int^{b^+(t)}_{b^-(t)}x^2\rho \,\dd x+\int^{b^+(t)}_{b^-(t)}\rho u^2 \,\dd x.
$$
Then we conclude \eqref{4.10} by using the Gr\"{o}nwall inequality. 
$\hfill\Box$

\subsection{Higher-order estimates of the density and the pressure}
In this section, we derive several higher-order estimates
for the density and the pressure.
To start, we analyze the behavior of density $\rho$  on the free boundary.
It follows from $\eqref{3.8}_1$ and \eqref{3.9} that
\begin{align*}
\frac{\dd}{\dd\tau}(\rho^{\alpha-\gamma}(\tau,M))=-\frac{\kappa(\alpha-\gamma)}{\varepsilon}.
\end{align*}
Then we have
\begin{align}\nonumber
\rho(\tau,M)=\rho_0(M)\,\Big(1+\frac{\kappa(\gamma-\alpha)}{\varepsilon}(\rho_0(M))^{\gamma-\alpha}\tau\Big)^{-\frac{1}{\g-\alpha}},
\end{align}
so that
\begin{align}\nonumber
\rho(t,b^+(t))=\rho_0(b)\,\Big(1+\frac{\kappa(\gamma-\alpha)}{\varepsilon}(\rho_0(b))^{\gamma-\alpha}t\Big)^{-\frac{1}{\g-\alpha}}.
\end{align}
We obtain
\begin{align}\label{4.29}
\rho(t,b^+(t))=\rho(t,b^-(t))=\rho_0(b)\Big(1+\frac{\kappa(\g-\alpha)}{\v}(\r_0(b))^{\gamma-\alpha}t\Big)^{-\frac{1}{\g-\alpha}}\leq \rho_0(b),
\end{align}
where we have used that $\rho_0(b)=\rho_0(-b).$
From \eqref{4.29}, it follows in particular that the third and fourth term
on the left-hand side (LHS) of \eqref{4.30} below is nonnegative.

\begin{lemma} \label{lem2.2}
Let $\alpha\geq0$ and $\rho^{\gamma}_0(\pm b)b\leq C_0$.
Then, for any given $T>0$, there exists $C>0$ such that, for any $t\in[0,T]$,
\begin{align}\label{4.30}
&\varepsilon^2\int^{b^+(t)}_{b^-(t)}(\rho^{2\alpha-3}\rho^2_x)(t,x)\,\dd x
+\varepsilon\kappa\gamma\int^t_0\int^{b^+(t)}_{b^-(t)}(\rho^{\alpha+\gamma-3}\rho^2_x)(s,x)\,\dd x\dd s\nonumber\\
&+\frac{\kappa\gamma}{\varepsilon}\int^{t}_0
\Big(\rho^{2\gamma-\alpha}(s,b^+(s))b^+(s)-\rho^{2\gamma-\alpha}(s,b^-(s))b^-(s)\Big)\,\dd s\nonumber\\
&+\kappa\Big(\rho^{\gamma}(t,b^+(t))b^+(t)-\rho^{\gamma}(t,b^-(t))b^-(t)\Big)\leq C .
\end{align}
\end{lemma}

\noindent{\bf Proof.}
It is direct to calculate from $\eqref{3.8}_2$ that
\begin{align}\label{4.31}
\Big(u+\frac{\varepsilon}{\alpha}(\rho^{\alpha})_{\xi}\Big)_{\tau}+P(\rho)_{\xi}-\lambda u-V+\rho(W\ast\rho)_{\xi}=0,
\end{align}
where we have used the fact that
\begin{align}\nonumber
\mu(\rho)\rho u_{\xi}=\rho^{1+\alpha}u_{\xi}=-\rho^{\alpha-1}\rho_{\tau}=-\frac{1}{\alpha}(\rho^{\alpha})_{\tau}.
\end{align}
Multiplying \eqref{4.31} by $u+\frac{\varepsilon}{\alpha}(\rho^{\alpha})_{\xi},$ we obtain
\begin{align}\label{4.32}
&\Big(\frac{(u+\frac{\varepsilon}{\alpha}(\rho^{\alpha})_{\xi})^2}{2}
  +\frac{\kappa\rho^{\gamma-1}}{\gamma-1}+\frac{1}{2}\rho W\ast\rho\Big)_{\tau}
   +\big(p(\rho)u\big)_{\xi}+\varepsilon\kappa\gamma\rho^{\alpha+\gamma-2}(\rho_{\xi})^2\nonumber\\
&+\big(u+\frac{\varepsilon}{\alpha}(\rho^{\alpha})_{\xi}\big)\big(-\lambda u-V+\rho(W\ast\rho)_{\xi}\big)=0.
\end{align}
Using $\eqref{3.8}_1$ and \eqref{3.9}, we have
\begin{align}\nonumber
\rho_{\tau}(\tau,M)=-\frac{\kappa}{\varepsilon}\rho^{\gamma+1-\alpha}(\tau,M),
\qquad \rho_{\tau}(\tau,0)=-\frac{\kappa}{\varepsilon}\rho^{\gamma+1-\alpha}(\tau,0).
\end{align}
Then we obtain
\begin{align}
\big(P(\rho)u\big)(\tau,M)-\big(P(\rho)u\big)(\tau,0)
&=\kappa\rho^{\gamma}(\tau,M)\frac{\dd b^+(\tau)}{\dd \tau}-\kappa\rho^{\gamma}(\tau,0)\frac{\dd b^-(\tau)}{\dd \tau}\nonumber\\
&=\kappa\big(\rho^{\gamma}(\tau,M)b^+(\tau)\big)_{\tau}
  -\kappa\big(\rho^{\gamma}(\tau,0)b^-(\tau)\big)_{\tau}\nonumber\\
&\quad+\kappa\gamma\Big(\rho^{\gamma-1}(\tau,0)\rho_{\tau}(\tau,0)b^{-}(\tau)-\rho^{\gamma-1}(\tau,M)\rho_{\tau}(\tau,M)b^{+}(\tau)\Big)\nonumber\\
&=\kappa\big(\rho^{\gamma}(\tau,M)b^+(\tau)\big)_{\tau}-\kappa\big(\rho^{\gamma}(\tau,0)b^-(\tau)\big)_{\tau}\nonumber\\
&\quad+\frac{\kappa\gamma}{\varepsilon}\Big(\rho^{2\gamma-\alpha}(\tau,M)b^+(\tau)-\rho^{2\gamma-\alpha}(\tau,0)b^-(\tau)\Big).\nonumber
\end{align}

Integrating \eqref{4.32} over $[0,M]$ yields
\begin{align}\label{4.34}
&\frac{\dd}{\dd\tau}\int^{M}_0\Big(\frac{1}{2}\big(u+\frac{\varepsilon}{\alpha}(\rho^{\alpha})_{\xi}\big)^2
  +\frac{\kappa\rho^{\gamma-1}}{\gamma-1}\Big)\,\dd \xi
  +\varepsilon\kappa\gamma\int^M_0\rho^{\alpha+\gamma-2}(\rho_{\xi})^2\,\dd \xi \nonumber\\
&\,\,\, +\kappa\Big(\rho^{\gamma}(\tau,M)b^+(\tau)-\rho^{\gamma}(\tau,0)b^-(\tau)\Big)_{\tau}
 +\frac{\kappa\gamma}{\varepsilon}\Big(\rho^{2\gamma-\alpha}(\tau,M)b^+(\tau)-\rho^{2\gamma-\alpha}(\tau,0)b^-(\tau)\Big)\nonumber\\
&=\int^M_0\big(u+\frac{\varepsilon}{\alpha}(\rho^{\alpha})_{\xi}\big)\big(-\lambda u-V+\rho(W\ast\rho)_{\xi}\big)\,\dd \xi.
\end{align}
For the last term in \eqref{4.34}, we have
\begin{align}
&\int^M_0\big(u+\frac{\varepsilon}{\alpha}(\rho^{\alpha})_{\xi}\big)\big(-\lambda u-V+\rho(W\ast\rho)_{\xi}\big)\,\dd \xi\nonumber\\
&\leq\int^M_0\big(u+\frac{\varepsilon}{\alpha}(\rho^{\alpha})_{\xi}\big)^2\,\dd \xi
   +\int^M_0\big(-\lambda u-V+\rho(W\ast\rho)_{\xi}\big)^2\,\dd \xi\nonumber\\
&\leq \int^M_0\big(u+\frac{\varepsilon}{\alpha}(\rho^{\alpha})_{\xi}\big)^2\,\dd \xi
   +C\int^M_0u^2\,\dd \xi+C\int^M_0 V^2\,\dd \xi+ C\int^M_0\big(\rho(W\ast\rho)_{\xi}\big)^2\,\dd \xi\nonumber\\
&\leq \int^M_0\big(u+\frac{\varepsilon}{\alpha}(\rho^{\alpha})_{\xi}\big)^2\,\dd \xi
  + C\int^{b^+(t)}_{b^-(t)}\rho\big(\partial_xW\ast\rho\big)^2\,\dd x+C,\nonumber
\end{align}
where we have used the fact that
\begin{align}
&\int^M_0 V^2\,\dd \xi=\int^{b^+(t)}_{b^-(t)}\rho V^2\,\dd x\nonumber\\
&=\int^{b^+(t)}_{b^-(t)}\rho(t,x)\Big(\int^{b^+(t)}_{b^-(t)}\varpi(x-y)\big(u(t,y)-u(t,x)\big)\rho(t,y)\,\dd y\Big)^2\,\dd x\nonumber\\
&\leq C\int^{b^+(t)}_{b^-(t)}\bigg(\rho(t,x)\Big(\int^{b^+(t)}_{b^-(t)}\big(\rho(t,y)+(\rho u^2)(t,y)\big)\,\dd y\Big)^2
+u^2(t,x)\bigg)\,\dd x
\leq C.\nonumber
\end{align}
Notice that
\begin{align}\label{4.36}
(\partial_xW\ast\rho)(t,x)
=\int^{b^+(t)}_{b^-(t)}\big(1-2H(x-y)+x-y\big)\rho(t,y)\,\dd y,
\end{align}
where $H(\cdot)$ is the Heaviside function.
We now rewrite the RHS of  \eqref{4.36} using 
\begin{align*}
\int^{b^+(t)}_{b^-(t)}\rho(t,y)\,\dd y=M,\qquad \int^{b^+(t)}_{b^-(t)}H(x-y)\rho(t,y)\,\dd y=\int^x_{b^-(t)}\rho(t,y)\,\dd y.
\end{align*}
We obtain
\begin{align}\nonumber
(\partial_xW\ast\rho)(t,x)=M-2\int^x_{b^-(t)}\rho(t,y)\,\dd y+xM-\int^{b^+(t)}_{b^-(t)}y\rho (t,y)\,\dd y.
\end{align}
It follows from \eqref{3.7} and \eqref{4.10} that
$$\Big|\int^{b^+(t)}_{b^-(t)}|x|\rho(t,x) \,\dd x\Big|\leq\int^{b^+(t)}_{b^-(t)}(x^2+1)\rho(t,x) \,\dd x\leq C.$$
Then we have
\begin{align}\nonumber
\int^{b^+(t)}_{b^-(t)}\rho(t,x)(\partial_xW\ast\rho)^2(t,x)\,\dd x
\leq\int^{b^+(t)}_{b^-(t)}\rho\Big(3M+|x|M+\int^{b^+(t)}_{b^-(t)}|y|\rho(t,y)\,\dd y\Big)^2\,\dd x\leq C.
\end{align}
\smallskip
Integrating \eqref{4.32} over $[0,\tau]$ leads to
\begin{align}\nonumber
&\int^M_0\Big(\frac{1}{2}\big(u+\frac{\varepsilon}{\alpha}(\rho^{\alpha})_{\xi}\big)^2
  +\frac{\kappa}{\gamma-1}\rho^{\gamma-1}\Big)\,\dd \xi
  +\varepsilon\kappa\gamma\int^{\tau}_0\int^M_0\rho^{\alpha+\gamma-2}(\rho_{\xi})^2\,\dd \xi \dd s\nonumber\\
&\,\,\,+\kappa\Big(\rho^{\gamma}(\tau,M)b^+(\tau)-\rho^{\gamma}(\tau,0)b^-(\tau)\Big)+\frac{\kappa\gamma}{\varepsilon}\int^{\tau}_0\Big(\rho^{2\gamma-\alpha}(s,M)b^+(s)-\rho^{2\gamma-\alpha}(s,0)b^-(s)\Big)\,\dd s\nonumber\\
&\leq C\int^M_0\Big(\frac{1}{2}\big(u_0+\frac{\varepsilon}{\alpha}(\rho^{\alpha}_0)_{\xi}\big)^2
+\frac{\kappa}{\gamma-1}\rho^{\gamma-1}_0\Big)\,\dd \xi
+\kappa\big(\rho^{\gamma}_0(M)+\rho^{\gamma}_0(0)\big)b+C.\nonumber
\end{align}
This completes the proof.
$\hfill\Box$

\medskip

Our next aim is to show that the domain $\Omega_T$ expands to the whole physical space $[0,T]\times\R$, that is,
$\displaystyle\inf_{t\in[0,T]}b(t)\rightarrow\infty$ as $b\rightarrow\infty$. We have the following result.

\begin{lemma}[\cite{He}, Lemma 2.3]\label{lem4.4}
Choose $0\leq\alpha<\gamma,$ $p>\frac{\gamma}{\gamma-\alpha},$ and $b:=\varepsilon^{-p}.$
Then, for any given $T>0$, there exists $\varepsilon_0>0$ such that, for $\v\in(0,\v_0]$,
\begin{align*}
\pm b^\pm(t)\geq \frac12 b
\qquad \mbox{for $t\in[0,T]$}.
\end{align*}
\end{lemma}

In particular, the proof of Lemma \ref{lem4.4} is independent of the nonlocal terms. Similar result for the general pressure law is proven later; see Lemma \ref{lem5.6}.

\begin{lemma}[\bf Higher integrability of the density]\label{lem4.5}
Let $(\r,u)$ be the smooth solution of \eqref{1.8} and \eqref{3.1}--\eqref{3.4},
and let the assumption of {\rm Lemma \ref{lem4.4}} hold.
Then, for any $K\Subset[b^-(t),b^+(t)]$ for any $t \in[0,T]$, there exists
$C(K)>0$ independent of $\v\in (0,1]$ such that
\begin{align}\label{4.39}
\int_0^T\int_K\r^{\g+1}(t,x)\,\dd x \dd t\leq C(K).
\end{align}
\end{lemma}

\noindent{\bf Proof.}  We divide the proof into two steps.

\smallskip
\noindent {\emph{Step 1.}}  For given $K\Subset[b^-(t),b^+(t)]$ for any $t\in[0,T]$,
there exist $r_1$ and $r_2$ such that $K\Subset (r_1,r_2)\Subset[b^-(t),b^+(t)]$.
Let $w(x)$ be a smooth function with $\text{supp}\,w\subseteq(r_1,r_2)$ and $w(x)=1$ for $x\in K$.
Multiplying $\eqref{1.8}_2$ by $w(x)$, we have
\begin{align}\label{4.40}
&(\r u w)_t+\big((\r u^2+P(\rho))w\big)_x\nonumber\\
&=(\rho u^2+P(\rho))w_x+\varepsilon(\rho^{\alpha}w u_x)_x-\varepsilon\rho^{\alpha}u_xw_x+\lambda \rho uw+\rho Vw-\rho\partial_xW\ast\rho w.
\end{align}
Integrating \eqref{4.40} over $[r_1,x)$ to obtain
\begin{align}\label{4.41}
(\rho u^2+P(\rho))w
&=\varepsilon\rho^{\alpha}wu_x+\int^{x}_{r_1}\Big((\rho u^2+P(\rho))w_y-\varepsilon
\rho^{\alpha}u_yw_y\Big)\,\dd y-\frac{\dd}{\dd t}\int^x_{r_1}\rho uw\,\dd y\nonumber\\
&\quad+\int^x_{r_1}\lambda \rho uw\,\dd y-\int^x_{r_1}\rho w \partial_x W\ast\rho \,\dd y+\int^x_{r_1}\rho w V\,\dd y.
\end{align}
Multiplying \eqref{4.41} by $\rho w$ and performing a direct calculation, we obtain
\begin{align}\label{4.42}
\r P(\rho) w^2&=\v \rho^{\alpha+1}w^2u_x-\Big(\rho w\int^x_{r_1}\rho uw\,\dd y\Big)_t-\Big(\rho uw\int^x_{r_1}\rho uw\,\dd y\Big)_x\nonumber\\
&\quad+\rho uw_x\int^x_{r_1}\rho uw\,\dd y+\rho w\int^x_{r_1}\Big((\rho u^2+P(\rho))w_y-\varepsilon\rho^{\alpha}u_yw_y\Big)\,\dd y\nonumber\\
&\quad+\lambda \rho w\int^x_{r_1}\rho u w\,\dd y-\rho w\int^x_{r_1}\rho\partial_xW\ast\rho w\,\dd y+\rho w\int^x_{r_1}\rho Vw\,\dd y:=\sum^{8}_{i=1}K_{i}.
\end{align}

\noindent{\emph{Step 2.}}  To estimate $K_i, i=1,\cdots, 8$, in \eqref{4.42}, we first notice that
\begin{align}
\int_{b^{-}(t)}^{b^+(t)}\rho|u|\, \dd x\leq \int_{b^-(t)}^{b^+(t)} (\rho+\rho u^2)\,\dd x\leq C.\nonumber
\end{align}
Then it follows from \eqref{3.7} and \eqref{4.1} that
\begin{align}
\left| \int_0^T\int_{r_1}^{r_2} K_2\, \dd x\dd t \right|&=\left| \int_0^T\int_{r_1}^{r_2} \Big(\rho w\int^x_{r_1}\rho uw\,\dd y\Big)_t\, \dd x\dd t \right|\nonumber\\
&\leq\left|\int^{r_2}_{r_1}\Big(\rho w\int^x_{r_1}\rho u w\,\dd y\Big)(T,x)\,\dd x\right|+\left|\int^{r_2}_{r_1}\Big(\rho w\int^x_{r_1}\rho u w\,\dd y\Big)(0,x)\,\dd x\right|\leq C.\nonumber
\end{align}
Similarly  to  \cite[Lemma 3.5]{Chen2021}, we obtain
\begin{align}
\left|\int_0^T\int_{r_1}^{r_2} K_1\,\dd x\dd t \right|&\leq C(r_1,r_2)
+\v\int^T_0\int^{r_2}_{r_1}\rho^{\gamma+1}w^2\,\dd x\dd t,\nonumber\\
\left| \int_0^T\int_{r_1}^{r_2}K_3\, \dd x\dd t \right|&=\left| \int_0^T\int_{r_1}^{r_2}
\Big(\rho uw\int^x_{r_1}\rho uw\,\dd y\Big)_x\, \dd x\dd t \right|=0,\nonumber\\
\left| \int_0^T\int_{r_1}^{r_2} K_4\, \dd x\dd t \right|&=\left| \int_0^T\int_{r_1}^{r_2} \Big(\rho uw_x\int^x_{r_1}\rho uw\,\dd y\Big)\, \dd x\dd t \right|\leq C.\nonumber
\end{align}
Using the fact that $\alpha\leq\gamma,$ we have
\begin{align}
\left| \int_0^T\int_{r_1}^{r_2} K_5\, \dd x\dd t \right|&=\left| \int_0^T\int_{r_1}^{r_2} (\rho uw\int^x_{r_1}\Big((\rho u^2+P(\rho))w_y-\v\rho^{\alpha}w_yu_y\,\dd y\Big)\, \dd x\dd t \right|\nonumber\\
&\leq C
  +\left|\v\int^T_0\int^{r_2}_{r_1}\rho w\int^x_{r_1}\rho^{\alpha}w_yu_y\,\dd x\dd t\right|\nonumber\\
&\leq C(r_1,r_2),\nonumber
\end{align}
where we have used the following estimate
\begin{align}
\left|\v\int_0^T\int_{r_1}^{r_2} \rho w\int^x_{r_1}\rho^{\alpha}w_yu_y\, \dd y\dd x\dd t \right|&\leq \v\int^T_0\int^{r_2}_{r_1}|\rho^{\alpha}w_xu_x|\,\dd x\dd t\nonumber\\
&\leq \Big(\v\int^T_0\int^{r_2}_{r_1}\rho^{\alpha}|u_x|^2\,\dd x \dd t\Big)^{\frac{1}{2}}\Big(\v\int^T_0\int^{r_2}_{r_1}\rho^{\alpha}|w_x|^2\,\dd x\dd t\Big)^{\frac{1}{2}}\nonumber\\
&\leq C\Big(\v\int^T_0\int^{r_2}_{r_1}\rho^{\alpha}\,\dd x\dd t\Big)^{\frac{1}{2}}
\leq C\Big(\v\int^T_0\int^{r_2}_{r_1}\rho^{\gamma}\,\dd x\dd t\Big)^{\frac{1}{2}}\nonumber\\
&\leq C(r_1,r_2).\nonumber
\end{align}
Moreover, we have
\begin{align}
\left|\int_0^T\int_{r_1}^{r_2} K_6\,\dd x\dd t \right|&\leq C(,\nonumber
\end{align}
and
\begin{align*}
&\int_0^T\int_{r_1}^{r_2} K_7\,\dd x\dd t \\
&\,\,\,=\int_0^T\int_{r_1}^{r_2} \rho(t,x) w(x)\Big(\int^x_{r_1}\rho(t,y) (\partial_xW\ast\rho)(t,y) w(y) \,\dd y\Big) \dd x\dd t \nonumber\\
&\,\,\,= \int_0^T\int_{r_1}^{r_2} \rho(t,x) w(x)\int^x_{r_1}\rho(t,y)
\Big( \int^{b^+(t)}_{b^-(t)}\big(1-2H(x-z)+x-z\big)\rho(t,z)\,\dd z\Big)w(y)\,\dd y\dd x\dd t \nonumber\\[2mm]
&\int_0^T\int_{r_1}^{r_2} K_8\,\dd x\dd t\nonumber\\
&\,\,\, =\int^T_0\int^{r_2}_{r_1}\rho(t,x)w(x)\int^x_{r_1}\rho(t,y)
\Big(\int_{b^{-}(t)}^{b^+(t)}\varpi(y-z)\big(u(z)-u(y)\big)\rho(t,z)\,\dd z\Big)w(y)\,\dd y\dd x\dd t,\nonumber\\
\end{align*}
which are clearly bounded by $C$.
Integrating \eqref{4.42} over $[0,T]\times[r_1, r_2],$ and collecting all the estimates in this step,
we can obtain \eqref{4.39}.
The proof
is completed.
$\hfill\Box$

\medskip
\section{Inviscid Limit for Polytropic Gases}\label{section4}
To prove Theorem \ref{thm3.3}, we intend to apply the compensated compactness framework
from Chen-Perepelitsa \cite{Perepelitsa}.
For clarity of presentation, we first focus on the polytropic case \eqref{pressure1}, while
the general pressure case will be discussed in \S 5 below.

We first explore some important properties of several special entropy entropy flux pairs.

\subsection{Choice of a special entropy and entropy flux pair}

Taking $\psi(s)=\frac12 {s}{|s|}$ in \eqref{weakentropy}, the corresponding entropy and entropy flux are represented as
\begin{align}\label{4.52}
\begin{cases}
\displaystyle\eta^{\#}(\r,\r u)=\f12 \r \int_{-1}^1 (u+\r^{\t} s) |u+\r^{\t}s| [1-s^2]_+^{\fb} \dd s,\\[3mm]
\displaystyle q^{\#}(\r, \r u)=\f12 \r \int_{-1}^1 (u+\theta\r^{\t}s)(u+\r^{\t} s) |u+\r^{\t}s| [1-s^2]_+^{\fb} \dd s.
\end{cases}
\end{align}
A direct calculation shows that
\begin{align}\label{4.53}
|\eta^{\#}(\r,\r u)|\leq C_{\g} \big(\r |u|^2+\r^{\g}\big), \qquad\,\, q^{\#}(\r,\r u)\geq C_{\g}^{-1} \big(\r |u|^3+\r^{\g+\t}\big),
\end{align}
where and whereafter $C_\g>0$ is a universal constant depending only on $\g>1$.
We regard $\eta^{\#}$ as a function of $(\rho,m)$ to obtain
\begin{align*}
\eta^{\#}_\r=\int_{-1}^1 \big(-\f12 u+(\t+\f12)\r^{\t} s\big)\,|u+\r^{\t}s| [1-s^2]_+^{\fb} \dd s,\qquad
\displaystyle \eta^{\#}_m=\int_{-1}^1|u+\r^{\t}s| [1-s^2]_+^{\fb} \dd s.
\end{align*}
It is direct to check that
\begin{align}\label{4.54}
|\eta^{\#}_m|\leq C_\g \big(|u|+\rho^\theta\big),\qquad\,\, |\eta^{\#}_\rho|\leq C_\g\big(|u|^2+\rho^{2\theta}\big).
\end{align}

\subsection{Higher integrability of the velocity}
The special entropy pair introduced in \S 4.1 allows us to derive a better estimate for
the integrability of the velocity vector field.

The following lemma is important to control the trace estimates for the higher integrability on the velocity; see \eqref{4.63}.
In fact, we have the boundary parts $(u \eta^{\#})(t,b^\pm(t))$ and $q^{\#}(t,b^\pm(t))$,
and it is impossible to have the uniform trace bound (independent of $\v$) for each of them.

\begin{lemma}[\cite{Chen2021}, Lemma 3.6]\label{lem2.5}
For the entropy pair defined in \eqref{4.52}, the following cancelation property holds:
\begin{align}\label{4.55}
|q^{\#}-u\eta^{\#}|\leq C_\g \big(\rho^{\gamma} |u|+ \rho^{\gamma+\theta} \big).
\end{align}
\end{lemma}

\medskip
\begin{lemma}\label{lem4.2}
Let $(\r,u)$ be the smooth solution of  \eqref{1.8} and \eqref{3.1}--\eqref{3.4},
and let the assumption of {\rm Lemma \ref{lem4.4}} hold.
Then, for any $(r_1,r_2)\Subset [b^-(t), b^+(t)]$, there exists $C(r_1,r_2)>0$ indepedent of $\v\in (0,1]$
such that
\begin{align}\label{4.56}
\int_0^T\int_{r_1}^{r_2}  \big(\rho|u|^3+\rho^{\gamma+\theta}\big)(t,x)\, \dd x \dd t\leq C(r_1,r_2).
\end{align}
\end{lemma}

\noindent{\bf Proof.}  We divide the proof into four steps.

\smallskip
\noindent{\emph{Step 1.}} Multiplying $\eqref{1.8}_1$ by $\eta^{\#}_\r$ and $\eqref{1.8}_2$ by $ \eta^{\#}_m$,
we have
\begin{align}\label{4.57}
& \eta^{\#}_t+q^{\#}_x = \eta^{\#}_m\, \Big(\v (\rho^{\alpha} u_x)_x+\lambda \rho u+\rho V-\rho \partial_xW\ast\rho\Big).
\end{align}
Using \eqref{2.6}, a direct calculation shows that
\begin{align}
\frac{\dd}{\dd t}\int_x^{b^+(t)}\eta^{\#}\,\dd y
& =\eta^{\#}(t,b^+(t)) \frac{\dd }{\dd t}b^+(t)+\int_x^{b^+(t)}\eta^{\#}_t(t,y)\,\dd y= (u\eta^{\#})(t, b^+(t))+\int_x^{b^+(t)} \eta^{\#}_t(t,y)\,\dd y.\nonumber
\end{align}
Integrating \eqref{4.57} over $[x,b^+(t))$,  we have
\begin{align}\label{4.59}
 q^{\#}(t,x) &=\Big( \int_x^{b^+(t)}\eta^{\#}(t,y)\,\dd y\Big)_t+\big(q^{\#}-u\eta^{\#}\big)(t,b^+(t))\nonumber\\
 &\quad-\v \int_x^{b^+(t)}  \eta^{\#}_m(\rho^{\alpha} u_y)_y\,\dd y-
 \lambda \int_x^{b^+(t)}  \eta^{\#}_m\rho u\,\dd y+\int_x^{b^+(t)} \eta^{\#}_m \rho\partial W\ast\rho\,\dd y-\int^{b^+(t)}_x\eta^{\#}\rho V\,\dd y\nonumber\\
 &=\sum^6_{i=1}I_i.
\end{align}
We now bound each term of the RHS of \eqref{4.59}.

\smallskip
\noindent{\emph{Step 2.}}  First, for the term involving the trace estimates (second term) in \eqref{4.59}, it follows from \eqref{4.55} and
Lemmas \ref{lem4.1}, \ref{lem2.2}, and \ref{lem2.5}
that
\begin{align}\label{4.60}
&\int_0^T\int^{r_2}_{r_1}
\big|(q^{\#}-u\eta^{\#})(t,b^+(t))\big|\,
\dd x\dd t\leq C(r_1,r_2)\int_0^T\big(\rho^{\gamma+\theta}(t,b^+(t))+(\rho^{\gamma}|u|)(t, b^+(t))\big)\,\dd t.
\end{align}
It follows from \eqref{4.29} that
\begin{align}\label{4.61}
&\int_0^T\big(\rho(t,b^+(t)))^{\gamma+\theta}\dd t
=\int^T_0\Big(\rho_0(b)\big(1
  +\frac{\kappa(\gamma-\alpha)}{\varepsilon}\rho^{\gamma-\alpha}_0(b)t\big)^{-\frac{1}{\gamma-\alpha}}\Big)^{\gamma+\theta}\,\dd t
  \leq (\rho_0(b))^{\gamma+
\theta}T\leq C.
\end{align}
Similar with argument in \cite{Chen2021,He}, we have
\begin{align}\label{4.62}
&\int_0^T\int^{r_2}_{r_1}
\big(\rho^{\gamma}|u|\big)(t,b^+(t))\,
\dd x\dd t\leq C\v^{\frac{p(\gamma-\alpha)-\gamma}{2\gamma}}
\leq C,
\end{align}
where $p>\frac{\gamma}{\gamma-\alpha}$.
Substituting \eqref{4.61}--\eqref{4.62} into \eqref{4.60}, we obtain
\begin{align}\label{4.63}
&\Big|\int^T_0\int^{r_2}_{r_1}I_2\,\dd x\dd t\Big|=\Big|\int^T_0\int^{r_2}_{r_1}
(q^{\#}-u\eta^{\#})(t,b^+(t))\,
\dd x\dd t\Big|\leq C(r_1,r_2).
\end{align}

For the first term of the RHS of \eqref{4.59}, using \eqref{4.1} and \eqref{4.30}, we see that, for $I_1,$
\begin{align}\label{4.64}
&\Big|\int^T_0\int^{r_2}_{r_1}I_1\,\dd x\dd t\Big|=\Big|\int_0^T\int_{r_1}^{r_2}\Big(\int_x^{b^+(t)}\eta^{\#}(\rho,\rho u)\, \dd y\Big)_t\,\dd x\dd t\Big| \nonumber\\
&\leq \Big|\int_{r_1}^{r_2}\int_{b^-(t)}^{b^+(t)}\eta^{\#}(\rho,\rho u)(T,y)\, \dd y\dd x\Big|
  +\Big|\int_{r_1}^{r_2}\int_{b^-(t)}^{b^+(t)}\eta^{\#}(\rho_0,\rho_0 u_0)\, \dd y\dd x\Big|\nonumber\\[2mm]
&\leq C(r_1,r_2).
\end{align}

\noindent{\emph{Step 3.}} For $I_3,$ we integrate by parts to obtain
\begin{align}
-\v\int_x^{b^+(t)}\eta^{\#}_m(\rho^{\alpha} u_y)_y\,\dd y    = &\,
-\v\Big(\eta^{\#}_m(t,b^+(t))\, (\rho^{\alpha} u_x)(t,b^+(t))-\eta^{\#}_m(t,x)(\rho^{\alpha} u_x)(t,x) \Big) \nonumber\\
&\,\,+\v\int_x^{b^+(t)}\rho^{\alpha} u_y(\eta^{\#}_{mu}u_y+\eta^{\#}_{m\rho}\rho_y)\,\dd y:= J_1+J_2.\label{4.65}
\end{align}
Now, we discuss $J_2$ first:
\begin{align*} 
|J_2|&=\Big|\v\int^{b^+(t)}_x\Big(\eta^{\#}_{m\rho}\rho^{\alpha}u_y\rho_y+\eta^{\#}_{mu}\rho^{\alpha}u^2_{y}\Big)\,\dd y\Big|\leq C\v\Big|\int^{b^+(t)}_x\rho^{\theta+\alpha-1}u_y\rho_y\,\dd y\Big|+\v\int^{b^+(t)}_x\rho^{\alpha}u^2_{y}\,\dd y\nonumber\\
&\leq \v\int^{b^+(t)}_x\rho^{\gamma+\alpha-3}\rho^2_y\,\dd y+\v\int^{b^+(t)}_x\rho^{\alpha}u^2_{y}\,\dd y,
\end{align*}
where we have used the fact that $|\eta^{\#}_{mu}|\leq C$ and $|\eta^{\#}_{m\r}|\leq C\r^{\t-1}$.
Then we obtain
\begin{align}\label{4.68}
\int^T_0\int^{r_2}_{r_1}|J_2|\,\dd x\dd t&\leq\varepsilon\int^T_0\int^{r_2}_{r_1}\int^{b^+(t)}_x\rho^{\gamma+\alpha-3}\rho^2_y\,\dd y\dd x\dd t+\varepsilon\int^T_0\int^{r_2}_{r_1}\int^{b^+(t)}_x\rho^{\alpha}u^2_y\,\dd y\dd x\dd t\nonumber\\
&\leq C(r_1,r_2).
\end{align}
For $J_1,$ we have
\begin{align}\label{4.69}
\varepsilon\Big|\int^T_0\int^{r_2}_{r_1} \eta^{\#}_m\rho^{\alpha}u_x\,\dd x\dd t\Big|
&\leq\varepsilon\int^T_0\int^{r_2}_{r_1}\rho^{\alpha}u^2_x\,\dd x\dd t+\varepsilon\int^T_0\int^{r_2}_{r_1}\rho^{\alpha}(|u|+\rho^{\theta})^2\,\dd x\dd t\nonumber\\
&\leq C
+\v\int^T_0\int^{r_2}_{r_1}\rho^{2\theta+\alpha}\,\dd x\dd t+\v\int^T_0\int^{r_2}_{r_1}\rho^{\alpha}u^2\,\dd x\dd t.
\end{align}
Similar to the argument as in \cite{Chen2021,He}, to control \eqref{4.69}, we need that
\begin{align}\label{4.70}
\v\int^T_0\int^{r_2}_{r_1}\rho^{2\theta+\alpha}\,\dd x\dd t
&=\v\int^T_0\int^{r_2}_{r_1}\rho^{\gamma-1+\alpha}\,\dd x\dd t
\leq C\int^T_0\int^{r_2}_{r_1}\rho^{\frac{\gamma-1}{2}}\,\dd x\dd t\nonumber\\
&\leq C\int^T_0\int^{r_2}_{r_1}\rho^{\gamma}\,\dd x\dd t+C(r_1,r_2)
\nonumber\\
&\leq C(r_1,r_2).
\end{align}

It follows from \eqref{4.1} and \eqref{4.30} that, for $\beta=\alpha+\frac{\gamma-1}{2},$
\begin{align}\label{4.71}
\v \rho^{\beta}(t,x)&=\v\Big(\rho^{\beta}(t,x)-\rho^{\beta}(t,b^-(t))+\rho^{\beta}(t,b^-(t))\Big)\nonumber\\
&\leq\v \beta\int^{b^+(t)}_{b^-(t)}\rho^{\beta-1}|\rho_x|\,\dd x+\v\rho^{\beta}_0(b)\nonumber\\
&\leq\beta\Big(\v^2\int^{b^+(t)}_{b^-(t)}\rho^{2\alpha-3}\rho^2_x \,\dd x\Big)^{\frac{1}{2}}\Big(\int^{b^+(t)}_{b^-(t)}\rho^{2(\beta-\alpha)+1}\,\dd x\Big)^{\frac{1}{2}}+\v\rho^{\beta}_0(b)\nonumber\\
&\leq C\Big(\int^{b^+(t)}_{b^-(t)}\rho^{2(\beta-\alpha)+1}\,\dd x\Big)^{\frac{1}{2}}
+C
\leq C.
\end{align}

Using \eqref{4.71}, we have
\begin{align}\label{4.72}
\v\int^T_0\int^{r_2}_{r_1}\rho^{\alpha}u^2\,\dd x\dd t&\leq \Big(\int^T_0\int^{r_2}_{r_1}\rho|u|^3\,\dd x\dd t\Big)^{\frac{2}{3}}\,\Big(\int^T_0\int^{r_2}_{r_1}\varepsilon^3\rho^{3\alpha-2}\,\dd x\dd t\Big)^{\frac{1}{3}}\nonumber\\
&\leq \Big(\int^T_0\int^{r_2}_{r_1}\rho|u|^3\,\dd x\dd t\Big)^{\frac{2}{3}}\,\Big(\int^T_0\int^{r_2}_{r_1}\varepsilon^3\rho^{3\beta}\,\dd x\dd t+C(r_1,r_2)\Big)^{\frac{1}{3}}\nonumber\\
&\leq C(r_1,r_2)
\Big(\int^T_0\int^{r_2}_{r_1}\rho|u|^3\,\dd x\dd t\Big)^{\frac{2}{3}},
\end{align}
where we have used the assumption: $\alpha\geq\frac{2}{3}.$
Inserting \eqref{4.70} and \eqref{4.72} into \eqref{4.69}, we find that, for $\delta>0,$
\begin{align}\label{4.73}
&\varepsilon\Big|\int^T_0\int^{r_2}_{r_1} \eta^{\#}_m\rho^{\alpha}u_x\,\dd x\dd t\Big|
\leq C(r_1,r_2)
\frac{1}{\delta}+\delta\int^T_0\int^{r_2}_{r_1}\rho|u|^3\,\dd x\dd t.
\end{align}

Using \eqref{4.54} and \eqref{3.3}, we have
$$
|\v(\eta^{\#}_m\rho^{\alpha}u_x)(t,b^+(t))|=\kappa|(\eta^{\#}_m\rho^{\gamma})(t,b^+(t))|
\leq C_{\gamma}\big(\rho^{\gamma}|u|+\rho^{\gamma+\theta}\big)(t,b^+(t)).
$$
Similar again to the argument as in \cite{Chen2021,He}, we obtain
\begin{align}\label{4.74}
\int^T_0\int^{r_2}_{r_1}|\v(\eta^{\#}_m\rho^{\alpha}u_x)(t,b^+(t))|\,\dd x\dd t&\leq C(r_1,r_2)
\int^T_0(\rho^{\gamma}|u|+\rho^{\gamma+\theta})(t,b^+(t))\,\dd t\nonumber\\
&\leq C(r_1,r_2).
\end{align}
Combining \eqref{4.73} and \eqref{4.74} yields
\begin{align}\label{4.75}
\Big|\int^T_0\int^{r_2}_{r_1}J_1\,\dd x\dd t\Big|&\leq \delta
\int^T_0\int^{r_2}_{r_1}\rho |u|^3\,\dd x\dd t
+C(r_1,r_2).
\nonumber\\
\end{align}
Inserting \eqref{4.75} and \eqref{4.68} into \eqref{4.65}, we obtain
\begin{align}\label{4.76}
\Big|\int^T_0\int^{r_2}_{r_1}I_3\,\dd x\dd t\Big|&\leq \delta
\int^T_0\int^{r_2}_{r_1}\rho  |u|^3\,\dd x\dd t
+C(r_1,r_2).
\end{align}

\smallskip
\noindent{\emph{Step 4.}} We also obtain the estimates for $I_4$ and $I_5:$
\begin{align*}
\Big|\int^T_0\int^{r_2}_{r_1}I_4\,\dd x\dd t\Big|
&=\Big|\int^T_0\int^{r_2}_{r_1}\int^{b^+(t)}_x\rho |u|(|u|+\rho^{\theta}) \,\dd y\dd x\dd t\Big|\nonumber\\
&=\Big|\int^T_0\int^{r_2}_{r_1}\rho u^2\,\dd x\dd t\Big|+\Big|\int^T_0\int^{r_2}_{r_1}\rho^{\gamma}\,\dd x\dd t\Big|
\leq C(r_1,r_2),
\end{align*}
\begin{align*}
&\Big|\int^T_0\int^{r_2}_{r_1}I_5\,\dd x\dd t\Big|\nonumber\\
&=\Big|\int^T_0\int^{r_2}_{r_1}\int^{b^+(t)}_x\rho \eta^{\#}_m\partial_xW\ast \rho \,\dd y\dd x\dd t\Big|\nonumber\\
&=\Big|\int^T_0\int^{r_2}_{r_1}\int^{b^+(t)}_x\rho\eta^{\#}_m\Big(M-2\int^{y}_{b^-(t)}\rho(t,z)\,\dd z+yM-\int^{b^+(t)}_{b^{-}(t)}z\rho(t,z) \,\dd z\Big)\,\dd y\dd x\dd t\Big|
\leq C(r_1,r_2).
\end{align*}
\begin{align}\label{4.782}
&\Big|\int^T_0\int^{r_2}_{r_1}I_6\,\dd x\dd t\Big|\nonumber\\
&=\Big|\int^T_0\int^{r_2}_{r_1}\int^{b^+(t)}_x\eta^{\#}_m\rho V\,\dd y\,\dd x\dd t\Big|\nonumber\\
&=\Big|\int^T_0\int^{r_2}_{r_1}\int^{b^+(t)}_x\eta^{\#}_m\rho
  \Big(\int^{b^+(t)}_{b^-(t)}\varpi(y-z)(u(y)-u(z))\rho(t,z)\,\dd z\Big)\,\dd y\,\dd x\dd t\Big|
  \leq C(r_1,r_2).
\end{align}

Combining estimates \eqref{4.63}--\eqref{4.64} with estimates \eqref{4.76}--\eqref{4.782},
we conclude \eqref{4.56} from \eqref{4.53} and \eqref{4.59}.
$\hfill\Box$

\medskip
\subsection{$H_{\rm loc}^{-1}(\mathbb{R}_+^2)$--Compactness}\label{section4.2}

In this section we use the uniform estimates obtained in \S 4.2 to prove the following key lemma, which states the
$H_{\rm loc}^{-1}(\mathbb{R}_+^2)-$compactness of the dissipation measures for the approximate solutions.

\begin{lemma}\label{lemma4.8}
Let $\alpha\in [\frac23, \gamma]$, and let
$(\eta^{\psi},q^{\psi})$ be the weak entropy pair generated by any $\psi\in C_0^2(\mathbb{R})$,
defined in \eqref{weakentropy}.
Then, for the solution sequence $(\rho^{\varepsilon},u^{\varepsilon})$ with
$m^{\varepsilon}=\rho^{\varepsilon}u^{\varepsilon}$ of CNSEs \eqref{1.8} and \eqref{3.1}--\eqref{3.4},
the following sequence
\begin{eqnarray}\nonumber
\eta^{\psi}(\rho^{\varepsilon},m^{\varepsilon})_t+q^{\psi}(\rho^{\varepsilon},m^{\varepsilon})_x
\qquad \mbox{is compact in $H_{\rm loc}^{-1}(\mathbb{R}_+\times \mathbb{R} )$.}
\end{eqnarray}
\end{lemma}

\noindent\textbf{Proof}. To prove this lemma, we first recall the following results
for the entropy pair
$(\eta^{\psi},q^{\psi})$ generated by $\psi\in C_0^2(\mathbb{R})$;
see also \cite{Chen6,Perepelitsa} for details.

For a $C^2$--function $\psi:\mathbb{R}\rightarrow\mathbb{R}$,
compactly supported on the interval $[a,b]$, we have
\begin{eqnarray}\nonumber
{\rm supp}(\eta^{\psi}),\,{\rm supp}(q^{\psi})\subset \left\{{(\rho,m)=(\rho,\rho
u)\,:\, u+\rho^{\theta}\geq a,\quad u-\rho^{\theta}\leq b}\right\}.
\end{eqnarray}
Furthermore, from \cite[Lemma 2.1]{Perepelitsa}, there exists $C_{\psi}>0$ such that, for any
$\rho\geq0$ and $u\in\mathbb{R}$, we have the following facts:
\begin{enumerate}
\item [(\rmnum{1})] For $\gamma\in(1,3]$,
\begin{align}\label{4.79}
|\eta^{\psi}(\rho,m)|+|q^{\psi}(\rho,m)|\leq C_{\psi}\rho.
\end{align}

\item [(\rmnum{2})] For $\gamma\in(3,\infty)$,
\begin{align}\label{4.80}
|\eta^{\psi}(\rho,m)|\leq C_{\psi}\rho,\quad \ |q^{\psi}(\rho,m)|\leq
C_{\psi}(\rho+\rho^{1+\theta}).
\end{align}

\item [(\rmnum{3})] If $\eta^{\psi}$ is considered as a function of $(\rho,m)$,
$m=\rho u$, then
\begin{align}\label{4.81}
|\eta^{\psi}_m(\rho,m)|\leq C_{\psi},\quad\, |\eta^{\psi}_{\rho}(\rho,m)|\leq C_{\psi}(1+\rho^{\theta}),
\end{align}
 and, if $\eta^{\psi}_m$ is considered as a
function of $(\rho, u)$, then
\begin{align}\label{4.82}
|\eta^{\psi}_{mu}(\rho,\rho
u)|+|\rho^{1-\theta}\eta^{\psi}_{m\rho}(\rho,\rho u)|\leq C_{\psi}.
\end{align}
\end{enumerate}
Now we are going to prove Lemma \ref{lemma4.8}.

A direct computation on $\eqref{1.8}_1\times\eta^{\psi}_{\rho}(\rho^{\varepsilon},m^{\varepsilon})
+\eqref{1.8}_2\times\eta^{\psi}_m(\rho^{\varepsilon},m^{\varepsilon})$
gives
\begin{align}\label{4.83}
\displaystyle\eta^{\psi}(\rho^{\varepsilon},m^{\varepsilon})_t+q^{\psi}(\rho^{\varepsilon},m^{\varepsilon})_x
=&\,\varepsilon\big(\eta^{\psi}_m(\rho^{\varepsilon},m^{\varepsilon})(\rho^{\varepsilon})^{\alpha}u_x^{\varepsilon}\big)_x
-\varepsilon\eta^{\psi}_{mu}(\rho^{\varepsilon},m^{\varepsilon})(\rho^{\varepsilon})^{\alpha}(u_x^{\varepsilon})^2\nonumber\\
&\,-\varepsilon\eta^{\psi}_{m\rho}(\rho^{\varepsilon},m^{\varepsilon})(\rho^{\varepsilon})^{\alpha}\rho_x^{\varepsilon}u_x^{\varepsilon}
 +\eta^{\psi}_m\big(\lambda\rho^{\v}u^{\v}+\rho^{\v}V-\rho^{\v}\partial_xW\ast\rho^{\v}\big).\nonumber
\end{align}
For any compact set $K\Subset[b^-(t),b^+(t)]$, using \eqref{4.82} and the Cauchy-Schwartz
inequality, we have
\begin{equation}
\begin{aligned}
\displaystyle\int_{0}^{T}\!\!\!\int_{K} & \big|\eta^{\psi}_{mu}(\rho^{\varepsilon},m^{\varepsilon})(\rho^{\varepsilon})^{\alpha}
(u_x^{\varepsilon})^2
+\eta^{\psi}_{m\rho}(\rho^{\varepsilon},m^{\varepsilon})(\rho^{\varepsilon})^{\alpha}\rho_x^{\varepsilon}u_x^{\varepsilon}\big|\,\dd x\dd t\\
&\displaystyle\leq
C_{\psi}\int_{0}^{T}\!\!\!\int_{K}(\rho^{\varepsilon})^{\alpha}(u_x^{\varepsilon})^2\, \dd x\dd t
+C_{\psi}\int_{0}^{T}\!\!\!\int_{K}(\rho^{\varepsilon})^{\alpha+\gamma-3}(\rho_x^{\varepsilon})^2\, \dd x\dd t
\leq C,
\nonumber\\[2mm]
\displaystyle\int_{0}^{T\!\!\!}\int_{K} & \big|\eta^{\psi}_m
  \big(\lambda\rho^{\v}u^{\v}+\rho^{\v} V-\rho^{\v}\partial_xW\ast\rho^{\v}\big)\big|\,\dd x\dd t\\
&\displaystyle\leq
C_{\psi}\int_{0}^{T}\!\!\!\int_{K}\rho^{\v}(u^{\v})^2+\rho^{\v} \,\dd x\dd t
+C_{\psi}\int_{0}^{T}\!\!\!\int_{K}|\rho^{\v}\partial_xW\ast\rho^{\v}| \,\dd x\dd t \leq C(K).
\nonumber
\end{aligned}
\end{equation}
This implies that
\begin{align}\label{4.86}
&-\varepsilon\eta^{\psi}_{mu}(\rho^{\varepsilon},m^{\varepsilon})(\rho^{\varepsilon})^{\alpha}(u_x^{\varepsilon})^2
-\varepsilon\eta^{\psi}_{m\rho}(\rho^{\varepsilon},m^{\varepsilon})(\rho^{\varepsilon})^{\alpha}\rho_x^{\varepsilon}u_x^{\varepsilon}
+\eta^{\psi}_m(\lambda\rho^{\v}u^{\v}+\rho^{\v}V-\rho^{\v}\partial_xW\ast\rho^{\v})\nonumber\\
&\,\,\mbox{is uniformly bounded in $L^1([0,T]\times K)$},
\end{align}
so that it is compact in
$W_{\rm loc}^{-1,p_1}(\mathbb{R}_+^2)$ for $1<p_1<2$.

If $2\alpha\leq\gamma+1,$ then
\begin{align}\label{4.87}
\v^{\frac{4}{3}}\int^T_0\int_{K}(\rho^{\v})^{2\alpha}\,\dd x\dd t\leq C(K)
\v^{\frac{4}{3}}.
\end{align}
If $2\alpha\geq\gamma+1,$ $\alpha<\gamma,$ we see from \cite{He} that
\begin{align}\label{4.88}
\v^{\frac{4}{3}}\int^T_0\int_{K}(\rho^{\v})^{2\alpha}\,\dd x\dd t
&\leq C(K)
\v^{\frac{1}{3}}\int^T_0\int_K(\rho^{\v})^{\alpha-\frac{\gamma-1}{2}}\,\dd x\dd t\nonumber\\
&\quad+C(K )\v^{\frac{1}{3}}\int^T_0\int_K(\rho^{\v})^{\gamma+1}\,\dd x\dd t.
\end{align}
It follows from \eqref{4.81} and \eqref{4.87}--\eqref{4.88} that
\begin{equation}\label{4.89}
\begin{aligned}
\displaystyle\int_{0}^{T}\int_{K}\Big(\varepsilon
\eta^{\psi}_m(\rho^{\varepsilon},m^{\varepsilon})
(\rho^{\varepsilon})^{\alpha}u_x^{\varepsilon}\Big)^{\frac{4}{3}}\, \dd x\dd t
&\displaystyle\leq\int_{0}^{T}\int_{K}\varepsilon^{\frac{4}{3}}
(\rho^{\varepsilon})^{^{\frac{4\alpha}{3}}}|u_x^{\varepsilon}|^{\frac{4}{3}}\,\dd x\dd t\\
&\displaystyle\leq C\varepsilon^{\frac{4}{3}}\int_{0}^{T}\int_{K}
(\rho^{\varepsilon})^{\alpha}|u_x^{\varepsilon}|^2\,\dd x\dd t+C\varepsilon^{\frac{4}{3}}\int_{0}^{T}\int_{K}
(\rho^{\varepsilon})^{2\alpha}\,\dd x\dd t\\
 &\displaystyle\leq C(K)\varepsilon^{\frac{1}{3}}+C\varepsilon^{\frac{4}{3}}\int_{0}^{T}\int_{K}
(\rho^{\varepsilon})^{\gamma+1}\,\dd x\dd t\\
&\displaystyle\leq C(K)\varepsilon^{\frac{1}{3}}\rightarrow  0\qquad \ \text{as $\varepsilon\rightarrow0^+$}.
\end{aligned}
\end{equation}
Then \eqref{4.86} and \eqref{4.89} yield
\begin{align}\label{4.90}
\eta^{\psi}(\rho^{\varepsilon},m^{\varepsilon})_t+q^{\psi}(\rho^{\varepsilon},m^{\varepsilon})_x
\qquad  \mbox{is compact in $W_{\rm loc}^{-1,p_2}(\mathbb{R}_+^2)$} \ \
\mbox{for some $1<p_2<2$}.
\end{align}

Furthermore, for $\gamma\in (1,3],$ using \eqref{4.79} and Lemma \ref{lem4.5}, we have
\begin{align}
\int^T_0\int_K\big(|\eta^{\psi}(\rho^{\varepsilon},m^{\varepsilon})|+|q^{\psi}(\rho^{\varepsilon},m^{\varepsilon})|\big)^{\gamma+1}\,\dd x\dd t
\leq C\int^T_0\int_K(\rho^{\v})^{\gamma+1}\,\dd x\dd t\leq C(K).
\nonumber
\end{align}
For $\gamma\in (3, \infty),$ using \eqref{4.80} and Lemma \ref{lem4.2}, we have
\begin{align}
\int^T_0\!\!\!\int_K\big(|\eta^{\psi}(\rho^{\varepsilon},m^{\varepsilon})|
  +|q^{\psi}(\rho^{\varepsilon},m^{\varepsilon})|\big)^{\frac{\gamma+\theta}{\theta+1}}\,\dd x\dd t
\leq C\int^T_0\!\!\!\int_K\big(|\rho^{\varepsilon}|^{\frac{\gamma+\theta}{\theta+1}}
  +|\rho^{\v}|^{\gamma+\theta}\big)\,\dd x\dd t\leq C(K).
  \nonumber
\end{align}
Thus, using the last two estimates, we obtain
\begin{align}
(\eta^{\psi}(\rho^{\varepsilon},m^{\varepsilon}),\, q^{\psi}(\rho^{\varepsilon},m^{\varepsilon}))
\qquad\, \mbox{is uniformly bounded in $L_{\rm loc}^{p_3}(\mathbb{R}_+^2)$} \ \
\mbox{for}\ p_3>2,\nonumber
\end{align}
where $p_3=\gamma+1>2$ when $\gamma\in(1,3]$, and
$p_3=\frac{\gamma+\theta}{1+\theta}>2$ when $\gamma\in(3,\infty)$. This implies that
\begin{align}\label{4.94}
\eta^{\psi}(\rho^{\varepsilon},m^{\varepsilon})_t+q^{\psi}(\rho^{\varepsilon},m^{\varepsilon})_x
\qquad\, \mbox{is uniformly bounded in $W_{\rm loc}^{-1, p_3}(\mathbb{R}_+^2)$ for $p_3>2$}.
\end{align}
Then, using \eqref{4.90}--\eqref{4.94} and the interpolation compactness
theorem $(${\it cf}. {\rm  \cite{Ding1985,Ding1989}$)$}, we conclude Lemma $\ref{lemma4.8}$.

\subsection{Proof of Theorem \ref{thm3.3} for the polytropic equation of state}\label{section4.3}

Recall the following compactness theorem established by Chen-Perepelitsa \cite{Perepelitsa}:
\begin{theorem}[Chen-Perepelitsa \cite{Perepelitsa}]\label{th4.9}
Let $(\eta^{\psi},q^{\psi})$ be a weak
entropy pair generated by $\psi\in C_0^2(\mathbb{R})$.
Assume that the sequence
$(\rho^{\varepsilon},u^{\varepsilon})(t,x)$ defined on
$\mathbb{R}_+\times\mathbb{R}$ with
$m^{\varepsilon}=\rho^{\varepsilon}u^{\varepsilon}$ satisfies the
following conditions{\rm :}
\begin{enumerate}
\item [(\rmnum{1})] For any $-\infty<r_1<r_2<\infty$ and
$T>0$,
\begin{align}
\int_{0}^{T}\int_{r_1}^{r_2}(\rho^\varepsilon)^{\gamma+1}\,\dd x\dd t\leq
C(r_1,r_2),\nonumber
\end{align}
where $C>0$ is  independent of $\varepsilon$.

\item [(\rmnum{2})] For any set $K\Subset\mathbb{R}$,
\begin{align}
\int_{0}^{T}\int_{K}\big((\rho^\varepsilon)^{\gamma+\theta}+\rho^\varepsilon|u^\varepsilon|^3\big)\,\dd x\dd t\leq
C(K),\nonumber
\end{align}
where $C(K)>0$ is independent of $\varepsilon$.

\item [(\rmnum{3})] The sequence of entropy dissipation measures
\begin{align}
\eta^{\psi}(\rho^{\varepsilon},m^{\varepsilon})
_t+q^{\psi}(\rho^{\varepsilon},m^{\varepsilon})_x \qquad\, \mbox{is
compact in $H_{\rm loc}^{-1}(\mathbb{R}_+^2)$}.\nonumber
\end{align}
\end{enumerate}
Then there exist both a subsequence $($still denoted$)$
$(\rho^{\varepsilon},m^{\varepsilon})(t,x)$ and a vector-valued function $(\rho, m)(t,x)$ such that
\begin{align*}
(\rho^{\varepsilon},m^{\varepsilon})\rightarrow (\rho, m)\qquad \ a.e. \text{ as $\varepsilon\rightarrow0^+$}.
\end{align*}
\end{theorem}

The uniform estimates and compactness properties obtained in \S\ref{section3}--\S\ref{section4} yields that,
for the sequence of solutions $(\rho^\varepsilon, m^\varepsilon)$ satisfying \eqref{1.8}, \eqref{3.1}--\eqref{3.4},
and the compensated compactness framework established in \cite{Perepelitsa} (see Theorem \ref{th4.9}),
there exist both a subsequence (still denoted) $(\rho^{\varepsilon},m^{\varepsilon})(t,x)$
and a vector-valued function $(\rho, m)(t,x)$ such that
\begin{align}
(\rho^{\varepsilon},m^{\varepsilon})\rightarrow (\rho, m)\qquad  \ a.e.\,\, (t,x)\in \R_+\times\R\quad
\text{as $\varepsilon\rightarrow0^+$}.\nonumber
\end{align}
Notice that
$$
|m|^{\frac{3(\gamma+1)}{\gamma+3}}\leq C\big(\rho|u|^3+\rho^{\gamma+1}\big).
$$
Then, using Lemma \ref{lem4.2} and \ref{lem4.5}, we obtain
\begin{align}
(\rho^{\v},m^{\v})\rightarrow(\rho,m) \qquad \text{ in } L^{q_1}_{\rm loc}(\R_+\times\R)\times L^{q_2}_{\rm loc}(\R_+\times\R),\nonumber
\end{align}
for $q_1\in[1,\gamma+1)$ and $q_2\in[1,\frac{3(\gamma+1)}{\gamma+3}).$

Using again Lemma \ref{lem4.2} and \ref{lem4.5}, we have
\begin{align}\label{4.101}
\eta^{\ast}(\rho^{\v},m^{\v})\rightarrow \eta^{\ast}(\rho,m)\qquad \text{ in $L^1_{\rm loc}(\R_+\times\R)$}\quad \text{ as $\varepsilon\rightarrow0^+$}.
\end{align}

Since $\eta^{\ast}$ is a positive convex function, we use
\eqref{3.7}, \eqref{4.1}, \eqref{4.101}, and Fatou's lemma to see that, for all $t_2\geq t_1\geq0,$
$$
\int^{t_2}_{t_1}\int_{\R}\big(\eta^{\ast}(\rho,m)+\frac{1}{2}\rho (W+\frac{1}{2})\ast \rho \big)(t,x)\,\dd x \dd t
\leq(t_2-t_1)\int_{\R}\big(\eta^{\ast}(\rho_0,m_0)+\frac{1}{2}\rho_0 (W+\frac{1}{2})\ast \rho_0\big)\,\dd x.
$$

Then, by the Lebesgue point theorem, we obtain
\begin{align}
\int_{\R}\big(\eta^{\ast}(\rho,m)+\frac{1}{2}\rho W\ast \rho\big)(t,x) \,\dd x
\leq \int_{\R}\big(\eta^{\ast}(\rho_0,m_0)+\frac{1}{2}\rho_0 W\ast \rho_0\big)(x)\,\dd x:=\mathcal{E}_0.\nonumber
\end{align}

We now prove that $(\rho,m)$ is an entropy solution of the Cauchy problem \eqref{1.1} and \eqref{pressure1}--\eqref{1.6}
for the polytropic case.

Let $\Psi(t,x)\in C^1(\R_+\times\R)$ be a function with compact support and
supp$\,\Psi(t,\cdot)\Subset(b^-(t),b^+(t))$ for any $t\in[0,T]$.
Then
$$
\int_{\R^2_+}\big(\rho^{\v}\Psi_t+\rho^{\v}u^{\v}\Psi_x\big)\,\dd x \dd t+\int_{\R}\rho^{\v}_0(x)\Psi(0,x)\,\dd x=0.
$$
By the Lebesgue dominated convergence theorem, through taking limit $\varepsilon\rightarrow 0^+,$ up to a subsequence, we have
$$
\int_{\R^2_+}\big(\rho\Psi_t+\rho u\Psi_x\big)\,\dd x\dd t+\int_{\R}\rho_0(x)\Psi(0,x)\,\dd x=0.
$$
Next, we consider the momentum equation.
First, we have
 \begin{align}\label{4.103}
&\int_{\R^2_+}\Big(\rho^{\v}u^{\v}\Psi_t+\big(\rho^{\v}(u^{\v})^2+P(\rho^{\v})\big)\Psi_x
+\rho^{\v}\big(\lambda u^{\v}+V +\partial_xW\ast\rho^{\v}\big)\Psi\Big)\,\dd x\dd t
+\int_{\R}(\rho^{\v}_0 u^{\v}_0)(x)\Psi(0,x)\,\dd x\nonumber\\
&=\v\int_{\R^2_+}(\rho^{\varepsilon})^{\alpha}u^{\varepsilon}_x\Psi_x\,\dd x\dd t.
\end{align}
Since
\begin{align*}
 \v\Big|\int_{\R^2_+}(\rho^{\varepsilon})^{\alpha}u^{\varepsilon}_x\Psi_x\,\dd x\dd t\Big|& \leq C\sqrt{\varepsilon}\Big(\int^T_0\int_K\v(\rho^{\v})^{\alpha}(u^{\v}_x)^2\,\dd x\dd t\Big)^{\frac{1}{2}}
  \Big(\int^T_0\int_K(\rho^{\v})^{\alpha}\,\dd x\dd t\Big)^{\frac{1}{2}}\nonumber\\
&\leq C(K)
\sqrt{\varepsilon}\rightarrow 0
\qquad \text{ as $\varepsilon\rightarrow0^+$}, 
\end{align*}
then it follow from  \eqref{4.103} that
\begin{align}
\int_{\R^2_+}\Big(m\Psi_t+\big(\frac{m^2}{\rho}+P(\rho)\big)\Psi_x\Big)\,\dd x\dd t
+\int_{\R}m_0(x)\Psi(0,x)\,\dd x=\int_{\R^2_+}(-\lambda m-\rho V +\rho \partial_xW\ast\rho)\Psi \,\dd x\dd t.\nonumber
\end{align}

The verification of entropy inequality is direct. Therefore, the proof of Theorem
\ref{thm3.3} and hence Theorem \ref{thm:merged} for the polytropic case is completed.
$\hfill\Box$

\section{Proof of Theorem \ref{thm3.3} for the General Pressure Law}\label{section5}
This section is dedicated to the essential improvements of the arguments in \S\ref{section3}--\S\ref{section4}
necessary to prove Theorem \ref{thm3.3} for the general pressure case.

\subsection{Properties of the general pressure law and the related internal energy}\label{section5.1}
In this section, we present some useful estimates
involving the general pressure $P(\r)$ with
\eqref{pressure2}--\eqref{pressure4}  and the corresponding internal energy $e(\rho).$

Denote
\begin{equation}\label{5.1}
	k(\rho):=\int_{0}^{\rho}\frac{\sqrt{P'(y)}}{y}\,\mathrm{d} y.
\end{equation}
By direct calculation, we recall the following asymptotic behaviors of $P(\rho)$, $e(\rho)$, and $k(\rho)$.
\begin{lemma}[\cite{Chen2023}, Lemma 3.1]\label{lemA.1}
The constant $\rho_{*}$  in \eqref{pressure3} can be chosen  small enough, and the constant $\rho^{*}$  in \eqref{pressure4}
 large enough, so that the following estimates hold{\rm :}
\begin{enumerate}
\item [(\rmnum{1})] When $\rho\in (0,\rho_{*}]$,
\begin{equation}\label{5.2}
\left\{\begin{aligned}
&\underline{\kappa}_{1}\rho^{\gamma_1}\leq P(\rho)\leq \bar{\kappa}_{1}\rho^{\g_1},\\
&\underline{\kappa}_{1}\gamma_1\rho^{\gamma_1-1}\leq P'(\rho)\leq \bar{\kappa}_{1}\g_1\rho^{\g_1-1},\\
&\underline{\kappa}_{1}\gamma_1(\gamma_1-1)\rho^{\gamma_1-2}\leq P''(\rho)\leq \bar{\kappa}_{1}\gamma_1(\gamma_1-1)\rho^{\g_1-2},
\end{aligned}
\right.
\end{equation}
and when $\rho\in [\rho^{*},\infty)$,
\begin{equation}\label{5.3}
\left\{\begin{aligned}
&\underline{\kappa}_{2}\rho^{\gamma_2}\leq P(\rho)\leq  \bar{\kappa}_{2}\rho^{\g_2},\\
&\underline{\kappa}_{2}\gamma_2\rho^{\gamma_2-1}\leq P'(\rho)\leq  \bar{\kappa}_{2}\g_2\rho^{\g_2-1},\\ &\underline{\kappa}_{2}\gamma_2(\gamma_2-1)\rho^{\gamma_2-2}\leq P''(\rho)\leq  \bar{\kappa}_{2}\gamma_2(\gamma_2-1)\rho^{\g_2-2},
\end{aligned}
\right.
\end{equation}
where we have denoted $\underline{\kappa}_{i}:=(1-\mathfrak{a}_0)\kappa_{i}$
and $\bar{\kappa}_{i}:=(1+\mathfrak{a}_0) \kappa_{i}$ with $\mathfrak{a}_0=\frac{3-\g_1}{2(\g_1+1)}$ and $i=1,2$.

\item [(\rmnum{2})] For $e(\rho)$ and $k(\rho)$, there exists
$C>0$ depending on
$(\gamma_1, \gamma_2,  \k_1,\k_2, \rho_{*}, \rho^{*})$ such that
\begin{align}
	&C^{-1}\rho^{\g_1-1}\leq e(\rho)\leq C\rho^{\gamma_1-1},\quad
	C^{-1}\rho^{\g_1-2}\leq e'(\rho)\leq C\rho^{\g_1-2}
	\quad\,\, \text{ for }\rho\in (0,\rho_{*}],\nonumber\\
	&C^{-1}\rho^{\g_2-1}\leq e(\rho)\leq C\rho^{\gamma_2-1},\quad
	C^{-1}\rho^{\g_2-2}\leq e'(\rho)\leq C\rho^{\g_2-2}
	\quad\,\, \text{ for }\rho\in [\rho^{*},\infty), \label{A.9-1}
\end{align}
and, for $i=0,1$,
\begin{align*}
&\qquad C^{-1}\rho^{\theta_{1}-i}\leq k^{(i)}(\rho)\leq C\rho^{\theta_{1}-i},\,\,
C^{-1}\rho^{\theta_{1}-2}\leq |k''(\rho)|\leq C\rho^{\theta_1-2}\,\,\,\, \text{for } \rho\in (0,\rho_{*}],\\
	&\qquad C^{-1}\rho^{\theta_{2}-i}\leq k^{(i)}(\rho)\leq C\rho^{\theta_{2}-i},
	\,\, C^{-1}\rho^{\theta_{2}-2}\leq |k''(\rho)|\leq C\rho^{\theta_2-2} \,\,\,\,\text{for } \rho\in [\rho^{*},\infty),
\end{align*}
where $\t_{1}=\frac{\g_1-1}{2}$ and $\t_2=\frac{\g_2-1}{2}$.
\end{enumerate}
\end{lemma}

\subsection{Uniform estimates of the approximate solutions}\label{section5.2}
We start with the basic energy estimate.

\begin{lemma}[\bf Basic energy estimate]\label{lem5.2}
For smooth solution $(\rho, u)(t,x)$ of  problem \eqref{1.8} and \eqref{3.1}--\eqref{3.4},
the following estimate holds{\rm :}
\begin{align}\label{5.4}
&\int^{b^+(t)}_{b^-(t)}\Big(\frac{1}{2}\rho u^2+\rho e(\rho)+\frac{1}{2}\rho W\ast\rho\Big)(t,x)\,\dd x
+\int_0^t\int_{b^-(s)}^{b^+(s)} (\v\rho^{\alpha}|u_x|^2-\lambda \rho u^2)(s,x)\,\dd x\dd  s\nonumber\\
&+\frac{1}{2}\int^t_0\int^{b^+(s)}_{b^-(s)}\int^{b^+(s)}_{b^-(s)}\varpi(x-y)\rho(x)\rho(y)|u(y)-u(x)|^2\,\dd y\dd x\dd s
= \mathcal{E}^{\varepsilon}_0\le C_0.
\end{align}
\end{lemma}

The proof is almost the same as Lemma \ref{lem4.1}, with \eqref{4.2} replaced by
\begin{align*}
-P(\rho)u_{\xi}=\frac{P(\rho)}{\rho^2}\rho_{\tau}=e'(\rho)\rho_{\tau}=e(\rho)_{\tau}.
\end{align*}

Using Lemma \ref{5.2}, we obtain
\begin{corollary}\label{cor5.3}
	It follows from \eqref{A.9-1} and {\rm Lemma \ref{lem5.2}} that
\begin{equation}\nonumber
\int_{b^-(t)}^{b^+(t)}\rho^{\g_2}(t,x)\,\mathrm{d}x
 \leq C
 \int_{b^-(t)}^{b^+(t)}\big(\rho+\r e(\r)\big)(t,x)\,\mathrm{d}x\leq C
 \qquad \text{ for $t\geq 0$}.
	\end{equation}	
\end{corollary}

\begin{lemma}[\bf Higher moment estimate]\label{lemma5.3}
The following estimate holds{\rm :}
\begin{align*}
\int^{b^+(t)}_{b^-(t)}x^2\rho(t,x) \,\dd x\leq C.
\end{align*}
\end{lemma}

\noindent{\bf Proof.} Following the same argument as in Lemma \ref{lemma4.2},
it suffices to show that
\begin{align*}
&\Big|\int^T_0\int^{b^+(s)}_{b^-(s)}P(\rho)(s,x)\,\dd x\dd s\Big|\\
&\leq \Big|\int^T_0\int^{b^+(s)}_{b^-(s)}\big(P(\rho)-\rho e(\rho)\big)(s,x)\,\dd x\dd s\Big|
+\Big|\int^T_0\int^{b^+(s)}_{b^-(s)}\big(\rho e(\rho)\big)(s,x)\,\dd x\dd s\Big|\\
&\leq\Big|\int^T_0\int^{b^+(s)}_{b^-(s)}\big(\rho^2 e'(\rho)-\rho e(\rho)\big)(s,x)\,\dd x\dd s\Big|+C
\leq C.
\end{align*}
$\hfill\Box$

For later use, we analyze the behavior of density $\rho$  on the free boundary. It follows from $\eqref{3.8}_1$ and \eqref{3.9} that
\begin{align*}
\rho_{\tau}(\tau,M)=-\frac{1}{\varepsilon}\lr{\frac{\rho P}{\mu(\rho)}}(\tau,M)\leq0.
\end{align*}
This yields that $\rho(\tau, M)\leq \rho_{0}(M)$.

In the Eulerian coordinates, it is equivalent to
\begin{equation}\label{5.8}
	\rho(t,b^+(t))\leq \rho_0(b).
\end{equation}
Moreover, noting that $\rho^{\gamma_1}_0(\pm b)b\leq C_0$ and $b\geq (\rho_{*})^{-\g_1}$,
we see that $\rho(t,b^+(t))\leq \rho_0(b)\leq \rho_{*}$ for all $t\geq 0$.
From  $\eqref{1.8}_1$ and \eqref{5.2}, there exists a positive constant $\tilde{C}$
depending only on $(\g_1, \k_1)$ such that
\begin{equation*}
	\rho_{\tau}(\tau,M)=-\frac{1}{\v}\, \big(\frac{\rho P}{\mu(\rho)}\big)(\tau,M)
	 \geq -\frac{\tilde{C}}{\v}\big(\rho(\tau,M)\big)^{\g_1+1-
\alpha},
\end{equation*}
which implies
\begin{equation*}
\rho(\tau,M)\geq \rho_{0}(M)
\Big(1+\frac{\tilde{C}(\g_1-\alpha)}{\v}\big(\rho_{0}(M)\big)^{\g_1-\alpha}\tau\Big)^{-\frac{1}{\g_1-\alpha}}.
\end{equation*}
Therefore, in the Eulerian coordinates,
\begin{equation}\label{5.9}
	\rho(t,b^+(t))\geq \rho_{0}(b)\Big(1+\frac{\tilde{C}(\g_1-\alpha)}{\v} (\rho_{0}(b))^{\g_1-\alpha}t\Big)^{-\frac{1}{\g_1-\alpha}}
 \qquad \mbox{for $t\geq 0$}.
\end{equation}

Notice that, in the Lagrangian coordinates, $\rho_0(M)=\rho_0(0)$ since $\rho_0(b)=\rho_0(-b)$ in the Eulerian coordinates.
Therefore, by the uniqueness of solutions of the ordinary differential equation:
$\rho_{\tau}(\tau,\cdot)=-\frac{1}{\v}\lr{\frac{\rho P}{\mu(\rho)}}(\tau,\cdot)$ with the same initial data, we conclude that
\begin{equation}\label{5.10}
 \rho(t,b^+(t))=\rho(t,b^-(t)).
 \end{equation}
 This implies that  $P(\rho(t,b^+(t)))=P(\rho(t,b^-(t)))$.
 This indicates, in particular, that the boundary term in \eqref{5.11} below is nonnegative.

\begin{lemma} \label{lem5.5}
Let
$\rho^{\gamma_1}_0(\pm b)b\leq C_0$.
Then, for any given $T>0$ and for any $t\in[0,T]$, the following holds{\rm :}
\begin{align}\label{5.11}
&\varepsilon^2\int^{b^+(t)}_{b^-(t)}\big(\rho^{2\alpha-3}\rho^2_x\big)(t,x)\,\dd x
+\varepsilon\int^t_0\int^{b^+(t)}_{b^-(t)}\big(P'(\rho)\rho^{\alpha-2}\rho^2_x\big)(s,x)\,\dd x\dd s\nonumber\\
&+\frac{1}{\varepsilon}\int^{t}_0\Big(\big(\frac{P'(\rho)P(\rho)\rho}{\mu(\rho)}\big)(s,b^+(s))b^+(s)
 -\big(\frac{P'(\rho)P(\rho)\rho}{\mu(\rho)}\big)(s,b^-(s))b^-(s)\Big)\,\dd s\nonumber\\
&+\big(P(\rho(t,b^+(t)))b^+(t)-P(\rho(t,b^-(t)))b^-(t)\big)\leq C.
\end{align}
\end{lemma}

\noindent{\bf Proof.}
 Similar with the argument in \eqref{4.32}, we obtain
\begin{align}\label{5.12}
&\Big(\frac{(u+\frac{\varepsilon}{\alpha}(\rho^{\alpha})_{\xi})^2}{2}+e(\rho)+\frac{1}{2}\rho W\ast\rho\Big)_{\tau}
+\big(P(\rho)u\big)_{\xi}+\varepsilon P'(\rho)\rho^{\alpha-1}(\rho_{\xi})^2\nonumber\\
&+\big(u+\frac{\varepsilon}{\alpha}(\rho^{\alpha})_{\xi}\big)\big(-\lambda u-V+\rho(W\ast\rho)_{\xi}\big)=0.
\end{align}
Using $\eqref{3.8}_1$ and \eqref{3.9}, we have
\begin{align}\nonumber
\rho_{\tau}(\tau,M)=-\frac{1}{\v}\Big(\frac{\rho P}{\mu(\rho)}\Big)(\tau,M),
\qquad \rho_{\tau}(\tau,0)=-\frac{1}{\v}\Big(\frac{\rho P}{\mu(\rho)}\Big)(\tau,0),
\end{align}
so that
\begin{align}
&\big(P(\rho)u\big)(\tau,M)-\big(P(\rho)u\big)(\tau,0)\nonumber\\
&=P(\rho)(\tau,M)\frac{\dd}{\dd \tau}b^+(\tau)-P(\rho)(\tau,0)\frac{\dd}{\dd \tau}b^-(\tau)\nonumber\\
&=\big(P(\rho)(\tau,M)b^+(\tau)\big)_{\tau}
  -\big(P(\rho)(\tau,0)b^-(\tau)\big)_{\tau}\nonumber\\
&\quad\, +\Big(P'(\rho)(\tau,0)\rho_{\tau}(\tau,0)b^{-}(\tau)-P'(\rho)(\tau,M)\rho_{\tau}(\tau,M)b^{+}(\tau)\Big)\nonumber\\
&=\big(P(\rho)(\tau,M)b^+(\tau)\big)_{\tau}-\big(P(\rho)(\tau,0)b^-(\tau))_{\tau}\nonumber\\
&\quad\, +\frac{1}{\v}\Big(\frac{P'(\rho)P(\rho)\rho}{\mu(\rho)}\Big)(\tau,M)b^+(\tau)
  -\frac{1}{\v}\Big(\frac{P'(\rho)P(\rho)\rho}{\mu(\rho)}\Big)(\tau,0)b^-(\tau).\nonumber
\end{align}

Integrating \eqref{5.12} over $[0,M]$, we arrive at
\begin{align}\label{5.14}
&\frac{\dd}{\dd \tau}\int^{M}_0\Big(\frac{1}{2}\big(u+\frac{\varepsilon}{\alpha}(\rho^{\alpha})_{\xi}\big)^2+e(\rho)\Big)\,\dd \xi
+\varepsilon\int^M_0P'(\rho)\rho^{\alpha-1}(\rho_{\xi})^2\,\dd \xi \nonumber\\
&\,\,\,+\big(P(\rho)(\tau,M)b^+(\tau)-P(\rho)(\tau,0)b^-(\tau)\big)_{\tau}\nonumber\\
&\,\,\,+\frac{1}{\v}\Big(\frac{P'(\rho)P(\rho)\rho}{\mu(\rho)}\Big)(\tau,M)b^+(\tau)
  -\frac{1}{\v}\Big(\frac{P'(\rho)P(\rho)\rho}{\mu(\rho)}\Big)(\tau,0)b^-(\tau)\nonumber\\
&=\int^M_0\big(u+\frac{\varepsilon}{\alpha}(\rho^{\alpha})_{\xi}\big)\big(-\lambda u-V+\rho(W\ast\rho)_{\xi}\big)\,\dd \xi.
\end{align}
Integrating \eqref{5.14} over $[0,\tau],$ we have
\begin{align}\label{5.142}
&\int^M_0\Big(\frac{1}{2}\big(u+\frac{\varepsilon}{\alpha}(\rho^{\alpha})_{\xi}\big)^2+e(\rho)\Big)\,\dd \xi
+\varepsilon\int^{\tau}_0\int^M_0P'(\rho)\rho^{\alpha-1}(\rho_{\xi})^2\,\dd \xi \dd s\nonumber\\
&\quad+\Big(P(\rho)(\tau,M)b^+(\tau)-P(\rho)(\tau,0)b^-(\tau)\Big)\nonumber\\
&\quad+\frac{1}{\varepsilon}\int^{\tau}_0\Big(\Big(\frac{P'(\rho)P(\rho)\rho}{\mu(\rho)}\Big)(\tau,M)b^+(\tau)
  -\Big(\frac{P'(\rho)P(\rho)\rho}{\mu(\rho)}\Big)(\tau,0)b^-(\tau)\Big)\,\dd s\nonumber\\
&\leq C\int^M_0\Big(\frac{1}{2}\big(u_0+\frac{\varepsilon}{\alpha}(\rho^{\alpha}_0)_{\xi}\big)^2
  +e(\rho_0)\Big)\,\dd \xi+C\big(P(\rho_0(M))+P(\rho_0(0))\big)\, b+C.
\end{align}
Using $b\geq\rho^{-\gamma_1}_{\ast},$ we see that $P(\rho_0(\pm b))\, b\leq C_0$ so that the second term
on the RHS of \eqref{5.142} is uniformly bounded.
This completes the proof.
$\hfill\Box$

\medskip
Motivated by \cite{Chen2021}, to take the limit: $b\rightarrow\infty$,
we need to make sure that domain $\Omega_T=[b^-(t),b^+(t)]$ can expand to the whole physical space $\R$.

\begin{lemma}\label{lem5.6}
Let $0\leq\alpha<\gamma_1,$ $p>\frac{\gamma_1}{\gamma_1-\alpha},$ and $b:=\varepsilon^{-p}$. Then, given any $T>0$,
there exists $\v_0=\v_0(\alpha,\gamma_1,\gamma_2,T)>0$
such that, for any $\v\in(0,\v_0]$,
$$
\pm b^\pm(t)\geq \frac12 b  \qquad\,\, \mbox{for $t\in[0,T]$}.
$$
\end{lemma}

\noindent{\bf Proof.}
We divide the proof into three steps:

\smallskip
\noindent {\emph{Step 1.}} Using \eqref{5.1} and \eqref{5.11}, we have
\begin{align}\label{5.16}
\v\rho^{\beta}(t,x)&=\v\big(\rho^{\beta}(t,x)-\rho^{\beta}(t,b^-(t))+\rho^{\beta}(t,b^-(t))\big)\nonumber\\
&\leq \v\beta\int^{b^+(t)}_{b^-(t)}\rho^{\beta-1}|\rho_x|\,\dd x+\v\rho^{\beta}_0(b)\nonumber\\
&\leq\beta\Big(\v^2\int^{b^+(t)}_{b^-(t)}\rho^{2\alpha-3}\rho^2_x\,\dd x\Big)^{\frac{1}{2}}\Big(\int^{b^+(t)}_{b^-(t)}\rho^{2(\beta-
\alpha)+1}\,\dd x\Big)^{\frac{1}{2}}+\v\rho^{\beta}_0(b)\nonumber\\
&\leq C(\mathcal{E}^{\v}_1)\Big(\int^{b^+(t)}_{b^-(t)}\rho^{2(\beta-
\alpha)+1}\,\dd x\Big)^{\frac{1}{2}}+C\le C,
\end{align}
where we have used \eqref{5.8} and $\beta=\alpha+\frac{\gamma_2-1}{2}.$

\smallskip
\noindent {\emph{Step 2.}} It follows from the boundary condition \eqref{3.2} that
\begin{equation}\label{5.17}
b^+(t)-b^-(t)=2b+\int^{t}_0\big(u(s,b^+(s))-u^-(s,b^-(s))\big)\,\dd s.
\end{equation}

Using \eqref{5.10}, we obtain
\begin{align}
&\big|u(t,b^+(t))-u(t,b^-(t))\big|\nonumber\\
&=\frac{1}{\rho^l(t,b^+(t))}\big|(\rho^lu)(t,b^+(t))-(\rho^lu)(t,b^-(t))\big|
=\frac{1}{\rho^l(t,b^+(t))}\Big|\int^{b^+(t)}_{b^-(t)}(\rho^l u)_x\,\dd x\Big|\nonumber\\
&=\frac{1}{\rho^l(t,b^+(t))}\Big|\int^{b^+(t)}_{b^-(t)}\big(l \rho^{l-1}\rho_xu+\rho^l u_x\big)\,\dd x\Big|\nonumber\\
&\leq\frac{1}{\rho^l(t,b^+(t))}\bigg\{l\Big(\int^{b^+(t)}_{b^-(t)}\v P'(\rho)\rho^{\alpha-2}\rho^2_x\,\dd x\Big)^{\frac{1}{2}}\Big(\int^{b^+(t)}_{b^-(t)}\varepsilon^{-1}\frac{\rho^{2l-\alpha}}{P'(\rho)}u^2\,\dd x\Big)^{\frac{1}{2}}\nonumber\\
&\qquad\qquad\qquad\quad +\Big(\int^{b^+(t)}_{b^-(t)}\v\rho^{\alpha}u^2_x\,\dd x\Big)^{\frac{1}{2}}
  \Big(\int^{b^+(t)}_{b^-(t)}\v^{-1}\rho^{2l-\alpha}\,\dd x\Big)^{\frac{1}{2}}\bigg\}.\label{5.18}
\end{align}
Taking $2l-(\alpha+\gamma_2)+1=1,$ $\it{i.e.,}$ $l=\frac{\alpha+\gamma_2}{2},$ so that $2l-\alpha=\gamma_2.$
We also have
\begin{align}\label{5.19}
\Big|\int^{b^+(t)}_{b^-(t)}\v^{-1}\frac{\rho^{2l-\alpha}}{P'(\rho)}u^2\,\dd x\Big|
&\leq\Big|\int^{b^+(t)}_{b^-(t)}\v^{-1}\frac{\rho^{2l-\alpha}}{P'(\rho)}u^2\mathbf{I}_{\{\rho\leq\rho_{\ast}\}}\,\dd x\Big|+\Big|\int^{b^+(t)}_{b^-(t)}\v^{-1}\frac{\rho^{2l-\alpha}}{P'(\rho)}u^2\mathbf{I}_{\{\rho_{\ast}\leq\rho\leq\rho^{\ast}\}}\,\dd x\Big|\nonumber\\
&\quad\,\,+\Big|\int^{b^+(t)}_{b^-(t)}\v^{-1}\frac{\rho^{2l-\alpha}}{P'(\rho)}u^2\mathbf{I}_{\{\rho\geq\rho^{\ast}\}}\,\dd x\Big|\nonumber\\
&\leq C(\rho^{\ast},\rho_{\ast}).
\end{align}
It follows from \eqref{5.18}--\eqref{5.19}, \eqref{3.7}, and \eqref{5.4} that
\begin{align}\nonumber
|u(t,b^+(t))-u(t,b^-(t))|
\leq\frac{C\v^{-\frac{1}{2}}}{\rho^l(t,b^+(t))}
 \bigg\{\Big(\int^{b^+(t)}_{b^-(t)}\v P'(\rho)\rho^{\alpha-2}\rho^2_x\,\dd x\Big)^{\frac{1}{2}}+\Big(\int^{b^+(t)}_{b^-(t)}\v\rho^{\alpha}u^2_x\,\dd x\Big)^{\frac{1}{2}}\bigg\}.
\end{align}
Using \eqref{5.1} and \eqref{5.11}, we obtain
\begin{align}\label{5.20}
&\int^t_0\big|u(s,b^+(s)-u(s,b^-(s))\big|\,\dd s\nonumber\\
&\leq C\v^{-\frac{1}{2}}\Big(\int^t_0(\rho(s,b^+(s)))^{-\alpha-\gamma_2}\,\dd s\Big)^{\frac{1}{2}}\nonumber\\
&\quad\,\times \bigg\{\Big(\int^t_0\int^{b^+(s)}_{b^-(s)}\v P'(\rho)\rho^{\alpha-2}\rho^2_x\,\dd x\dd s\Big)^{\frac{1}{2}}
+\Big(\int^t_0\int^{b^+(s)}_{b^-(s)}\v\rho^{\alpha}u^2_x\,\dd x\dd s\Big)^{\frac{1}{2}}\bigg\}\nonumber\\
&\leq C(\v^{-\frac{1}{2}}\Big(\int^t_0(\rho(s,b^+(s)))^{-\alpha-\gamma_2}\,\dd s\Big)^{\frac{1}{2}}.
\end{align}
Since $\rho_0(b)\leq Cb^{-\frac{1}{\gamma_1}},$ we take $b\gtrsim(\frac{1}{\v})^{\frac{\gamma_1}{\gamma_1-\alpha}}$ to obtain
\begin{align}\label{5.21}
\frac{\tilde{C}(\gamma_1-\alpha)}{\v}(\rho_0(b))^{\gamma_1-\alpha}\leq\frac{C\tilde{C}(\gamma_1-\alpha)}{\v}b^{-\frac{\gamma_1-\alpha}{\gamma_1}}\leq C.
\end{align}
It follows from \eqref{5.9} and \eqref{5.21} that, for $0\leq s\leq T,$
\begin{align*}
\rho(s,b^+(s))\geq \rho_{0}(b)\Big(1+\frac{\tilde{C}(\g_1-\alpha)}{\v} (\rho_{0}(b))^{\g_1-\alpha}s\Big)^{-\frac{1}{\g_1-\alpha}}\geq\frac{1}{C}(1+T)^{-\frac{1}{\gamma_1-\alpha}}\rho_0(b).
\end{align*}
Then we have
\begin{align}\label{5.23}
\int^t_0(\rho(s,b^+(s)))^{-\alpha-\gamma_2}\,\dd s\leq C(1+T)^{\frac{\alpha+\gamma_2}{\gamma_1-\alpha}}(\rho_0(b))^{-\alpha-\gamma_2}.
\end{align}

Take $b=\v^{-p}$ with $p>\frac{\gamma_1}{\gamma_1-\alpha}$. For any given $T\geq1$,
choose $\v_1:=\big(C_1(1+T)^{\frac{\alpha+\gamma_2}{2(\gamma_1-\alpha)}}\big)^{-\frac{2\gamma_1}{p(\gamma_1-\alpha)-\gamma_1}},$
where $C_1\geq1$ is a large constant depending only on $C_0$.
Then, for any $\v\in(0,\v_1],$ it follows from \eqref{5.20} and \eqref{5.23} that
\begin{align*}
\int^t_0\big|u(s,b^+(s))-u(s,b^-(s))\big|\,\dd s
&\leq C_1\v^{-\frac{1}{2}}(1+T)^{\frac{\alpha+\gamma_2}{2(\gamma_1-\alpha)}}(\rho_0(b))^{-\frac{\alpha+\gamma_2}{2}}\nonumber\\
&=C_1(1+T)^{\frac{\alpha+\gamma_2}{2(\gamma_1-\alpha)}}\v^{-\frac{1}{2}-\frac{p(\alpha+\gamma_2)}{2\gamma_1}}\nonumber\\
&\leq b\, C_1(1+T)^{\frac{\alpha+\gamma_2}{2(\gamma_1-\alpha)}}\v^{p-\frac{1}{2}-\frac{p(\alpha+\gamma_2)}{2\gamma_1}}\nonumber\\
&\leq b\, (1+T)^{\frac{\alpha+\gamma_2}{2(\gamma_1-\alpha)}}\v^{\frac{p(\gamma_1-\alpha)-\gamma_1}{2\gamma_1}}_1\leq b.\nonumber\\
\end{align*}

Using \eqref{5.17}, we conclude
that
\begin{equation}\label{5.24}
b\leq b^+(t)-b^-(t)\leq3b.
\end{equation}

\smallskip
\noindent {\emph{Step 3.}} There exists $x_0(t)\in(b^-(t),b^+(t))$ such that
$$
(\rho u^2)(t,x_0(t))\,\big(b^+(t)-b^-(t)\big)
=\int^{b^+(t)}_{b^-(t)}(\rho u^2)(t,x)\,\dd x\leq C.
$$
Using \eqref{5.24}, we have
\begin{align}\label{5.25}
(\rho u^2)(t,x_0(t))\leq  C(\mathcal{E}^{\v}_0)\,\big(b^+(t)-b^-(t)\big)^{-1}\leq  C)b^{-1}.
\end{align}

It follows from \eqref{3.2} that
\begin{align}\label{5.26}
b^+(t)=b+\int^t_0u(s,b^+(s))\,\dd s=b+\int^t_0\frac{1}{\rho^l(s,b^+(s))}(\rho^lu)(s,b^+(s))\,\dd s,
\end{align}
where $l=\frac{\alpha+\gamma_2}{2}.$
Similar to those as in \eqref{5.18}, we use  \eqref{5.16} and \eqref{5.25} to obtain
\begin{align}\label{5.27}
(\rho^lu)(s,b^+(s))&=(\rho^lu)(s,b^+(s))-(\rho^lu)(s,x_0(s))+(\rho^l u)(s,x_0(s))\nonumber\\
&=\int^{b^+(s)}_{b^-(s)}\big(l\rho^{l-1}\rho_xu+\rho^lu_x\big)\,\dd x+(\rho^lu)(s,x_0(s))\nonumber\\
&\leq\Big(\int^{b^+(s)}_{b^-(s)}P'(\rho)\rho^{\alpha-2}\rho^2_x\,\dd x\Big)^{\frac{1}{2}}\Big(\int^{b^+(s)}_{b^-(s)}\frac{\rho^{2l-\alpha}}{P'(\rho)}u^2\,\dd x\Big)^{\frac{1}{2}}\nonumber\\
&\quad\,+\Big(\int^{b^+(s)}_{b^-(s)}\rho^{\alpha}u^2_x\,\dd x\Big)^{\frac{1}{2}}\Big(\int^{b^+(s)}_{b^-(s)}\rho^{2l-\alpha}\,\dd x\Big)^{\frac{1}{2}}
+C
b^{-\frac{1}{2}}\,(\rho(s,x_0(s)))^{\frac{\alpha+\gamma_2-1}{2}}\nonumber\\
&\leq C
\v^{-\frac{1}{2}}\bigg\{\Big(\v\int^{b^+(s)}_{b^-(s)}P'(\rho)\rho^{\alpha-2}\rho^2_x\,\dd x\Big)^{\frac{1}{2}}+\Big(\v\int^{b^+(s)}_{b^-(s)}\rho^{\alpha}u^2_x\,\dd x\Big)^{\frac{1}{2}}\bigg\}
+C
\v^{-1}b^{-\frac{1}{2}}.
\end{align}

Using \eqref{5.4}, \eqref{5.11}, \eqref{5.23}, and \eqref{5.27}, we have
\begin{align}\label{5.28}
&\Big|\int^t_0\frac{1}{\rho^l(s,b^+(s))}(\rho^lu)(s,b^+(s))\,\dd s\Big|\nonumber\\
&\leq C\v^{-\frac{1}{2}}\Big(\int^t_0(\rho(s,b^+(s)))^{-(\alpha+\gamma_2)}\,\dd s\Big)^{\frac{1}{2}}\nonumber\\
&\quad\,\,\,\times
\bigg\{\Big(\int^t_0\int^{b^+(s)}_{b^-(s)}P'(\rho)\rho^{\alpha-2}\rho^2_x\,\dd x \dd s\Big)^{\frac{1}{2}}
+\Big(\int^t_0\int^{b^+(s)}_{b^-(s)}\v\rho^{\alpha}u^2_x\,\dd x\dd s\Big)^{\frac{1}{2}}\bigg\}\nonumber\\
&\quad +C\v^{-1}b^{-\frac{1}{2}}\int^t_0(\rho(s,b^+(s)))^{-\frac{\alpha+\gamma_2}{2}}\,\dd s\nonumber\\
&\leq C
  \bigg\{\v^{-\frac{1}{2}}\Big(\int^t_0(\rho(s,b^+(s)))^{-(\alpha+\gamma_2)}\,\dd s\Big)^{\frac{1}{2}}
    +\v^{-1}b^{-\frac{1}{2}}\int^t_0(\rho(s,b^+(s)))^{-\frac{\gamma_2+\alpha}{2}}\,\dd s\bigg\}\nonumber\\
&\leq b\,C
\Big(\v^{-\frac{1}{2}}b^{-1}(1+T)^{\frac{\alpha+\gamma_2}{2(\gamma_1-\alpha)}}
  +\v^{-1}b^{-\frac{3}{2}}(1+T)^{\frac{3\gamma_2-\alpha}{2(\gamma_1-\alpha)}}\Big)\, b^{\frac{\alpha+\gamma_2}{2\gamma_1}}\nonumber\\
&\leq b\, C_2(1+T)^{\frac{3\gamma_2-\alpha}{2(\gamma_1-\alpha)}}\v^{\frac{3}{2}p-1-p\frac{\alpha+\gamma_2}{2\gamma_1}}\nonumber\\
&\leq b\, C_2(1+T)^{\frac{3\gamma_2-\alpha}{2(\gamma_1-\alpha)}}\v^{\frac{p(\gamma_1-\alpha)-\gamma_1}{2\gamma_1}}
\leq \frac{1}{2}b
\end{align}
for any $\v\in(0,\v_2]$ with
$\v_2:=\big(2C_2(1+T)^{\frac{3\gamma_2-\alpha}{2(\gamma_1-\alpha)}}\big)^{-\frac{2\gamma_1}{p(\gamma_1-\alpha)-\gamma_1}}.$

Take $\tilde{C_0}:=2\max\{C_1,C_2\}$ and
$\v_0:=\big(\tilde{C}_0(1+T)^{\frac{3\gamma_2-\alpha}{2(\gamma_1-\alpha)}}\big)^{-\frac{2\gamma_1}{p(\gamma_1-\alpha)-\gamma_1}}.$
Then, for any $\v\in(0,\v_0]$, we obtain from \eqref{5.26} and \eqref{5.28} that
$
b^+(t)\geq\frac{1}{2}b.
$

Similarly, it can be proved that $b^-(t)\leq-\frac{1}{2}b.$
This completes the proof.
$\hfill\Box$

\begin{lemma}[\bf Higher integrability of the density]\label{lemma5.6}
Let $(\r,u)$ be the smooth solution of \eqref{1.8} and \eqref{3.1}--\eqref{3.4}.
Then,  under the assumption of {\rm Lemma \ref{lem5.6}},
\begin{align}\label{5.30}
\int_0^T\int_K (\r P(\rho))(t,x)\,\dd x \dd t\leq C(K)
\qquad\mbox{for any $K\Subset[b^-(t),b^+(t)]\,\,$ for any $t \in[0,T]$}.
\end{align}
\end{lemma}

\noindent{\bf Proof.}  We divide the proof into two steps.

\smallskip
\noindent {\emph{Step 1.}}  For given $K\Subset[b^-(t),b^+(t)]$ for any $t\in[0,T]$,
there exist $r_1$ and $r_2$ such that $K\Subset (r_1,r_2)\Subset[b^-(t),b^+(t)]$.
Let $w(x)$ be a smooth, compactly supported function with $\text{supp}\,w\subseteq(r_1,r_2)$ and $w(x)=1$ for $x\in K$.
Multiplying $\eqref{1.8}_2$ by $w(x)$, we have
\begin{align}\label{5.31}
&(\r u w)_t+\big((\r u^2+P(\rho))w\big)_x\nonumber\\
&=(\rho u^2+P(\rho))w_x+\varepsilon(\rho^{\alpha}w u_x)_x-\varepsilon\rho^{\alpha}u_xw_x+\lambda \rho uw+\rho Vw-\rho\partial_xW\ast\rho w.
\end{align}
Integrating \eqref{5.31} over $[r_1,x)$ to obtain
\begin{align}\label{5.32}
(\rho u^2+P(\rho))w
&=\varepsilon\rho^{\alpha}wu_x+\int^{x}_{r_1}\Big((\rho u^2+P(\rho))w_y-\varepsilon
\rho^{\alpha}u_yw_y\Big)\,\dd y-\frac{\dd}{\dd t}\int^x_{r_1}\rho uw\,\dd y\nonumber\\
&\quad+\int^x_{r_1}\lambda \rho uw\,\dd y-\int^x_{r_1}\rho \partial_x W\ast\rho w\,\dd y+\int^x_{r_1}\rho Vw\,\dd y .
\end{align}
Multiplying \eqref{5.32} by $\rho w$ and performing a direct calculation, we have
\begin{align}\label{5.33}
\r P(\rho) w^2&=\v \rho^{\alpha+1}w^2u_x-\Big(\rho w\int^x_{r_1}\rho uw\,\dd y\Big)_t-\Big(\rho uw\int^x_{r_1}\rho uw\,\dd y\Big)_x\nonumber\\
&\quad+\rho uw_x\int^x_{r_1}\rho uw\,\dd y+\rho w\int^x_{r_1}\Big((\rho u^2+P(\rho))w_y-\varepsilon\rho^{\alpha}u_yw_y\Big)\,\dd y\nonumber\\
&\quad+\lambda \rho w\int^x_{r_1}\rho u w\,\dd y-\rho w\int^x_{r_1}\rho\partial_xW\ast\rho w\,\dd y+\rho w\int^x_{r_1}\rho Vw\,\dd y\nonumber\\
&=\sum^{8}_{i=1}K_{i}.
\end{align}

\noindent {\emph{Step 2.}} To estimate $K_i, i=1,\cdots, 8$, in \eqref{5.33}, we first notice that
\begin{align}\label{5.34}
\int_{b^{-}(t)}^{b^+(t)}\rho|u|\, \dd x\leq \int_{b^-(t)}^{b^+(t)} (\rho+\rho u^2)\,\dd x\leq C.
\end{align}
Then it follows from \eqref{5.34} that
\begin{align}\label{5.35}
\left| \int_0^T\int_{r_1}^{r_2} K_2\, \dd x\dd t \right|&=\left| \int_0^T\int_{r_1}^{r_2} \Big(\rho w\int^x_{r_1}\rho uw\,\dd y\Big)_t\, \dd x\dd t \right|\nonumber\\
&\leq\left|\int^{r_2}_{r_1}\Big(\rho w\int^x_{r_1}\rho u w\dd y\Big)(T,x)\,\dd x\right|+\left|\int^{r_2}_{r_1}\Big(\rho w\int^x_{r_1}\rho u w\dd y\Big)(0,x)\,\dd x\right|\nonumber\\[1mm]
&\leq C.
\end{align}

For $K_1$, we have
\begin{align*}
\left|\int_0^T\int_{r_1}^{r_2} K_1\,\dd x\dd t \right|&=\left|\int^T_0\int^{r_2}_{r_1}\rho^{\alpha+1}\omega^2u_x\,\dd x\dd t\right|\nonumber\\
&\leq\v\int^T_0\int^{r_2}_{r_1}\rho^{\alpha}|u_x|^2\omega^2\,\dd x\dd t+\v\int^T_0\int^{r_2}_{r_1}\rho^{\alpha+2}w^2\,\dd x\dd t\nonumber\\
&\leq C
+\v\int^T_0\int^{r_2}_{r_1}\rho^{\alpha+2}w^2\,\dd x\dd t.
\end{align*}	
We now estimate
$\int_{0}^{T}\int_{r_1}^{r_2}\v\rho^{\alpha+2}\omega^2\,\mathrm{d}x\mathrm{d}t$:
For any fixed $t\in [0,T]$,
denoting
$$
A(t):=\{x\in [r_1,r_2]\,:\, \rho(t,x)\geq \rho^{*}\},
$$
then it follows from \eqref{3.7} that $|A(t)|\leq C(r_1,r_2,\rho^{*})M$.
For any $x\in A(t)$, let $x_{0}$ be the closest point to $r$ so that
$\rho(t,x_{0})=\rho^{*}$ with $|x-x_{0}|\leq |A(t)|\leq C(r_1,r_2,\rho^{*})M$.
Then, for any smooth function $f(\rho)$,
\begin{align}
	\sup_{x\in A(t)}f(\rho(t,x))\omega^{2}(x)&\leq f(\rho(t,x_0))\omega^2(x_0)
     +\Big|\int_{x_0}^{x}\partial_y\big(f(\rho(t,y))\omega^{2}(y)\big)\,\mathrm{d}y\Big|\nonumber\\
	&\leq C(\|\omega\|_{C^1})|f(\rho^{*})|
 +\int_{A(t)}\big|\partial_y\big(f(\rho(t,y))\omega^{2}(y)\big)\big|\,\mathrm{d}y.\nonumber
\end{align}
Recalling  \eqref{5.3} and \eqref{A.9-1},
we notice that $P(\rho)\cong \r^{\g_2}$ and $e(\rho)\cong \r^{\g_2-1}$ for any $x\in A(t)$.
Then
\begin{align*}
	&\v\int_{0}^{T}\int_{r_1}^{r_2}\rho^{\alpha+2}\omega^2\,\mathrm{d}x\mathrm{d}t\nonumber\\
	&= \v\int_{0}^{T}\int_{r_1}^{r_2}\rho^{\alpha+2}{\bf I}_{\{\rho\leq \rho^{\ast}\}}\omega^2\,\mathrm{d}x\mathrm{d}t
	+ \v\int_{0}^{T}\int_{r_1}^{r_2}\rho^{\alpha+2}
     {\bf I}_{\{\rho\geq \rho^{\ast}\}}\omega^2\,\mathrm{d}x\mathrm{d}t\nonumber\\
	&\leq C(\rho^{*})
 + C
 \, \v\int_{0}^{T}\Big(\int_{r_1}^{r_2}\rho e(\rho)\dd x\Big)
      \sup_{x\in A(t)}  \Big(\frac{\rho^{\alpha+1}}{e(\rho)}\omega^2\Big)\mathrm{d}t\nonumber\\
	&\leq C(\rho^{*})+ C\, \v\int_{0}^{T}\int_{A(t)}\Big\vert\Big(\frac{\rho^{\alpha+1}}{e(\rho)}\omega^2\Big)_{x}\Big\vert
     \,\mathrm{d}x\mathrm{d}t\nonumber\\
	&\leq  C\,\v\int_{0}^{T}\int_{A(t)}\Big(\big(\frac{(\alpha+1)\rho^{\alpha}}{e(\rho)}-\frac{\rho^{\alpha-1}P(\rho)}{e(\rho)^2}\big)
      |\rho_{x}|\omega^2+\frac{\rho^{\alpha+1}}{e(\rho)}\omega|\omega_{x}|\Big)\,\mathrm{d}x\mathrm{d}t
      +C(\rho^{*}).
\end{align*}
A direct calculation shows that
\begin{align*}
	&\int_{0}^{T}\int_{A(t)}\v\Big(\frac{(\alpha+1)\rho^{\alpha}}{e(\rho)}-\frac{\rho^{\alpha-1}P(\rho)}{e(\rho)^2}\Big)
      |\rho_{x}|\,\omega^2\,\mathrm{d}x\mathrm{d}t\nonumber\\
	&\leq \int_{0}^{T}\int_{A(t)}\v\frac{P'(\rho)}{\rho^{2-\alpha}}|\rho_{x}|^2\omega^2\,\mathrm{d}x\mathrm{d}t
    +\int_{0}^{T}\int_{A(t)}\v\Big(\frac{(\alpha+1)\rho^{\alpha}}{e(\rho)}-\frac{\rho^{\alpha-1}P(\rho)}{e(\rho)^2}\Big)^2\frac{\rho^{2-\alpha}}{P'(\rho)}\omega^2\,\mathrm{d}x\mathrm{d}t
	\nonumber\\
	&\leq C+\int_{0}^{T}\int_{A(t)}\v\rho^{5+\alpha-3\gamma_2}\omega^2\,\mathrm{d}x\mathrm{d}t\nonumber\\
&\leq
    \begin{cases}
	C+\v\int^T_0\int^{r_2}_{r_1}(\rho^{\beta}+1)\,\dd x\dd t \quad &\text{if } 5+\alpha-3\gamma_2\leq\beta,\\[1mm]
		C
  +\frac{\v}{2}\int^T_0\int^{r_2}_{r_1}\rho^{\alpha+2}\omega^2\,\dd x\dd t\,\,\,&\text{if }5+\alpha-3\gamma_2\geq\beta.
	\end{cases}
\end{align*}
We also have
\begin{align*}
&\int^T_0\int_{A(t)}\v\frac{\rho^{\alpha+1}}{e(\rho)}\omega|\omega_x|\,\dd x\dd t\nonumber\\
&\leq\int^T_0\int^{r_2}_{r_1}\v\rho^{\alpha-\gamma_1+2}\omega\,\dd x\dd t
\leq\frac{\v}{2}\int^T_0\int^{r_2}_{r_1}\rho^{\alpha+2}\omega^2\,\dd x\dd t+C(r_1,r_2,\rho^{\ast}).
\end{align*}

Notice that
\begin{align*}
&\left| \int_0^T\int_{r_1}^{r_2}K_3\, \dd x\dd t \right|=\left| \int_0^T\int_{r_1}^{r_2} \Big(\rho uw\int^x_{r_1}\rho uw\dd y\Big)_x\, \dd x\dd t \right|=0,\\
&\left| \int_0^T\int_{r_1}^{r_2} K_4\, \dd x\dd t \right|=\left| \int_0^T\int_{r_1}^{r_2} \rho uw_x\Big(\int^x_{r_1}\rho uw\dd y\Big)\, \dd x\dd t \right|\leq C.
\end{align*}

Using the fact that $\alpha\leq\gamma_2,$ we have
\begin{align*}
\left| \int_0^T\int_{r_1}^{r_2} K_5\, \dd x\dd t \right|
&=\left|\int_0^T\int_{r_1}^{r_2} \Big(\rho uw\int^x_{r_1}\big((\rho u^2+P(\rho))w_y-\v\rho^{\alpha}w_yu_y\big)\dd y\Big)\,
\dd x\dd t \right|\nonumber\\
&\leq C
+\left|\v\int^T_0\int^{r_2}_{r_1}\rho w\Big(\int^x_{r_1}\rho^{\alpha}w_yu_y\,\dd y\Big)\dd x\dd t\right|\nonumber\\
&\leq C(r_1,r_2),
\end{align*}
where we have used the fact that
\begin{align*}
\left|\v\int_0^T\int_{r_1}^{r_2} \rho w\Big(\int^x_{r_1}\rho^{\alpha}w_yu_y\, \dd y\Big)\dd x\dd t \right|
&\leq \v\int^T_0\int^{r_2}_{r_1}|\rho^{\alpha}w_xu_x|\,\dd x\dd t\nonumber\\
&\leq \Big(\v\int^T_0\int^{r_2}_{r_1}\rho^{\alpha}|u_x|^2\,\dd x \dd t\Big)^{\frac{1}{2}}\Big(\v\int^T_0\int^{r_2}_{r_1}\rho^{\alpha}|w_x|^2\,\dd x\dd t\Big)^{\frac{1}{2}}\nonumber\\
&\leq C
\Big(\v\int^T_0\int^{r_2}_{r_1}\rho^{\alpha}\,\dd x\dd t\Big)^{\frac{1}{2}}\nonumber\\
&\leq C
\Big(\v\int^T_0\int^{r_2}_{r_1}(\rho^{\gamma_2}+1)\,\dd x\dd t\Big)^{\frac{1}{2}}\nonumber\\
&\leq C(r_1,r_2).
\end{align*}

We also have
\begin{align}
\left|\int_0^T\int_{r_1}^{r_2} K_6\,\dd x\dd t \right|&\leq C,
\nonumber
\end{align}
\begin{align*}
&\left|\int_0^T\int_{r_1}^{r_2} K_7\,\dd x\dd t \right|\nonumber\\
&=\left|\int_0^T\int_{r_1}^{r_2} \rho(t,x) w(x)\Big(\int^x_{r_1}\rho(t,y) \partial_xW\ast\rho w(y) \,\dd y\Big)\dd x\dd t \right|\nonumber\\
&= \left|\int_0^T\int_{r_1}^{r_2} \rho(t,x) w(x)\Big(\int^x_{r_1}\rho(t,y) \big(\int^{b^+(t)}_{b^-(t)}(1-2H(x-z)+x-z)\rho(z)\,\dd z\big)w(y)\,\dd y\Big)\dd x\dd t \right|
\leq C,
\nonumber
\end{align*}
and
\begin{align}
&\left|\int_0^T\int_{r_1}^{r_2} K_8\,\dd x\dd t \right|\nonumber\\
&=\left|\int^T_0\int^{r_2}_{r_1}\rho(t,x)w(x)\int^x_{r_1}\rho(t,y)\Big(\int_{b^{-}(t)}^{b^+(t)}\varpi(y-z)(u(z)-u(y))\rho(t,z)\,\dd z\Big)w(y)\,\dd y\dd x\dd t\right|
\leq C.\label{5.452}
\end{align}

Integrating \eqref{5.33} over $[0,T]\times[r_1, r_2],$ and utilizing \eqref{5.35}--\eqref{5.452},
we can obtain \eqref{5.30}. This completes the proof.
$\hfill\Box$

\smallskip
\begin{corollary}\label{cordensity}
Under the assumptions of {\rm Lemma \ref{lem5.6}}, it follows from {\rm  Lemma \ref{lemma5.6}} and \eqref{5.3} that
\begin{equation}\nonumber
\int_{b^-(t)}^{b^+(t)}\rho^{\g_2+1}(t,x)\,\mathrm{d}x
 \leq C
 \int_{b^-(t)}^{b^+(t)}(\rho+\r P(\r))(t,x)\,\mathrm{d}x\leq C(r_1,r_2)
 \qquad \text{ for $t\geq 0$}.
	\end{equation}	
\end{corollary}

\medskip
\subsection{A special entropy pair}
Compared with the polytropic gas case in \cite{Chen2021},
there is no explicit formula of the entropy kernel for the general pressure law
\eqref{pressure2}--\eqref{pressure4}
so that we have to analyze the entropy equation \eqref{2.7} carefully to obtain several desired estimates.
In order to obtain the higher integrability of the velocity, we use a special entropy pair constructed in \cite{Chen2023}
such that $\rho|u|^3$ can be controlled by the entropy flux.
Indeed, such a special entropy $\hat{\eta}(\rho,u)$ is constructed as
\begin{equation*}
\hat{\eta}(\rho,u)=\begin{cases}
 \frac{1}{2}\rho u^2+\rho e(\rho)\quad&\text{for }u\ge k(\rho),\\
-\frac{1}{2}\rho u^2-\rho e(\rho)\,\,\, &\text{for }u\le -k(\rho),
\end{cases}
\end{equation*}
for $k(\rho)=\int_{0}^{\rho}\frac{\sqrt{P'(y)}}{y}\,\mathrm{d}y$
and, in the intermediate region $-k(\rho)\le u\le k(\rho)$,
$\hat{\eta}(\rho,u)$ is the unique solution of the Goursat problem
of the entropy equation \eqref{2.7}:
\begin{equation}\label{5.47}
	\left\{\begin{aligned}
		\dis&\eta_{\rho\rho}-k'(\rho)^2\eta_{uu}=0 \qquad \mbox{for $-k(\rho)\le u\le k(\rho)$},\\
  		\dis&\eta(\rho,u)\vert_{u=\pm k(\rho)}=\pm \big(\frac{1}{2}\rho u^2+\rho e(\rho)\big){\rm .}
	\end{aligned}
	\right.
\end{equation}

Let us recall the following lemma, which will be used in the proof of Lemma \ref{lem5.9}.
\begin{lemma}[\cite{Chen2023}, Lemma 4.1]\label{lem5.8}
The Goursat problem \eqref{5.47} admits a unique solution $\hat{\eta}\in C^2(\R_{+} \times \R)$ such that
 \begin{itemize}
\item [(\rmnum{1})]
$|\hat{\eta}(\rho,u)|\leq C(\rho|u|^2+\rho^{\gamma(\rho)})$ for $(\rho,u)\in \R_{+}\times \R$,
where $\gamma(\rho)=\gamma_1$ if $\rho\in [0,\rho_{*}]$ and $\gamma(\rho)=\gamma_2$ if $\rho\in (\rho_{*},\infty)$.

\smallskip
\item [(\rmnum{2})] If $\hat{\eta}$ is regarded as a function of $(\r, u)$,
\begin{align*}
\qquad |\hat{\eta}_{\rho}(\rho,u)|\leq C(|u|^2+\rho^{2\theta(\rho)}),\quad|\hat{\eta}_{u}(\rho,u)|\leq C(\rho|u|+\rho^{\theta(\rho)+1})
\qquad\text{for }(\rho,u)\in \R_{+}\times \R,
\end{align*}
and, if $\hat{\eta}$ is regarded as a function of $(\r, m)$,
\begin{align*}
\qquad |\hat{\eta}_{\rho}(\r,m)|\leq C(|u|^2+\rho^{2\t(\r)}),\quad|\hat{\eta}_{m}(\r,m)|\leq C(|u|+\rho^{\t(\r)})
\qquad\,\text{for } (\rho,m)\in \R_{+}\times \R,
\end{align*}
 where $\theta(\rho):=\frac{\gamma(\rho)-1}{2}$.

 \smallskip
\item [(\rmnum{3})] If $\hat{\eta}_{m}$ is regarded as a function of $(\r,u)$,
\begin{align*}
\qquad |\hat{\eta}_{m\rho}(\rho,u)|\leq C\rho^{\theta(\rho)-1},\quad |\hat{\eta}_{mu}(\rho,u)|\leq C,
\end{align*}
 and, if $\hat{\eta}_{m}$ is regarded as a function of $(\r,m)$,
$$
|\hat{\eta}_{m\rho}(\rho,m)|\leq C\rho^{\theta(\rho)-1},\quad |\hat{\eta}_{mm}(\rho,m)|\leq C\rho^{-1}.
$$

\item[(\rmnum{4})]  If $\hat{q}$ is the corresponding entropy flux determined by \eqref{entropypair},
then $\hat{q}\in C^2(\R_{+}\times \R)$ and
\begin{align*}
&\hat{q}(\rho,u)=\frac{1}{2}\rho |u|^3\pm \rho u(e(\rho)+\rho e'(\rho))\qquad\quad\text{for }\pm u\geq k(\rho),\\
&|\hat{q}(\rho,u)|\leq C\r^{\g(\rho)+\t(\r)}\qquad\qquad\qquad\qquad\,\,\,\,\,\,\text{for }|u|<k(\rho),\\
&\hat{q}(\rho,u)\geq \frac{1}{2}\rho |u|^3\qquad\qquad\qquad\qquad\qquad\quad\,\,\,\text{for } |u|\geq k(\rho),\\
&|\hat{q}-u\hat{\eta}|\leq C(\rho^{\g(\rho)}|u|+\rho^{\g(\rho)+\t(\rho)})\qquad\qquad\text{for } (\rho,u)\in \R_{+}\times \R.
\end{align*}
\end{itemize}
\end{lemma}

We are now ready to prove the better integrability of the velocity.
\begin{lemma}[\bf Higher integrability of the velocity]\label{lem5.9}
Let $(\r,u)$ be the smooth solution of \eqref{1.8} and \eqref{3.1}--\eqref{3.4}.
Then, under the assumption of {\rm Lemma \ref{lem5.6}},
\begin{align*}
\int_0^T\int_{r_1}^{r_2}  \rho|u|^3(t,x) \,\dd x \dd t\leq C(r_1,r_2)
\end{align*}
for any $(r_1,r_2)\Subset [b^-(t), b^+(t)]$.
\end{lemma}

\noindent{\bf Proof.}  Multiplying $\eqref{1.8}_1$ by $\hat{\eta}_\r$ and $\eqref{1.8}_2$ by $ \hat{\eta}_m$,
we have
\begin{align}\label{5.49}
& \hat{\eta}_t+\hat{q}_x = \hat{\eta}_m\, \Big(\v (\rho^{\alpha} u_x)_x+\lambda \rho u+\rho V-\rho \partial_xW\ast\rho\Big).
\end{align}
A direct calculation shows that
\begin{align}\label{5.50}
\frac{\dd}{\dd t}\int_x^{b^+(t)}\hat{\eta}\,\dd y
& =\hat{\eta}(t,b^+(t)) \frac{\dd}{\dd t}b^+(t)+\int_x^{b^+(t)}\hat{\eta}_t(t,y)\,\dd y= (u\hat{\eta})(t, b^+(t))+\int_x^{b^+(t)} \hat{\eta}_t(t,y)\,\dd y.
\end{align}
Integrating \eqref{5.49} over $[x,b^+(t))$,  we have
\begin{align}\label{5.51}
 \hat{q}(t,x) &=\Big( \int_x^{b^+(t)}\hat{\eta}(t,y)\,\dd y\Big)_t+\big(\hat{q}-u\hat{\eta}\big)(t,b^+(t))\nonumber\\
 &\quad-\v \int_x^{b^+(t)}  \hat{\eta}_m(\rho^{\alpha} u_y)_y\,\dd y-
 \lambda \int_x^{b^+(t)}  \hat{\eta}_m\rho u\,\dd y+\int_x^{b^+(t)}\hat{ \eta}_m \rho\partial W\ast\rho\,\dd y-\int^{b^+(t)}_x\hat{\eta}_m\rho V\,\dd y\nonumber\\
 &:=\sum^6_{i=1}I_i.
\end{align}

\smallskip
We now estimate each term $I_i, i=1,\cdots, 6$, in \eqref{5.51}.
First,
for $I_2$ involving the trace estimates in \eqref{5.51}, it follows from \cite[Lemma 3.6]{Chen2021}, Lemma \ref{lem5.8}
that
\begin{align}\label{5.52}
&\int_0^T\int^{r_2}_{r_1}
\big|(\hat{q}-u\hat{\eta})(t,b^+(t))\big|\,
\dd x\dd t\leq C(r_1,r_2)\int_0^T\big(\rho^{\gamma_1+\theta_1}(t,b^+(t))+(\rho^{\gamma_1}|u|)(t, b^+(t)\big)\,\dd t.
\end{align}

It follows from \eqref{5.8} that
\begin{align}\label{5.53}
&\int_0^T\big(\rho(t,b^+(t))^{\gamma_1+\theta_1}\,\dd t=\int^T_0(\rho_0(b))^{\gamma_1+\theta_1}\,\dd t\leq\rho^{\gamma_1+\theta_1}_{\ast}T\leq C.
\end{align}

It is noted from \eqref{5.4}, \eqref{5.8}, and \eqref{5.11} that
\begin{align}\label{5.54}
&\int_0^T\int^{r_2}_{r_1}
\big(\rho^{\gamma_1}|u|\big)(t,b^+(t))\,
\dd x\dd t\nonumber\\
&=\int^T_0\rho^{\frac{\gamma_1-\alpha}{2}}(t,b^+(t))|(\rho^{\frac{\gamma_1+\alpha}{2}}u)(t,b^+(t))|\,\dd t\nonumber\\
&\leq(\rho_0(b))^{\frac{\gamma_1-\alpha}{2}}\int^T_0\bigg\{\Big(\int^{b^+(t)}_{b^-(t)}P'(\rho)\rho^{\alpha-2}\rho^2_x\,\dd x\Big)^{\frac{1}{2}}+\Big(\int^{b^+(t)}_{b^-(t)}\rho^{\alpha}u^2_x\,\dd x\Big)^{\frac{1}{2}}+\v^{-1}b^{-\frac{1}{2}}\bigg\}\,\dd t\nonumber\\
&\leq C
(\rho_0(b))^{\frac{\gamma_1-\alpha}{2}}\nonumber\\
&\quad\,\,\times \bigg\{\v^{-\frac{1}{2}}b^{-\frac{1}{2}}T+\Big(\varepsilon\int^T_0\int^{b^+(t)}_{b^-(t)}P'(\rho)\rho^{\alpha-2}\rho^2_x\,\dd x\dd t\Big)^{\frac{1}{2}}(\frac{T}{\v})^{\frac{1}{2}}
+\Big(\v\int^T_0\int^{b^+(t)}_{b^-(t)}\rho^{\alpha}u^2_x\,\dd x\dd t\Big)^{\frac{1}{2}}(\frac{T}{\v})^{\frac{1}{2}}\bigg\}\nonumber\\
&\leq C
(\rho_0(b))^{\frac{\gamma_1-\alpha}{2}}
  \big(\v^{-1}b^{-\frac{1}{2}}+\v^{-\frac{1}{2}}\big)
  \leq C
b^{-\frac{\gamma_1-\alpha}{2\gamma_1}}
  \big(\v^{-1}b^{-\frac{1}{2}}+\v^{-\frac{1}{2}}\big)\nonumber\\
&\leq C
\v^{\frac{\gamma_1-\alpha}{2\gamma_1}p-\frac{1}{2}}
\leq C,
\end{align}
where $p>\frac{\gamma_1}{\gamma_1-\alpha}.$
Substituting \eqref{5.53}--\eqref{5.54} into \eqref{5.52}, we obtain
\begin{align}\label{5.55}
&\Big|\int^T_0\int^{r_2}_{r_1}I_2\,\dd x\dd t\Big|=\Big|\int^T_0\int^{r_2}_{r_1}
(\hat{q}-u\hat{\eta})(t,b^+(t))\,
\dd x\dd t\Big|\leq C(r_1,r_2).
\end{align}
For $I_1$, using Lemma \ref{lem5.8}, we have
\begin{align}\label{5.56}
&\Big|\int^T_0\int^{r_2}_{r_1}I_1\,\dd x\dd t\Big|=\Big|\int_0^T\int_{r_1}^{r_2}\Big(\int_x^{b^+(t)}\hat{\eta}(\rho,\rho u)\, \dd y\Big)_t\,\dd x\dd t\Big| \nonumber\\
&\leq \Big|\int_{r_1}^{r_2}\int_{b^-(t)}^{b^+(t)}\hat{\eta}(\rho,\rho u)(T,y)\, \dd y\dd x\Big|
  +\Big|\int_{r_1}^{r_2}\int_{b^-(t)}^{b^+(t)}\hat{\eta}(\rho_0,\rho_0 u_0) \,\dd y\dd x\Big|
  \leq C(r_1,r_2).
\end{align}

For $I_3,$ we integrate by parts to obtain
\begin{align}\label{5.57}
-\v\int_x^{b^+(t)}\hat{\eta}_m(\rho^{\alpha} u_y)_y\,\dd y
&=-\v\Big(\hat{\eta}_m(t,b^+(t))\, (\rho^{\alpha} u_x)(t,b^+(t))-\hat{\eta}_m(t,x)(\rho^{\alpha} u_x)(t,x) \Big) \nonumber\\
&\quad\,\,+\v\int_x^{b^+(t)}\rho^{\alpha} u_y(\hat{\eta}_{mu}u_y+\hat{\eta}_{m\rho}\rho_y)\,\dd y:= J_1+J_2.
\end{align}
We first estimate $J_2$:
\begin{align*}
|J_2|=&\Big|\v\int^{b^+(t)}_x\Big(\hat{\eta}_{m\rho}\rho^{\alpha}u_y\rho_y+\hat{\eta}_{mu}\rho^{\alpha}u^2_{y}\Big)\,\dd y\Big|\leq C\v\Big|\int^{b^+(t)}_x\rho^{\theta(\rho)+\alpha-1}u_y\rho_y\,\dd y\Big|+\v\int^{b^+(t)}_x\rho^{\alpha}u^2_{y}\,\dd y\nonumber\\
&\leq \v\int^{b^+(t)}_xP'(\rho)\rho^{\alpha-2}\rho^2_y\,\dd y+\v\int^{b^+(t)}_x\frac{\rho^{\alpha+2\theta(\rho)}}{P'(\rho)}u^2_y\,\dd y+\v\int^{b^+(t)}_x\rho^{\alpha}u^2_{y}\,\dd y\nonumber\\
&\leq \v\int^{b^+(t)}_xP'(\rho)\rho^{\alpha-2}\rho^2_y\,\dd y+\v\int^{b^+(t)}_x\rho^{\alpha}u^2_{y}\,\dd y,
\end{align*}
where we have used the fact that
$|\hat{\eta}_{mu}|\leq C$, $|\hat{\eta}_{m\r}|\leq C\r^{\theta(\rho)-1}$,
and
\begin{align*}
&\Big|\v\int^{b^+(t)}_x\frac{\rho^{\alpha+2\theta(\rho)}}{P'(\rho)}u^2_y\,\dd y\Big|\nonumber\\
&\leq\v\int^{b^+(t)}_x\frac{\rho^{\alpha+2\theta(\rho)}}{P'(\rho)}u^2_y{\bf I}_{\{\rho\leq \rho_{\ast}\}}\,\dd y+\v\int^{b^+(t)}_x\frac{\rho^{\alpha+2\theta(\rho)}}{P'(\rho)}u^2_y{\bf I}_{\{\rho_{\ast}\leq\rho\leq \rho^{\ast}\}}\,\dd y+\v\int^{b^+(t)}_x\frac{\rho^{\alpha+2\theta(\rho)}}{P'(\rho)}u^2_y{\bf I}_{\{\rho\geq \rho^{\ast}\}}\,\dd y\nonumber\\
&\leq\v C\int^{b^+(t)}_x\rho^{\alpha+2\theta(\rho)-(\gamma_1-1)}u^2_y\,\dd y+\v C\int^{b^+(t)}_x\rho^{\alpha+2\theta(\rho)-(\gamma_2-1)}u^2_y\,\dd y
\leq \v C\int^{b^+(t)}_x\rho^{\alpha}u^2_y\,\dd y.
\end{align*}
Then we have
\begin{align}\label{5.61}
\int^T_0\int^{r_2}_{r_1}|J_2|\,\dd x\dd t&\leq\varepsilon\int^T_0\int^{r_2}_{r_1}\int^{b^+(t)}_xP'(\rho)\rho^{\alpha-2}\rho^2_y\,\dd y\dd x\dd t+\varepsilon\int^T_0\int^{r_2}_{r_1}\int^{b^+(t)}_x\rho^{\alpha}u^2_y\,\dd y\dd x\dd t\nonumber\\
&\leq C(r_1,r_2).
\end{align}
For $J_1,$ we have
\begin{align}
&\varepsilon\Big|\int^T_0\int^{r_2}_{r_1} \hat{\eta}_m\rho^{\alpha}u_x\,\dd x\dd t\Big|\nonumber\\
&\leq\varepsilon\int^T_0\int^{r_2}_{r_1}\rho^{\alpha}u^2_x\,\dd x\dd t
 +\varepsilon\int^T_0\int^{r_2}_{r_1}\rho^{\alpha}\big(|u|+\rho^{\theta(\rho)}\big)^2\,\dd x\dd t\nonumber\\
&\leq C
+\v\int^T_0\int^{r_2}_{r_1}\rho^{2\theta(\rho)+\alpha}\,\dd x\dd t+\v\int^T_0\int^{r_2}_{r_1}\rho^{\alpha}u^2\,\dd x\dd t.\label{5.62}
\end{align}
To control \eqref{5.62}, we need
\begin{align}\label{5.63}
\v\int^T_0\int^{r_2}_{r_1}\rho^{2\theta(\rho)+\alpha}\,\dd x\dd t
&=\v\int^T_0\int^{r_2}_{r_1}\rho^{\gamma_2-1+\alpha}\,\dd x\dd t\leq C(\mathcal{E}^{\v}_0,\mathcal{E}^{\v}_1)\int^T_0\int^{r_2}_{r_1}\rho^{\frac{\gamma_2-1}{2}}\,\dd x\dd t\nonumber\\
&\leq C
\int^T_0\int^{r_2}_{r_1}\rho^{\gamma_2}\,\dd x\dd t+C(r_1,r_2)
\leq C(r_1,r_2).
\end{align}
Then we have
\begin{align}\label{5.64}
\v\int^T_0\int^{r_2}_{r_1}\rho^{\alpha}u^2\,\dd x\dd t
&\leq \Big(\int^T_0\int^{r_2}_{r_1}\rho|u|^3\,\dd x\dd t\Big)^{\frac{2}{3}}
\,\Big(\int^T_0\int^{r_2}_{r_1}\varepsilon^3\rho^{3\alpha-2}\,\dd x\dd t\Big)^{\frac{1}{3}}\nonumber\\
&\leq \Big(\int^T_0\int^{r_2}_{r_1}\rho|u|^3\,\dd x\dd t\Big)^{\frac{2}{3}}
\,\Big(\int^T_0\int^{r_2}_{r_1}\varepsilon^3\rho^{3\beta}\,\dd x\dd t+C(r_1,r_2)\Big)^{\frac{1}{3}}\nonumber\\
&\leq C(r_1,r_2)
\Big(\int^T_0\int^{r_2}_{r_1}\rho|u|^3\,\dd x\dd t\Big)^{\frac{2}{3}},
\end{align}
where we have used $\alpha\geq\frac{2}{3}.$
Inserting \eqref{5.63} and \eqref{5.64} into \eqref{5.62}, we see that, for $\delta>0$,
\begin{align*}
&\varepsilon\Big|\int^T_0\int^{r_2}_{r_1} \hat{\eta}_m\rho^{\alpha}u_x\,\dd x\dd t\Big|
\leq C(r_1,r_2)
\frac{1}{\delta}+\delta\int^T_0\int^{r_2}_{r_1}\rho|u|^3\,\dd x\dd t.
\end{align*}
Using Lemma \ref{lem5.8} and \eqref{3.3}, we obtain
$$
|\v(\hat{\eta}_m\rho^{\alpha}u_x)(t,b^+(t))|=|(\hat{\eta}_mP)(t,b^+(t))|
\leq C_{\gamma}\big(\rho^{\gamma_1}|u|+\rho^{\gamma_1+\theta_1}\big)(t,b^+(t)).
$$
Similar again to the argument as in \cite{Chen2021,He} yields
\begin{align}\label{5.66}
\int^T_0\int^{r_2}_{r_1}|\v(\hat{\eta}_m\rho^{\alpha}u_x)(t,b^+(t))|\,\dd x\dd t&\leq C(r_1,r_2)
\int^T_0\big(\rho^{\gamma_1}|u|+\rho^{\gamma_1+\theta_1}\big)(t,b^+(t)))\,\dd t\nonumber\\
&\leq C(r_1,r_2).
\end{align}
Combining \eqref{5.72} and \eqref{5.66}, we have that
\begin{align}\label{5.67}
\Big|\int^T_0\int^{r_2}_{r_1}J_1\,\dd x\dd t\Big|&\leq \delta
\int^T_0\int^{r_2}_{r_1}\rho |u|^3\,\dd x\dd t+C(r_1,r_2).
\end{align}
Inserting \eqref{5.61} and \eqref{5.67} into \eqref{5.57}, we obtain
\begin{align}\label{5.68}
\Big|\int^T_0\int^{r_2}_{r_1}I_3\,\dd x\dd t\Big|&\leq \delta
\int^T_0\int^{r_2}_{r_1}\rho  |u|^3\,\dd x\dd t+C(r_1,r_2).
\end{align}
We also obtain the estimates for $I_i, i=4,5,6$:
\begin{align*}
\Big|\int^T_0\int^{r_2}_{r_1}I_4\,\dd x\dd t\Big|
&=\Big|\int^T_0\int^{r_2}_{r_1}\int^{b^+(t)}_x\rho |u|\big(|u|+\rho^{\theta(\rho)}\big) \,\dd y\dd x\dd t\Big|\nonumber\\
&=\Big|\int^T_0\int^{r_2}_{r_1}\rho u^2\,\dd x\dd t\Big|+\Big|\int^T_0\int^{r_2}_{r_1}\rho^{\gamma_2}\,\dd x\dd t\Big|
\leq C(r_1,r_2),
\end{align*}
\begin{align*}
\Big|\int^T_0\int^{r_2}_{r_1}I_5\,\dd x\dd t\Big|&=\Big|\int^T_0\int^{r_2}_{r_1}\int^{b^+(t)}_x\rho \hat{\eta}_m\partial_xW\ast \rho \,\dd y\dd x\dd t\Big|\nonumber\\
&=\Big|\int^T_0\int^{r_2}_{r_1}\int^{b^+(t)}_x\rho\hat{\eta}_m\Big(M-2\int^{y}_{b^-(t)}\rho(t,z)\,\dd z+yM-\int^{b^+(t)}_{b^{-}(t)}z\rho(z) \,\dd z\Big)\,\dd y\dd x\dd t\Big|\nonumber\\
&\leq C(r_1,r_2),
\end{align*}
and
\begin{align}\label{5.702}
\Big|\int^T_0\int^{r_2}_{r_1}I_6\,\dd x\dd t\Big|&=\Big|\int^T_0\int^{r_2}_{r_1}\int^{b^+(t)}_x\hat{\eta}_m\rho V\,\dd y\,\dd x\dd t\Big|\nonumber\\
&=\Big|\int^T_0\int^{r_2}_{r_1}\int^{b^+(t)}_x\hat{\eta}_m\rho \Big(\int^{b^+(t)}_{b^-(t)}\varpi(y-z)(u(y)-u(z))\rho(t,z)\,\dd z\Big)\,\dd y\,\dd x\dd t\Big|\nonumber\\
&\leq C(r_1,r_2).
\end{align}

Combining estimates \eqref{5.55}--\eqref{5.56} and \eqref{5.68}--\eqref{5.702} leads to
$$
\int_{0}^{T}\int_{r_1}^{r_2}\hat{q}\,\mathrm{d}x\mathrm{d}t\leq C(r_1, r_2),
$$
which gives
\begin{align}\label{5.71}
\int_{0}^{T}\int_{[r_1, r_2]\cap \{x: |u|\geq k(\rho)\}}\rho |u|^3\,\mathrm{d}x\mathrm{d}t
	&\leq 2\int_{0}^{T}\int_{[r_1, r_2]\cap \{x: |u|\geq k(\rho)\}}\hat{q}\,\mathrm{d}x\mathrm{d}t\nonumber\\
 &=2\int_{0}^{T}\int_{r_1 }^{r_2}\hat{q}\,\mathrm{d}x\mathrm{d}t
       -2\int_{0}^{T}\int_{[r_1, r_2]\cap \{x: |u|<k(\rho)\}}\hat{q}\,\mathrm{d}x\mathrm{d}t\nonumber\\
	&\leq C(r_1, r_2)
 +C\int_{0}^{T}\int_{r_1}^{r_2}(\rho+\rho^{\g_2+1})\,\mathrm{d}x\mathrm{d}t\nonumber\\
 &\leq C(r_1, r_2)
\end{align}
by using Lemma \ref{lem4.1}.  On the other hand, we have
\begin{align}\label{5.72}
\int_{0}^{T}\int_{[r_1, r_2]\cap \{x:\,|u|\leq  k(\rho)\}}\rho |u|^3\,\mathrm{d}x\mathrm{d}t
 &\leq C\int_{0}^{T}\int_{r_1}^{r_2}\rho^{\g(\rho)+\t(\rho)}\,\mathrm{d}x\mathrm{d}t
 \nonumber\\	&
 \leq C\int_{0}^{T}\int_{r_1}^{r_2}\big(\rho+\rho P(\rho)\big)\,\mathrm{d}x\mathrm{d}t
 \leq C.
\end{align}
Combining \eqref{5.71} with \eqref{5.72}, we obtain that
$\int_{0}^{T}\int_{r_1}^{r_2}\rho |u|^3\,\mathrm{d}x\mathrm{d}t\leq C(r_1, r_2)$.
This completes the proof.
$\hfill\Box$

\medskip
\subsection{$W^{-1,p}_{\rm loc}(\mathbb{R}_+^2)-$Compactness}\label{section5.4}
In this section, we use the uniform estimates obtained
in \S 5.3 to prove the following key lemma, which
states the
$W^{-1,p}_{\rm loc}(\mathbb{R}_+^2)-$compactness of entropy dissipation measures for the approximate
solution sequence.

\begin{lemma}\label{lemm5.10}
Let $\frac23\leq\alpha\leq\gamma_2$, and let $(\eta^{\psi},q^{\psi})$ be a weak entropy pair
generated by $\psi\in C_0^2(\mathbb{R})$ defined in \eqref{2.8}--\eqref{2.9}.
Then, for the solutions $(\rho^{\varepsilon},u^{\varepsilon})$ with
$m^{\varepsilon}=\rho^{\varepsilon}u^{\varepsilon}$ of CNSEs \eqref{1.8} and \eqref{3.1}--\eqref{3.4},
\begin{align}\label{5.73}
\eta^{\psi}(\rho^{\varepsilon},m^{\varepsilon})_t+q^{\psi}(\rho^{\varepsilon},m^{\varepsilon})_x
\qquad \mbox{is compact in $W^{-1,p}_{\rm loc}(\mathbb{R}_+\times \mathbb{R})\,$ for any $p\in[1,2)$}.
\end{align}
\end{lemma}

\noindent
\textbf{Proof}. To prove this lemma, we first
recall the entropy pair
$(\eta^{\psi},q^{\psi})$ generated by $\psi\in C_0^2(\mathbb{R})$
({\it cf.} \cite{Chen2023}).
Given any $\psi\in C_{0}^2(\R)$, a regular weak entropy pair $(\eta^{\psi},\,q^{\psi})$ is given by
\begin{align}\label{5.74}
	\eta^{\psi}(\r,u)=\int_{\R} \psi(s)\,\chi(\r,u,s)\, \mathrm{d}s,\qquad
	q^{\psi}(\r,u)=\int_{\R}\psi(s)\,\sigma(\r,u,s) \,\mathrm{d}s,
\end{align}
and satisfies the following properties
(\cite{Chen2023}; also see \cite{Chen-LeFloch-2000,Chen-LeFloch-2003,Schrecker-Schulz-2019,Schrecker-Schulz-2020}):
There exists a constant $C_{\psi}>0$ depending only on $\rho^{*}$ and $\psi$
such that
\begin{enumerate}
\item[(i)]
		$|\eta^{\psi}(\r,u)|+
  |q^{\psi}(\rho,u)|\leq C_{\psi}\rho\, $ for all $\rho\in [0,2\rho^{*}]$.
\item[(ii)] If $\eta^{\psi}$ is regarded as a function of $(\rho,m)$, then
\begin{align*}
|\eta_{m}^{\psi}(\rho,m)|+|\r\eta_{mm}^{\psi}(\rho,m)|\leq C_{\psi},\qquad
|\eta_{\r}^{\psi}(\r,m)|\leq C_{\psi}(1+\r^{\t_1}).
\end{align*}
\item[(iii)] If $\eta_m^{\psi}$ is regarded as a function of $(\r,u)$, then
\begin{equation}\nonumber
		|\eta_{mu}^{\psi}(\r,u)|+|\rho^{1-\t_1}\eta_{m\rho}^{\psi}(\rho,u)|\leq C_{\psi}.
\end{equation}
\end{enumerate}

A direct computation on $\eqref{1.8}_1\times\eta^{\psi}_{\rho}(\rho^{\varepsilon},m^{\varepsilon})
+\eqref{1.8}_2\times\eta^{\psi}_m(\rho^{\varepsilon},m^{\varepsilon})$
gives
\begin{align}
\displaystyle\eta^{\psi}(\rho^{\varepsilon},m^{\varepsilon})_t+q^{\psi}(\rho^{\varepsilon},m^{\varepsilon})_x
=&\,\varepsilon\Big(\eta^{\psi}_m(\rho^{\varepsilon},m^{\varepsilon})(\rho^{\varepsilon})^{\alpha}u_x^{\varepsilon}\Big)_x
-\varepsilon\eta^{\psi}_{mu}(\rho^{\varepsilon},m^{\varepsilon})(\rho^{\varepsilon})^{\alpha}(u_x^{\varepsilon})^2\nonumber\\
&\,-\varepsilon\eta^{\psi}_{m\rho}(\rho^{\varepsilon},m^{\varepsilon})(\rho^{\varepsilon})^{\alpha}\rho_x^{\varepsilon}u_x^{\varepsilon}+\eta^{\psi}_m(\lambda\rho^{\v}u^{\v}+\rho^{\v}V-\rho^{\v}\partial_xW\ast\rho^{\v}).\nonumber
\end{align}

Let $K\Subset[b^-(t),b^+(t)]$ be compact. Using properties (ii)-(iii) of the weak entropy pair $(\eta^{\psi},q^{\psi})$
and
the Cauchy-Schwartz
inequality, we have
\begin{equation*}
\begin{aligned}
&\displaystyle\varepsilon\int_{0}^{T}\int_{K}|\eta^{\psi}_{mu}(\rho^{\varepsilon},m^{\varepsilon})(\rho^{\varepsilon})^{\alpha}
(u_x^{\varepsilon})^2
+\eta^{\psi}_{m\rho}(\rho^{\varepsilon},m^{\varepsilon})(\rho^{\varepsilon})^{\alpha}\rho_x^{\varepsilon}u_x^{\varepsilon}|\,\dd x\dd t\\
&\displaystyle\leq
C_{\psi}\varepsilon\int_{0}^{T}\int_{K}(\rho^{\varepsilon})^{\alpha}(u_x^{\varepsilon})^2 \,\dd x\dd t
+C_{\psi}\varepsilon\int_{0}^{T}\int_{K}(\rho^{\varepsilon})^{\theta_2-1+\alpha}(\rho_x^{\varepsilon})^2 \,\dd x\dd t
\leq C(K),
\end{aligned}
\end{equation*}
and
\begin{equation*}
\begin{aligned}
&\displaystyle\int_{0}^{T}\int_{K}\Big|\eta^{\psi}_m\Big(\lambda\rho^{\v}u^{\v}+\rho^{\v}V-\rho^{\v}\partial_xW\ast\rho^{\v}\Big)\Big| \,\dd x\dd t\\
&\displaystyle\leq
C_{\psi}\int_{0}^{T}\int_{K}\rho^{\v}(u^{\v})^2+\rho^{\v} \,\dd x\dd t
+C_{\psi}\int_{0}^{T}\int_{K}|\rho^{\v}\partial_xW\ast\rho^{\v}| \,\dd x\dd t
\leq C(K).
\end{aligned}
\end{equation*}
This implies that
\begin{eqnarray*}
-\varepsilon\eta^{\psi}_{mu}(\rho^{\varepsilon},m^{\varepsilon})(\rho^{\varepsilon})^{\alpha}(u_x^{\varepsilon})^2
-\varepsilon\eta^{\psi}_{m\rho}(\rho^{\varepsilon},m^{\varepsilon})(\rho^{\varepsilon})^{\alpha}\rho_x^{\varepsilon}u_x^{\varepsilon}+\eta^{\psi}_m(\lambda\rho^{\v}u^{\v}+\rho^{\v}V-\rho^{\v}\partial_xW\ast\rho^{\v}),
\end{eqnarray*}
is uniformly bounded in $L^1([0,T]\times K)$, and thus it is compact in
$W_{\rm loc}^{-1,p_1}(\mathbb{R}_+^2)$, for $1<p_1<2$.

If $2\alpha\leq\gamma_2+1,$ then we obtain that
\begin{align}\label{5.75}
\v^{\frac{4}{3}}\int^T_0\int_{K}(\rho^{\v})^{2\alpha}\,\dd x\dd t\leq C(K)\v^{\frac{4}{3}}.
\end{align}
If $2\alpha\geq\gamma_2+1,$ $\alpha<\gamma_2,$ it yields that,
\begin{align}\label{5.75b}
\v^{\frac{4}{3}}\int^T_0\int_{K}(\rho^{\v})^{2\alpha}\,\dd x\dd t
&\leq C(K)\v^{\frac{1}{3}}\int^T_0\int_K(\rho^{\v})^{\alpha-\frac{\gamma_2-1}{2}}\,\dd x\dd t\nonumber\\
&\quad+C(K)\v^{\frac{1}{3}}\int^T_0\int_K(\rho^{\v})^{\gamma_2+1}\,\dd x\dd t.
\end{align}
Moreover, it follows from \eqref{5.75} and \eqref{5.75b} that
\begin{equation*}
\begin{aligned}
\displaystyle\int_{0}^{t}\int_{K}\Big(\varepsilon
\eta^{\psi}_m(\rho^{\varepsilon},m^{\varepsilon})
(\rho^{\varepsilon})^{\alpha}u_x^{\varepsilon}\Big)^{\frac{4}{3}}\, \dd x\dd t
&\displaystyle\leq\int_{0}^{t}\int_{K}\varepsilon^{\frac{4}{3}}
(\rho^{\varepsilon})^{^{\frac{4\alpha}{3}}}|u_x^{\varepsilon}|^{\frac{4}{3}}\,\dd x\dd t\\
&\displaystyle\leq C\varepsilon^{\frac{4}{3}}\int_{0}^{t}\int_{K}
(\rho^{\varepsilon})^{\alpha}|u_x^{\varepsilon}|^2\,\dd x\dd t+C\varepsilon^{\frac{4}{3}}\int_{0}^{t}\int_{K}
(\rho^{\varepsilon})^{2\alpha}\,\dd x\dd t\\
&\displaystyle\leq C(K)\varepsilon^{\frac{1}{3}}\rightarrow  0\qquad\text{ as $\varepsilon\rightarrow0^+$}.
\end{aligned}
\end{equation*}
Then, \eqref{4.86} and \eqref{4.89} yield that
\begin{eqnarray}\label{5.77}
\eta^{\psi}(\rho^{\varepsilon},m^{\varepsilon})_t+q^{\psi}(\rho^{\varepsilon},m^{\varepsilon})_x
\ \mbox{is compact in $W_{\rm loc}^{-1,p_2}(\mathbb{R}_+^2)$} \ \
\mbox{for some}\ p_2\in(1,2).
\end{eqnarray}
which also implies that
\begin{equation}\label{5.78}
	\eta^{\psi}(\rho^{\varepsilon},m^{\varepsilon})_t+q^{\psi}(\rho^{\varepsilon},m^{\varepsilon})_x\qquad \text{is compact in
 $W_{\mathrm{loc}}^{-1,p}(\R_{+}^2)$ with $1\leq p\leq p_{2}$}.
\end{equation}

On the other hand, for $\gamma_2\in (1,3),$ using property (i) of the weak entrpy pair $(\eta^{\psi}, q^{\psi})$,
we have
\begin{eqnarray*}
\int^T_0\int_K\big(|\eta^{\psi}(\rho^{\varepsilon},m^{\varepsilon})|^{\gamma_2+1}+|q^{\psi}(\rho^{\varepsilon},m^{\varepsilon})|^2\big)\,\dd x\dd t
\leq C_{\psi}\int^T_0\int_K(\rho^{\v})^{\gamma_2+1}\,\dd x\dd t\leq C(K),
\end{eqnarray*}
so that $\eta^{\psi}(\rho^{\varepsilon},m^{\varepsilon})$ and $q^{\psi}(\rho^{\varepsilon},m^{\varepsilon})$
are uniformly bounded in $L_{\rm loc}^{2}(\mathbb{R}_+^2)$. This yields
that
\begin{eqnarray}\label{5.81}
\eta^{\psi}(\rho^{\varepsilon},m^{\varepsilon})_t+q^{\psi}(\rho^{\varepsilon},m^{\varepsilon})_x
\qquad \mbox{is uniformly bounded in $W_{\rm loc}^{-1, 2}(\mathbb{R}_+^2)$}.
\end{eqnarray}

Then the interpolation compactness theorem ({\it cf.} \cite{Chen2, Ding1985, Ding1989}) indicates that,
 for $p_{2}>1, p_{1} \in\left(p_{2}, \infty\right]$, and $p_{0} \in\left[p_{2}, p_{1}\right)$,
\begin{align}
\big(\text{compact set of }W_{\mathrm{loc}}^{-1, p_{2}}(\mathbb{R}_{+}^{2})\big)
 \cap \big(\text{bounded set of }W_{\mathrm{loc}}^{-1, p_{1}}(\mathbb{R}_{+}^{2})\big)
 \subset \big(\text{compact set of }W_{\mathrm{loc}}^{-1, p_{0}}(\mathbb{R}_{+}^{2})\big),\nonumber
\end{align}
which is a generalization of Murat's lemma in \cite{F. Murat,L. Tartar}.
Combining this theorem for $1<p_{2}<2$ and $p_{1}=2$ with the facts in \eqref{5.81} and \eqref{5.77}, we conclude that
\begin{equation}\label{6.73-2}
\eta^{\psi}(\rho^{\varepsilon},m^{\varepsilon})_t+q^{\psi}(\rho^{\varepsilon},m^{\varepsilon})_x\qquad
 \text{is compact in $W_{\mathrm{loc}}^{-1,p}(\R_{+}^2)$ with $p_2\leq p<2$}.
\end{equation}
Combining \eqref{6.73-2} with \eqref{5.78}, we conclude \eqref{5.73}.
$\hfill\square$

\subsection{Proof of Theorem \ref{thm3.3}}\label{section5.5}

Recall the following $L^p$ compensated compactness theorem established by Chen-Huang-Li-Wang-Wang \cite{Chen2023}{\rm :}
\begin{theorem}[\cite{Chen2023}, Lemma 2.2]\label{thm5.11}
Let $(\rho^{\v},m^{\v})(t,x)=(\rho^\v, \rho^\v u^\v)(t,x)$ be a sequence of measurable functions
with $\rho^{\v}\geq 0$ a.e. on $\R_{+}^2$
satisfying the following two conditions{\rm :}
\begin{enumerate}
\item [(\rmnum{1})] For any $T>0$ and $K\Subset\mathbb{R}_+$, there exists $C(K)>0$ independent of $\v$ such that
$$
\int_0^T\int_K\big((\rho^\v)^{\g_2+1}+\rho^{\v}|u^{\v}|^3\big)\,\mathrm{d}x\mathrm{d}t\le C(K).
$$
\item [(\rmnum{2})]
For any entropy pair $(\eta^\psi,q^\psi)$ defined in \eqref{5.74} with
any smooth function $\psi(s)$ of compact support on $\mathbb{R}$,	
$$
\eta^{\psi}(\rho^{\varepsilon},m^{\varepsilon})_t+q^{\psi}(\rho^{\varepsilon},m^{\varepsilon})_x
\qquad\mbox{is compact in $ W_{\mathrm{loc}}^{-1,1}(\mathbb{R}_{+}^{2})$}.
$$
\end{enumerate}
Then there exists a subsequence $($still denoted$)$ $(\rho^{\v},m^{\v})(t,x)$ and a vector-valued
function $(\rho,m)(t,x)$ such that, as $\v\to 0^+$,
\begin{equation*}
\begin{aligned}
&\rho^{\v}(t,x)\to \rho(t,x)~~ \mbox{in}~L^{q_1}_{\rm loc}(\R_{+}^2)\qquad\,\,\,\,\text{for }q_1\in [1,\g_2+1),\\
& m^{\v}(t,x)\to m(t,x)~~ \mbox{in}~L^{q_2}_{\rm loc}(\R_{+}^2)\qquad \text{for }q_2\in [1,\frac{3(\g_2+1)}{\g_2+3}),
\end{aligned}
\end{equation*}
where $L_{\mathrm{loc}}^{p}(\mathbb{R}_{+}^{2})$ represents $L^{p}([0, T] \times K)$ for any $T>0$
and compact set $K \Subset \mathbb{R}_+$.
\end{theorem}

The uniform estimates and compactness properties obtained in \S\ref{section5.1}--\S\ref{section5.4} yields that
the sequence of solutions
$(\rho^\varepsilon, m^\varepsilon)$ of problem \eqref{1.8} and \eqref{3.1}--\eqref{3.4} satisfies
the compensated compactness framework established in \cite{Chen2023} so that there exist both a subsequence (still denoted)
$(\rho^{\varepsilon},m^{\varepsilon})$ and a vector function $(\rho, m)$ such that
\begin{eqnarray*}
(\rho^{\varepsilon},m^{\varepsilon})\rightarrow (\rho, m)\qquad\,\, a.e.\,\, (t,x)\in \R_+\times\R\,\,\,
\text{ as $\varepsilon\rightarrow0^+$}.
\end{eqnarray*}

Since
$$
|m^{\v}|^{\frac{3(\gamma_2+1)}{\gamma_2+3}}\leq C\big(\rho^{\v}|u^{\v}|^3+(\rho^{\v})^{\gamma_2+1}\big),
$$
it follows from Lemma \ref{lem5.9} and Corollary \ref{cordensity} that
\begin{eqnarray*}
(\rho^{\v},m^{\v})\rightarrow(\rho,m)\qquad  \text{ in $L^{q_1}_{\rm loc}(\R_+\times\R)\times L^{q_2}_{\rm loc}(\R_+\times\R)$}
\end{eqnarray*}
for $q_1\in[1,\gamma_2+1)$ and $q_2\in[1,\frac{3(\gamma_2+1)}{\gamma_2+3}).$

Moreover, We have
\begin{align*}
\eta^{\ast}(\rho^{\v},m^{\v})\rightarrow \eta^{\ast}(\rho,m)\qquad
\text{ in $L^1_{\rm loc}(\R_+\times\R)$\quad as $\varepsilon\rightarrow0^+$}.
\end{align*}

Then it is direct to check that $(\rho,m)$ is a finite-energy entropy solution of the Cauchy problem \eqref{1.1}.
Therefore, the proof of Theorem
\ref{thm3.3} and hence Theorem \ref{thm:merged} is mow completed for the general pressure case.
$\hfill\Box$

\appendix
\section{Construction of the Approximate Initial Data Sequences}\setcounter{equation}{0}\label{SecA}

We now construct the approximate initial data sequences $(\rho^\v_0(x),\rho^\v_0u^{\v}_0(x))$ with desired estimates,
regularity, and boundary compatibility for the polytropic case.
For the general pressure law case, the construction arguments are similar.

\smallskip
\subsection{The initial density}
Let  $J(x)$ be the standard mollification function and
$\dis J_{\delta}(x):=\frac{1}{\delta}J(\frac{x}{\delta})$ for $\delta\in(0,1)$.
For later use, we take $\delta=\v^{\frac12}$ and define  $\tilde{\rho}_{0}^{\v}(x)$ as
\begin{align}\label{A.1}
\tilde{\rho}_{0}^{\v}(x)
:=\Big(\int_{\mathbb{R}} \big((\rho_0\mathbf{I}_{[-b+1,b-1]})(x-y)\big)^{\alpha-\frac{1}{2}}J_{\sqrt{\varepsilon}}(y)\,\dd y
+\v e^{-x^2}\Big)^\frac{2}{2\alpha-1},
\end{align}
where we have denoted $\mathbf{I}_{[-b+1,b-1]}$ to be the characteristic function of $\{x\in\mathbb{R}\,:\, -b+1\leq|x|\leq b-1\}$.
Since $\alpha>\frac{2}{3}$, then $\tilde{\rho}_{0}^\v(x)\geq  \v^{\frac{2}{2\alpha-1}} e^{-\frac{2}{2\alpha-1}x^2}>0$.

\begin{lemma}\label{lemA.2}
Let $q\in \{1, \gamma, 2\alpha-1\}$ and $\frac{2}{3}<\alpha\leq 1$. Then $\tilde{\rho}_{0}^\v(x)$ satisfy the following properties{\rm :}
There exists $C_0>0$ independent of $\v$, but may depend on $\mathcal{E}_0, M, M_2, \gamma$, and $\alpha$, such that
\begin{align}
&\tilde{\rho}^{\v}_0(b)=\tilde{\rho}^{\v}_0(-b)\leq C_0 b^{-\frac{1}{\gamma}}, \label{A.388}\\[1.5mm]
&\|\tilde{\rho}_0^\v\|_{L^q}\leq C_0\big(\|\rho_0\|_{L^q}+ \v^{\frac{2}{2\alpha-1}}\big) \qquad\,\mbox{for $\v\in(0,1]$},\label{A.2}\\[1.5mm]
&\lim_{\v\rightarrow0+} \big(\|\tilde{\rho}_0^\v-\rho_0\|_{L^q}
 +\|(\tilde{\rho}^{\v}_0)^{\alpha-\frac{1}{2}}-(\rho_0)^{\alpha-\frac{1}{2}}\|_{L^{\frac{2q}{2\alpha-1}}}\big)=0,\label{A.3}\\[1.5mm]
&\v^2\int_{\R}\Big|\big((\tilde{\rho}^{\v}_0(x))^{\alpha-\frac{1}{2}}\big)_x\Big|^2\dd x
  \rightarrow0 \qquad\mbox{as $\v\rightarrow0^+$},\label{A.4}\\[1.5mm]
&\int_{\R} x^2|\tilde{\rho}^{\v}_0(x)|\,\dd x\leq C_0\Big(\int_{\R} x^2|\rho_0(x)|\,\dd x+\v\Big),\label{A.5}\\[1.5mm]
&\lim_{\v\rightarrow0+} \int_{\R} x^2|\tilde{\rho}_0^\v(x)-\rho_0(x)|\,\dd x= 0,\label{A.6}\\[1.5mm]
 &\int_{\R}\big|\tilde{\rho}^{\v}_0(x)(W\ast\tilde{\rho}^{\v}_0)(x)-\rho_0(x)(W\ast\rho_0)(x)\big|\,\dd x\rightarrow0 \qquad\text{ as $\varepsilon\rightarrow0^+$}.\label{A.62}\\[1.5mm]
 &\int_{\R}\tilde{\rho}^{\v}_0(x)(W\ast\tilde{\rho}^{\v}_0)(x) \,\dd x\leq C_0\big(1+\v^{\beta})\qquad
 \text{ for $\beta:=\frac{2}{2\alpha-1}\in[2,6)$}.\label{A.63}
\end{align}
\end{lemma}

\noindent{\bf Proof.} We divide the proof into six steps. In the proof below,
for simplicity, 
$C>0$ is a universal constant independent of $\v$, which may depend on $\mathcal{E}_0, M, M_2,\gamma$, and $\alpha$, and 
may be different at each occurrence so that $C_0>0$ can be chosen eventually, depending only on this $C>0$ for the statement in this lemma.

\medskip
\noindent{\emph{Step 1.}} The proof of \eqref{A.2}--\eqref{A.4} are similar as in \cite{Chen2021}, so we omit them for brevity.
Here, we remark that $\frac{2}{3}<\alpha\leq1,$ then $\beta:=\frac{2}{2\alpha-1}\in [2,\,6)$. Since in the proof of \eqref{A.4}, using Young's inequality, one has
\begin{align}\label{A.4222}
\v^2\int_{\R}\Big|\big((\tilde{\rho}^{\v}_0(x))^{\alpha-\frac{1}{2}}\big)_x\Big|^2\,\dd x&\leq  C\v\Big(\|\rho_0(x)\|^{2\alpha-1}_{L^{2\alpha-1}}+1\Big) \nonumber\\
&\leq C\v\Big(\int_{\R} (1+x^2)\rho_0(x)\,\dd x+\int_{\R}\frac{1}{(1+x^2)^{\frac{2\alpha-1}{2-2\alpha}}}\,\dd x+1\Big)\nonumber\\
&\leq C_0\v.
\end{align}

\smallskip
\noindent{\emph{Step 2.}} We first prove \eqref{A.388}, which is needed in the proof of Lemmas \ref{lem2.2} and \ref{lem4.5} by \eqref{A.1} since $\v<1$:
\begin{align*}
\tilde{\rho}^{\v}_0(b)=\tilde{\rho}^{\v}_0(-b)=(\v e^{-b^2})^{\frac{2}{2\alpha-1}}\leq C b^{-\frac{1}{\gamma}}.
\end{align*}
For $b=\v^{-p}$ and $0<\v<1,$ it is clear that there exists a positive constant C such that
$$
(\v e^{-\v^{-2p}})^{\frac{2}{2\alpha-1}}\leq C\v^{\frac{p}{\gamma}},
$$
where $C$ is independent of $\v\in (0,1]$.

\smallskip
\noindent{\emph{Step 3.}} We are now ready to prove \eqref{A.5}. By Young's inequality, we have
\begin{align}\label{A.7}
\tilde{\rho}_{0}^{\v}(x)
&=\Big(\int_{\mathbb{R}} \big((\rho_0\mathbf{I}_{[-b+1,b-1]})(x-y)\big)^{\alpha-\frac{1}{2}}J_{\sqrt{\varepsilon}}(y)\,\dd y
+\v e^{-x^2}\Big)^\frac{2}{2\alpha-1}\nonumber\\
&\leq 2^{\frac{2}{2\alpha-1}}\bigg(\Big(\int_{\mathbb{R}}
\big((\rho_0\mathbf{I}_{[-b+1,b-1]})(x-y)\big)^{\alpha-\frac{1}{2}}J_{\sqrt{\varepsilon}}(y)\,\dd y \Big)^\frac{2}{2\alpha-1}
+(\v e^{-x^2})^{\frac{2}{2\alpha-1}}\bigg)\nonumber\\
&\leq 2^{\frac{2}{2\alpha-1}}\Big(\int_{\mathbb{R}} J_{\sqrt{\varepsilon}}(y)\rho_0(x-y)\,\dd y +(\v e^{-x^2})^{\frac{2}{2\alpha-1}}\Big).
\end{align}
Therefore, using \eqref{A.7}, we obtain
$$
x^2\tilde{\rho}_{0}^{\v}(x)\leq Cx^2\Big(\int_{\mathbb{R}} J_{\sqrt{\varepsilon}}(y)\rho_0(x-y)\,\dd y +(\v e^{-x^2})^{\frac{2}{2\alpha-1}}\Big),
$$
so that
$$
\int_{\R} x^2|\tilde{\rho}_{0}^{\v}(x)|\, \dd x
\leq C\Big(\int_{\mathbb{R}}\int_{\mathbb{R}} J_{\sqrt{\varepsilon}}(y)x^2\rho_0(x-y)\,\dd y\dd x +\int_{\mathbb{R}} x^2(\v e^{-x^2})^{\frac{2}{2\alpha-1}}\,\dd x\Big).
$$

Notice that
\begin{align}\label{A.8}
\int_{\mathbb{R}^2} J_{\sqrt{\varepsilon}}(y)x^2\rho_0(x-y)\,\dd y\dd x
&=\int_{\mathbb{R}^2} J_{\sqrt{\varepsilon}}(x-y)x^2\rho_0(y)\,\dd y\dd x\nonumber\\
&=\int_{\mathbb{R}}\Big(\int_{\mathbb{R}}x^2J_{\sqrt{\varepsilon}}(x-y) \,\dd x\Big)\rho_0(y)\,\dd y.
\end{align}
For the integrand in \eqref{A.8}, we have
$$
\int_{\R}x^2J_{\sqrt{\varepsilon}}(x-y)\,\dd x=\int_{\R}(s+y)^2J_{\sqrt{\varepsilon}}(s)\,\dd s
=\int_{\R}(s^2+2sy+y^2)J_{\sqrt{\varepsilon}}(s)\,\dd s.
$$
Since
$$
\int_{\R}sJ_{\sqrt{\varepsilon}}(s)\,\dd s=0,\qquad
\int_{\R}s^2J_{\sqrt{\v}}(s)\,\dd s=\int_{\R}s^2\frac{1}{\sqrt{\varepsilon}}J(\frac{s}{\sqrt{\varepsilon}})\,\dd s
=\v\int_{\R}t^2J(t)\,\dd t\leq C\v,
$$
then, for \eqref{A.8}, we have
$$
\int_{\R^2}J_{\sqrt{\varepsilon}}(y)x^2\rho_0(x-y)\,\dd y\dd x\leq C\Big(\int_{\R}y^2\rho_0(y)\,\dd y+\v\Big).
$$

Therefore, we obtain \eqref{A.5}:
$$
\int_{\R} x^2|\tilde{\rho}^{\v}_0(x)|\,\dd x
\leq C\Big(\int_{\R} x^2|\rho_0(x)|\,\dd x+\v+\v^{\frac{2}{2\alpha-1}}\Big)
\leq C\Big(\int_{\R} x^2|\rho_0(x)|\,\dd x+\v\Big).
$$

\medskip
\noindent{\emph{Step 4.}} We now prove \eqref{A.6} by three three steps.

\medskip
\noindent{\emph{Step 4(i).}} Since $M_2<\infty$, then $x^2 \rho_0(x)\in L^1$.  By the standard properties of mollifier $J_{\sqrt{\v}}$,
we have
$$
\int_{\R} x^2\big|\big((\rho_0\mathbf{I}_{[-b+1,b-1]})\ast J_{\sqrt{\v}}\big)(x)-x^2\rho_0(x)\big|\,\dd x \rightarrow 0
\qquad\text{ as $\varepsilon\rightarrow 0^+$.}
$$

\smallskip
\noindent{\emph{Step 4(ii).}} We write
\begin{align*}
&x^2\big(\big((\rho_0\mathbf{I}_{[-b+1,b-1]})\ast J_{\sqrt{\v}})(x)-\rho_0(x)\big)\\
&=x^2\big(\big((\rho_0\mathbf{I}_{[-b+1,b-1]})\ast J_{\sqrt{\v}})(x)-\big((x^2\rho_0\mathbf{I}_{[-b+1,b-1]})\ast J_{\sqrt{\v}}\big)(x)\\
&\quad+\big((x^2\rho_0\mathbf{I}_{[-b+1,b-1]})\ast J_{\sqrt{\v}}\big)(x)-x^2\rho_0(x).
\end{align*}
Then
\begin{align*}
&\int_{\R} x^2\big|\big((\rho_0\mathbf{I}_{[-b+1,b-1]})\ast J_{\sqrt{\v}})(x)-\rho_0(x)\big|\,\dd x\nonumber\\
&\leq \int_{\R} x^2\big|\big((\rho_0\mathbf{I}_{[-b+1,b-1]})\ast J_{\sqrt{\v}}\big)(x)
   -\big((x^2\rho_0\mathbf{I}_{[-b+1,b-1]})\ast J_{\sqrt{\v}}\big)(x)\big|\,\dd x\nonumber\\
&\quad+\int_{\R} x^2\big|\big((\rho_0\mathbf{I}_{[-b+1,b-1]})\ast J_{\sqrt{\v}})(x)-\rho_0(x)\big|\,\dd x:=I_1+I_2.
\end{align*}
It follows from Step 4(i) that $I_2\rightarrow0$ as $\varepsilon\rightarrow 0^+$.

Now, we prove that $I_1\rightarrow0$  as $\varepsilon\rightarrow 0^+$.
We write
\begin{align*}
I_1&=\int_{\R}\Big|x^2\int_{\R}(\rho_0\mathbf{I}_{[-b+1,b-1]})(y)J_{\sqrt{\v}}(x-y)\,\dd y-\int_{\R}y^2(\rho_0\mathbf{I}_{[-b+1,b-1]})(y)J_{\sqrt{\v}}(x-y)\,\dd y\Big|\,\dd x\nonumber\\
&=\int_{\R}\Big|\int_{\R}(x^2-y^2)(\rho_0\mathbf{I}_{[-b+1,b-1]})(y)J_{\sqrt{\v}}(x-y)\,\dd y\Big|\,\dd x\nonumber\\
&\leq\int_{\R}\Big(\int_{\R}|x^2-y^2|J_{\sqrt{\v}}(x-y)\,\dd x\Big)|\rho_0(y)|\,\dd y.
\end{align*}
When $|y|<1,$ we have
\begin{align*}
\int_{\R}|x^2-y^2|J_{\sqrt{\v}}(x-y)\,\dd x&=\int_{|x-y|<\sqrt{\v}}|x^2-y^2|J_{\sqrt{\v}}(x-y)\,\dd x\nonumber\\
&\leq C\int_{|x-y|<\sqrt{\v}}|x-y|(|x|+|y|)J_{\sqrt{\v}}(x-y)\,\dd x\leq C\sqrt{\v}.
\end{align*}
When $|y|\geq1,$ we have
\begin{align*}
\int_{\R}|x^2-y^2|J_{\sqrt{\v}}(x-y)\,\dd x
&=|y|^2\int_{\R}\big|1-(\frac{x}{y})^2\big|J_{\sqrt{\v}}(x-y)\,\dd x\nonumber\\
&=|y|^2\int_{\R}\big|1-(1+\frac{s}{y})^2\big|J_{\sqrt{\v}}(s)\,\dd s\nonumber\\
&=|y|^2\int_{|s|<\sqrt{\v}}\big|1-(1+\frac{s}{y})^2\big|J_{\sqrt{\v}}(s)\,\dd s
\leq C\sqrt{\v}|y|^2.
\end{align*}
Since $M_2<\infty$, we conclude that
\begin{align*}
\int_{\R} x^2\big|\big((\rho_0\mathbf{I}_{[-b+1,b-1]})\ast J_{\sqrt{\v}}\big)(x)-\rho_0(x))\big|\,\dd x\rightarrow0
\qquad\text{ as $\varepsilon\rightarrow0^+$}.
\end{align*}

\smallskip
\noindent{\emph{Step 4(iii).}} Using the properties of the mollifier, we have
\begin{align*}
&\big\|\big((\rho_0\mathbf{I}_{[-b+1,b-1]})^{\frac{1}{\beta}}\ast J_{\sqrt{\v}}+\v e^{-x^2}\big)^{\beta}
  -(\rho_0\mathbf{I}_{[-b+1,b-1]})\ast J_{\sqrt{\v}}\big\|_{L^1}\nonumber\\
&\leq\big\|\big((\rho_0\mathbf{I}_{[-b+1,b-1]})^{\frac{1}{\beta}}\ast J_{\sqrt{\v}}+\v e^{-x^2}\big)^{\beta}
-\rho_0\mathbf{I}_{[-b+1,b-1]}\big\|_{L^1}\nonumber\\
&\quad+\big\|(\rho_0\mathbf{I}_{[-b+1,b-1]})-(\rho_0\mathbf{I}_{[-b+1,b-1]})\ast J_{\sqrt{\v}}\big\|_{L^1}
\rightarrow0 \qquad \text{ as $\varepsilon\rightarrow0^+$}.
\end{align*}

We now prove
\begin{align*}
\int_{\R}x^2\big|\big((\rho_0\mathbf{I}_{[-b+1,b-1]})^{\frac{1}{\beta}}\ast J_{\sqrt{\v}})(x)
+\v e^{-x^2}\big)^{\beta}-\rho_0(x)\big|\,\dd x\rightarrow 0\qquad \text{ as $\varepsilon\rightarrow0^+$}.
\end{align*}
Notice first that
\begin{align}\label{A.188}
&\int_{\R} x^2\big|\big(((\rho_0\mathbf{I}_{[-b+1,b-1]})^{\frac{1}{\beta}}\ast J_{\sqrt{\v}})(x)+\v e^{-x^2}\big)^{\beta}
 -\rho_0(x)\big|\,\dd x\nonumber\\
&\leq\int_{\R}x^2\big|\big(((\rho_0\mathbf{I}_{[-b+1,b-1]})^{\frac{1}{\beta}}\ast J_{\sqrt{\v}})(x)+\v e^{-x^2}\big)^{\beta}
  -\big((\rho_0\mathbf{I}_{[-b+1,b-1]})\ast J_{\sqrt{\v}}\big)(x)\big|\,\dd x^1\nonumber\\
&\quad+\int_{\R}x^2\big|\big((\rho_0\mathbf{I}_{[-b+1,b-1]})\ast J_{\sqrt{\v}}\big)(x)-\rho_0(x)\big|\, \dd x
\rightarrow0\qquad \,\, \text{ as $\varepsilon\rightarrow0^+$}.
\end{align}
As in
Step 4(ii),
we can prove for the second term on the right-hand side of \eqref{A.188} that
$$
\int_{\R} x^2\big|\big((\rho_0\mathbf{I}_{[-b+1,b-1]})\ast J_{\sqrt{\v}}\big)(x)-\rho_0(x)\big|\,\dd x\rightarrow 0
\qquad \text{ as $\varepsilon\rightarrow0^+$}.
$$
For the first term in \eqref{A.188}, we have
\begin{align*}
&\big((\rho_0\mathbf{I}_{[-b+1,b-1]})^{\frac{1}{\beta}}\ast J_{\sqrt{\v}}+\v e^{-x^2}\big)^{\beta}
-(\rho_0\mathbf{I}_{[-b+1,b-1]})\ast J_{\sqrt{\v}}\nonumber\\
&=\big((\rho_0\mathbf{I}_{[-b+1,b-1]})^{\frac{1}{\beta}}\ast J_{\sqrt{\v}}+\v e^{-x^2}\big)^{\beta}
-\big((\rho_0\mathbf{I}_{[-b+1,b-1]})^{\frac{1}{\beta}}\ast J_{\sqrt{\v}}\big)^{\beta}\nonumber\\
&\quad+\big((\rho_0\mathbf{I}_{[-b+1,b-1]})^{\frac{1}{\beta}}\ast J_{\sqrt{\v}}\big)^{\beta}
 -(\rho_0\mathbf{I}_{[-b+1,b-1]})\ast J_{\sqrt{\v}}.
\end{align*}

Notice that
\begin{align*}
&\Big|\big((\rho_0\mathbf{I}_{[-b+1,b-1]})^{\frac{1}{\beta}}\ast J_{\sqrt{\v}}+\v e^{-x^2}\big)^{\beta}
 -\big((\rho_0\mathbf{I}_{[-b+1,b-1]})^{\frac{1}{\beta}}\ast J_{\sqrt{\v}}\big)^{\beta}\Big|\nonumber\\
&\leq\beta\v e^{-x^2}\big((\rho_0\mathbf{I}_{[-b+1,b-1]})^\frac{1}{\beta}\ast J_{\sqrt{\v}}+\v e^{-x^2}\big)^{\beta-1},
\end{align*}
which leads directly to
\begin{align*}
\lim_{\v\rightarrow0+}
\int_{\R} x^2\big|\big(((\rho_0\mathbf{I}_{[-b+1,b-1]})^{\frac{1}{\beta}}\ast J_{\sqrt{\varepsilon}})(x)+\v e^{-x^2}\big)^{\beta}
-\big(((\rho_0\mathbf{I}_{[-b+1,b-1]})^{\frac{1}{\beta}}\ast J_{\sqrt{\v}})(x)\big)^{\beta}\big|\,\dd x=0.
\end{align*}
Thus, it remains to prove
$$
\lim_{\v\rightarrow0+}
\int_{\R} x^2\big|\big(((\rho_0\mathbf{I}_{[-b+1,b-1]})^{\frac{1}{\beta}}\ast J_{\sqrt{\varepsilon}})(x)\big)^{\beta}
-\rho_0(x)\big|\,\dd x=0.
$$

Since $M_2<\infty$, then, for any $\varepsilon>0$, there exists $N>1$ such that
\begin{align*}
\int_{|x|>N-1}x^2\rho_0(x)\,\dd x<\varepsilon.
\end{align*}
Using the H\"{o}lder inequality and the fact that $\int_{\R}J_{\sqrt{\varepsilon}}(y)\,\dd y=1$ for $0<\v<1$,
we have
\begin{align*}
&\int_{|x|>N}x^2\Big(\int_{\R}(\rho_0\mathbf{I}_{[-b+1,b-1]})^{\frac{1}{\beta}}(x-y)J_{\sqrt{\v}}(y)\,\dd y\Big)^{\beta}\,\dd x\nonumber\\
&\leq \int_{|x|>N}x^2\int_{\R}\rho_0(x-y)J_{\sqrt{\v}}(y)\,\dd y\dd x=\int_{\R}\Big(\int_{|x|>N}x^2\rho_0(x-y)\,\dd x\Big)J_{\sqrt{\v}}(y)\,\dd y\nonumber\\
&=\int_{\R}\Big(\int_{|y+\xi|>N}(y^2+\xi^2+2y\xi)\rho_0(\xi)\,\dd \xi\Big)J_{\sqrt{\v}}(y)\,\dd y\nonumber\\
&\leq2\int_{|y|<\sqrt{\v}}\Big(\int_{|y+\xi|>N}(\v^2+\xi^2)\rho_0(\xi)\,\dd\xi\Big)J_{\sqrt{\v}}(y)\,\dd y< 2(\v^2+\v)\rightarrow 0
\qquad \text{ as $\varepsilon\rightarrow0^+$}.\nonumber
\end{align*}
Therefore, it suffices to show
\begin{align}\label{AA}
\lim_{\v\rightarrow0+}\int_{|x|\leq N}x^2\Big|\Big(\int_{\R}(\rho_0\mathbf{I}_{[-b+1,b-1]})^{\frac{1}{\beta}}(x-y)J_{\sqrt{\v}}(y)\,\dd y\Big)^{\beta}-\rho_0(x)\Big|\,\dd x=0.
\end{align}

Note that
 \begin{align*}
&\int_{|x|\leq N}x^2\Big|\Big(\int_{\R}(\rho_0\mathbf{I}_{[-b+1,b-1]})^{\frac{1}{\beta}}(x-y)J_{\sqrt{\v}}(y)\,\dd y\Big)^{\beta}-\rho_0(x)\Big|\,\dd x\nonumber\\
&\leq N^2\int_{|x|<N}\Big(\int_{|y|<\sqrt{\v}}\rho_0^{\frac{1}{\beta}}(x-y)J_{\sqrt{\v}}(y)\,\dd y\ \Big)^{\beta}-\rho_0(x)\Big|\,\dd x.
\end{align*}
Using H\"{o}lder's inequality, we have
\begin{align*}
\Big(\int_{|y|<\sqrt{\varepsilon}}\rho_0^{\frac{1}{\beta}}(x-y)J_{\sqrt{\v}}(y)\,\dd y\Big)^{\beta}
&\leq \int_{|y|<\sqrt{\v}}\rho_0(x-y)J_{\sqrt{\v}}(y)\,\dd y\nonumber\\
&\leq \sup_{0<\varepsilon<1}\int_{|y|<\sqrt{\v}}\rho_0(x-y)J_{\sqrt{\v}}(y)\,\dd y:=\Phi^{\ast}\rho_0(x).
\end{align*}

Notice that
$$
\Phi^{\ast}\rho_0(x)\leq C\|J\|_{L^1}\mathcal{M}\rho_0(x),
$$
where $\mathcal{M}$ is the Hardy-Littlewood maximal operator and is of strong $(p,p)$ type for $1<p<\infty$.
Then, using Lebesgue's dominated convergence theorem and
$$
\lim_{\v\rightarrow0^+}\big((\rho_0\mathbf{I}_{[-b+1,b-1]})^{\frac{1}{\beta}}\ast J_{\sqrt{\v}})(x)
=\rho_0^{\frac{1}{\beta}}(x),
$$
we conclude \eqref{AA}, which implies \eqref{A.6}.

\medskip
\noindent{\emph{Step 5.}} We now prove \eqref{A.62}. We first note that
 \begin{align*}
&\int_{\R}\Big|\tilde{\rho}^{\v}_0(x)(W\ast\tilde{\rho}^{\v}_0)(x)-\rho_0(x)(W\ast\rho_0)(x)\Big|\,\dd x\nonumber\\
&\leq \int_{\R}|\tilde{\rho}^{\v}_0(x)-\rho_0(x)||(W\ast\tilde{\rho}^{\v}_0)(x)|\,\dd x
+\int_{\R}\rho_0(x)|(W\ast(\tilde{\rho}_0-\rho))(x)|\,\dd x:=I+II.
\end{align*}
Since
\begin{align*}
&|(W\ast\tilde{\rho}^{\v}_0)(x)|=\Big|\int_{\R}W(x-y)\tilde{\rho}^{\v}_0(x)\,\dd x\Big|
\leq \int_{\R}\big(x^2+y^2+\frac{1}{2}\big)\tilde{\rho}^{\v}_0(y)\,\dd y,\\
&\tilde{\rho}^{\v}_0(x)=\Big(\int_{\R}((\rho_0\mathbf{I}_{[-b+1,b-1]})(x-\bar{y}))^{\alpha-\frac{1}{2}}J_{\sqrt{\v}}(\bar{y})\,\dd \bar{y} +\v e^{-x^2}\Big)^{\beta}\nonumber\\
&\qquad\,\,\leq C\Big(\int_{\R}(\rho_0\mathbf{I}_{[-b+1,b-1]})^{\frac{1}{\beta}}(x-\bar{y})J_{\sqrt{\v}}(\bar{y})\,\dd \bar{y}\Big)^{\beta}
+C\big(\v e^{-x^2}\big)^{\beta},
\end{align*}
we have
\begin{align}\label{A.27}
&\int_{\R}\big(x^2+y^2+\frac{1}{2}\big)\tilde{\rho}^{\v}_0(y)\,\dd y\nonumber\\
&=(x^2+\frac{1}{2})\int_{\R}\tilde{\rho}^{\v}_0(y)\,\dd y+\int_{\R}y^2\tilde{\rho}^{\v}_0(y)\,\dd y\nonumber\\
&\leq C(x^2+\frac{1}{2})\Big(\int_{\R}\Big(\int_{\R}(\rho_0\mathbf{I}_{[-b+1,b-1]})^{\frac{1}{\beta}}(y-\bar{y})J_{\sqrt{\v}}(\bar{y})\,\dd \bar{y}\Big)^{\beta}\,\dd y+\int_{\R}(\v e^{-y^2})^{\beta}\,\dd y\Big)\nonumber\\
&\quad+ C\int_{\R}\Big(\int_{\R}|y|^{\frac{2}{\beta}}(\rho_0\mathbf{I}_{[-b+1,b-1]})^{\frac{1}{\beta}}(y-\bar{y})J_{\sqrt{\v}}(\bar{y})\,\dd \bar{y}\Big)^{\beta}\,\dd y+ C\int_{\R}y^2(\v e^{-y^2})^{\beta}\,\dd y\nonumber\\
&\leq C(x^2+\frac{1}{2})\Big(\int_{\R}\Big(\int_{\R}\rho_0(y-\bar{y})(J_{\sqrt{\v}}(\bar{y}))^{\beta}\,\dd y\Big)^{\frac{1}{\beta}}\,\dd \bar{y}\Big)^{\beta}+\v^{\beta}\Big)\nonumber\\
&\quad+ C\Big(\int_{\R}\Big(\int_{\R}|y|^2\rho_0(y-\bar{y})(J_{\sqrt{\v}}(\bar{y}))^{\beta}\,\dd y\Big)^{\frac{1}{\beta}}\,\dd \bar{y}\Big)^{\beta}+C \v^{\beta}\nonumber\\
&\leq CM^{\beta}(x^2+\frac{1}{2})+C\v^{\beta}(x^2+\frac{1}{2})+C\int_{\R}\Big(|\bar{y}|^{\frac{2}{\beta}}M+M_2\Big)J_{\sqrt{\v}}(\bar{y})\,\dd \bar{y}\nonumber\\
&\leq C(x^2+1).
\end{align}

Using \eqref{A.27}, we have
\begin{align}
I&=\int_{\R}|\tilde{\rho}^{\v}_0(x)-\rho_0(x)||(W\ast\tilde{\rho}^{\v}_0)(x)|\,\dd x\nonumber\\
&\leq C\int_{\R}(x^2+1)|\tilde{\rho}^{\v}_0(x)-\rho_0(x)|\,\dd x\rightarrow 0\qquad \text{ as $\varepsilon\rightarrow0^+$}.\label{A.28}
\end{align}

Since
\begin{align*}
\big|(W\ast(\tilde{\rho}^{\v}_0(x)-\rho_0))(x)\big|
&=\Big|\int_{\R}W(x-y)(\tilde{\rho}^{\v}_0(y)-\rho_0(y))\,\dd y\Big|\nonumber\\
&\leq \int_{\R}(x^2+y^2+1)|\tilde{\rho}^{\v}_0(y)-\rho_0(y)|\,\dd y\nonumber\\
&=(x^2+1)\int_{\R}|\tilde{\rho}^{\v}_0(y)-\rho_0(y)|\,\dd y+\int_{\R}y^2|\tilde{\rho}^{\v}_0(y)-\rho_0(y)|\,\dd y,
\end{align*}
then
\begin{align}\label{A.29}
II&=\int_{\R}\rho_0(x)\big|(W\ast(\tilde{\rho}_0-\rho))(x)\big|\,\dd x\nonumber\\
&\leq\int_{\R}(x^2+1)\rho_0(x)\,\dd x\int_{\R}|\tilde{\rho}^{\v}_0(y)-\rho_0(y)|\,\dd y+\int_{\R}\rho_0(x)\,\dd x+
\int_{\R}y^2|\tilde{\rho}^{\v}_0(y)-\rho_0(y)|\,\dd y\nonumber\\
&=(M_2+M)\|\tilde{\rho}^{\v}_0-\rho_0\|_{L^1}+M\|y^2(\tilde{\rho}^{\v}_0(y)-\rho_0(y))\|_{L^1}\rightarrow 0
\qquad \text{ as $\varepsilon\rightarrow0^+$}.
\end{align}

Combining \eqref{A.28} and \eqref{A.29}, we conclude \eqref{A.62}.

\medskip
\noindent{\emph{Step 6.}} We now prove \eqref{A.63}.
For $\tilde{W}:=W+\frac{1}{2}$, then $0\leq\tilde{W}(x-y)\leq x^2+y^2+1.$
Define
$$
\hat{\rho}_0(x):=\Big(\int_{\R}(\rho_0\mathbf{I}_{[-b+1,b-1]})^{\alpha-\frac{1}{2}}(x-y)J_{\sqrt{\v}}(y)\,\dd y\Big)^{\beta}.
$$
Then we have
\begin{align}\label{A.291}
\tilde{\rho}^{\v}_0(x)\leq C\big(\hat{\rho}_0(x)+(\v e^{-x^2})^{\beta}\big).
\end{align}
It follows from \eqref{A.291} that
\begin{align}\label{A.30}
&\int_{\R}\tilde{\rho}^{\v}_0(x)(W\ast\tilde{\rho}^{\v}_0)(x) \,\dd x
\leq\int_{\R}\tilde{\rho}^{\v}_0(x)(\tilde{W}\ast\tilde{\rho}^{\v}_0)(x) \,\dd x\nonumber\\
&\leq C\int_{\R}(\hat{\rho}_0(x)+\v^{\beta}e^{-\beta x^2})\Big(\int_{\R}\tilde{W}(x-y)\hat{\rho}_0(y)\,\dd y+\v^{\beta}\int_{\R}\tilde{W}(x-y)e^{-\beta y^2}\,\dd y\Big)\,\dd x\nonumber\\
&=C\int_{\R}\hat{\rho}_0(x)\int_{\R}\tilde{W}(x-y)\hat{\rho}_0(y)\,\dd y \dd x+C\v^{\beta}\int_{\R}\hat{\rho}_0(x)\int_{\R}\tilde{W}(x-y)e^{\beta y^2}\,\dd y \dd x\nonumber\\
&\quad+C\v^{\beta}\int_{\R}e^{-\beta x^2}\int_{\R}\tilde{W}(x-y)\hat{\rho}_0(y)\,\dd y\dd x+C\v^{2\beta}\int_{\R}e^{-\beta x^2}\int_{\R}\tilde{W}(x-y)e^{-\beta y^2}\,\dd y\dd x\nonumber\\
&=\sum^4_i\mathcal{T}_i.
\end{align}

Using H\"{o}lder's inequality, we have
\begin{align*}
\hat{\rho}_0(x)&=\Big(\int_{\R}(\rho\mathbf{I}_{[-b+1,b-1]})^{\alpha-\frac{1}{2}}(x-y)J^{\frac{1}{\beta}}_{\sqrt{\v}}(y)J^{\frac{1}{\beta'}}_{\sqrt{\v}}(y)\,\dd y\Big)^{\beta}\nonumber\\
&\leq\Big(\int_{\R}\rho_0(x-y)J_{\sqrt{\v}}(y)\,\dd y\Big)\Big(\int_{\R}J_{\sqrt{\v}}(y)\,\dd y\Big)^{\frac{\beta}{\beta'}}\leq\int_{\R}\rho_0(x-y)J_{\sqrt{\v}}(y)\,\dd y,
\end{align*}
where we have used $\frac{1}{\beta}+\frac{1}{\beta'}=1,$ and the property of mollifier: $\int_{\R}J_{\sqrt{\v}}(y)\,\dd y=1.$
Then
$$
\int_{\R}\hat{\rho}_0(x)\,\dd x\leq\int_{\R}\int_{\R}\rho_0(x-y)J_{\sqrt{\v}}(y)\,\dd x\dd y=M,
$$
and
\begin{align*}
\int_{\R}x^2\hat{\rho}_0(x)\,\dd x
&\leq\int_{\R}x^2\int_{\R}\rho_0(x-y)J_{\sqrt{\v}}(y)\,\dd y\dd x=\int_{\R}\Big(\int_{\R}x^2\rho_0(x-y)\,\dd x\Big)J_{\sqrt{\v}}(y)\,\dd y\nonumber\\
&=\int_{\R}\Big(\int_{\R}(y+\bar{x})^2\rho_0(\bar{x})\,\dd \bar{x}\Big)J_{\sqrt{\v}}(y)\,\dd y\nonumber\\
&\leq2\int_{\R}(y^2M+M_2)J_{\sqrt{\v}}(y)\,\dd y\leq C,
\end{align*}
where we have used $\int_{\R}y^2J_{\sqrt{\v}}(y)\,\dd y\leq\v.$

Now, we estimate \eqref{A.30} term by term.
\begin{align*}
\mathcal{T}_1&\leq C\int_{\R}\hat{\rho}_0(x)\int_{\R}(x^2+y^2+1)\hat{\rho}_0(y)\,\dd y\dd x\nonumber\\
&\leq CM\int_{\R}\hat{\rho}_0(x)(x^2+1)\,\dd x+C\int_{\R}\hat{\rho}_0(x)M_2\,\dd x\leq C(M,M_2).\\
\mathcal{T}_2&\leq C\v^{\beta}\int_{\R}\hat{\rho}_0(x)\int_{\R}(x^2+y^2+1)e^{-\beta y^2}\,\dd y\dd x\leq C\v^{\beta}\int_{\R}\hat{\rho}_0(x)(x^2+1)\,\dd x\leq C\v^{\beta}.\\
\mathcal{T}_3&\leq
C\v^{\beta}\int_{\R}e^{-\beta x^2}\int_{\R}(x^2+y^2+1)\hat{\rho}_0(y)\,\dd y\dd x\leq C\v^{\beta}\int_{\R}e^{-\beta x^2}(x^2M+M_2+M)\,\dd x\leq C\v^{\beta}.\\
\mathcal{T}_4&\leq
C\v^{2\beta}\int_{\R}e^{-\beta x^2}\int_{\R}(x^2+y^2+1)e^{-\beta x^2}\,\dd y\dd x\leq C\v^{2\beta}.
\end{align*}
Combining these estimates together, we conclude \eqref{A.63}.
$\hfill\Box$

\medskip
Now we define $\bar{\rho}_{0}^{\v}(x)$ by
\begin{align}\label{A.40}
\bar{\rho}_{0}^{\v}(x)
&:=\tilde{\rho}_{0}^{\v}(x)\,\mathbf{I}_{[-b,b]}(x).
\end{align}

\begin{remark}
We have multiplied by a cut-off function $\mathbf{I}_{[-b,b]}(x)$ in \eqref{A.40}
such that $(W\ast\bar{\rho}_{0}^{\v})
(x)$ is meaningful.
Otherwise, when $|y|>b$, $\bar{\rho}_{0}^{\v}(y)=b^{-\frac{1}{\gamma}}$
so that $(W\ast\bar{\rho}_{0}^{\v})(x)\sim\int_{|y|>b}W(x-y)\,\dd y \, b^{-\frac{1}{\gamma}}=\infty$.
\end{remark}

Then we can obtain the following lemma.
\begin{lemma}\label{lemA.6}
Let $b=\v^{-p}, p>\frac{\gamma}{\gamma-\alpha},$  and $q\in\{1,\gamma,2\alpha-1\}.$
The smooth function $\bar{\rho}_{0}^{\v}(x)$ defined in \eqref{A.40} satisfies the
following properties: There exists $C_0>0$ independent of $\v$ but may depend on $(\mathcal{E}_0, M, M_2, \gamma,\alpha)$ such that
\begin{align*}
&\bar{\rho}^{\v}_0(b)=\bar{\rho}^{\v}_0(-b)\leq C_0 b^{-\frac{1}{\gamma}},\\
&\int^{b}_{-b} |\tilde{\rho}_0^{\v,b}(x)-\bar{\rho}_0^\v(x)|^q\,\dd x
 +\int^{b}_{-b} \Big| (\tilde{\rho}_{0}^{\v}(x))^{\alpha-\frac{1}{2}}
-(\bar{\rho}_{0}^{\v}(x))^{\alpha-\frac{1}{2}}\Big|^{\frac{2q}{2\alpha-1}} \,\dd x\rightarrow0 \qquad\mbox{as $\v\rightarrow 0^+$},\\
&\int^b_{-b}\bar{\rho}^{\v}_0(x)\,\dd x\rightarrow\int_{\R}\rho_0(x)\,\dd x=M\qquad
\text{ as } \v\rightarrow 0^+,\\
&\v^2\int^b_{-b}\Big|\big((\bar{\rho}_0^{\v}(x))^{\alpha-\frac{1}{2}}\big)_x\Big|^2 \,\dd x
\leq C_0\v ,
\\
&\int^{b}_{-b} x^2|\tilde{\rho}_0^{\v}(x)-\bar{\rho}_0^\v(x)|\,\dd x\rightarrow 0
\qquad\,\,\mbox{as $\v\rightarrow 0^+,$}\\
&\int^b_{-b}|\bar{\rho}^{\v}_0(x)W\ast\bar{\rho}^{\v}_0(x)-\tilde{\rho}^{\v}_0(x)W\ast\tilde{\rho}^{\v}_0(x)|\,\dd x\rightarrow 0 \qquad\,\,\text{ as $\varepsilon\rightarrow0^+$}. 
\end{align*}
\end{lemma}

\smallskip
Since $\dis \int^b_{-b} \tilde{\rho}_{0}^{\v}(x)\,\dd x\neq M$ in general, we define
\begin{align}\label{A.23-1}
\rho_0^{\v}(x):=\frac{M }{\int^b_{-b} \bar{\rho}_{0}^{\v}(x)\,\dd x}
  \bar{\rho}_{0}^{\v}(x).
\end{align}
Combining \eqref{A.23-1} with Lemmas \ref{lemA.2} and \ref{lemA.6}, we have

\begin{lemma}
Let $q \in\{1,\gamma,2\alpha-1\}$ and $\frac{2}{3}<\alpha\leq1.$
The smooth function $\rho_{0}^{\v}(x)$ defined in \eqref{A.23-1} satisfies the following properties: There exists $C_0>0$
independent of $\v\in (0,1]$ but may depend on $(\mathcal{E}_0, M, M_2, \gamma,\alpha)$ such that
\begin{align*}
&(\rho^{\v}_0(b))^{\gamma}\, b=(\rho^{\v}_0(-b))^{\gamma}\, b\leq C_0,\\
	&\int^b_{-b} \big|\rho_0^{\v}(x)-\rho_0(x)\big|^q\, \dd x+\int^b_{-b} \Big| (\rho_{0}^{\v}(x))^{\alpha-\frac{1}{2}}
	-(\rho_{0}(x))^{\alpha-\frac{1}{2}}\Big|^{\frac{2q}{2\alpha-1}} \,\dd x\rightarrow0 \qquad\mbox{as $\varepsilon\rightarrow0^+$},\\
	&\v^2\int^b_{-b}\Big|\big((\rho^{\v}_0(x))^{\alpha-\frac{1}{2}}\big)_x\Big|^2 \dd x
	\leq C_0\v,
\\
&\int^{b}_{-b} x^2|\rho_0^{\v}(x)-\rho_0(x)|\,\dd x\rightarrow 0\qquad\mbox{as $\v\rightarrow 0^+,$}\\
&\int^b_{-b}|\rho^{\v}_0(x)(W\ast\rho^{\v}_0)(x)-\rho_0(x)(W\ast\rho_0)(x)|\,\dd x\rightarrow 0
 \qquad\text{ as $\varepsilon\rightarrow0^+$}. 
	\end{align*}
\end{lemma}

\medskip

\subsection{The initial velocity} Next, we construct the approximate initial data for the velocity. We first define
\begin{align}
&\bar{u}_{0}^{\v}(x)
:=\frac{1}{\sqrt{\rho_0^{\v}(x)}}\int_{\mathbb{R}}\Big(\frac{m_0\mathbf{I}_{[-b+2,b-2]}}{\sqrt{\rho_0}}\Big)
(x-y)\,J_{\sqrt{\varepsilon}}(y)\,\dd y,\label{A.31}
\end{align}
for which the following properties are valid.
\begin{lemma}\label{lemA.5}
 Let $\frac{2}{3}<\alpha\leq1$. Then $\bar{u}^\v_0(x)$ defined in \eqref{A.31}  satisfies
\begin{align*}
&{\rm supp}\, \bar{u}^{\v}_0(x) \subset \{x\in \R\, :\,-b+1\leq x\leq b-1\},\\[1.5mm]
&\int_{\mathbb{R}}\rho_0^\v(x) |\bar{u}_0^\v(x)|^2\,\dd x
\equiv\int_{\mathbb{R}}\frac{|m_0(x)|^2}{\rho_0(x)}\,\dd x
  \qquad\,\,\,\mbox{for any $\v\in (0,1]$}, \\[1.5mm]
&\lim_{\v\rightarrow0^+} \|\rho^\v_0 \bar{u}^\v_0-m_0\|_{L^1(\mathbb{R})}=0.
\end{align*}
\end{lemma}

The proof of Lemma \ref{lemA.5} is similar to that in \cite[Lemma A.8]{Chen2023} and \cite[Lemma 5.4]{He},
so we omit it for brevity.

\smallskip
We still need to modify $\bar{u}_{0}^{\v}(x)$ so that it satisfies the stress-free boundary condition \eqref{2.5}
at $x=\pm b$. To this end, we define a cut-off function $S=S(z)\in C^\infty(\mathbb{R})$ that satisfies
\begin{align}\label{A.19-1}
S(z)=\begin{cases}
0\quad   &\mbox{if $z\in(-\infty,0)$},\\
1\quad  &\mbox{if $z\in(2,\infty)$},\\
\mbox{monotonically  increasing}  &\mbox{if $z\in [0,2]$}.
\end{cases}
\end{align}
 Then we define
\begin{equation}\label{A.37}
	u_{0}^{\v}(x):=\bar{u}_{0}^{\v}(x)
	-\frac{1}{\varepsilon}\,S(4(x-(b-\frac{1}{2})))
	\int^b_x\frac{P(\rho^{\v}_0(z))}{\mu(\rho^{\v}_0(z))}\,\dd z-\frac{1}{\varepsilon}\,S(-4(x-(b-\frac{1}{2})))
	\int^{-b}_{x}\frac{P(\rho^{\v}_0(z))}{\mu(\rho^{\v}_0(z))}\,\dd z.
\end{equation}

It is direct to check that  $u_{0}^{\v}(x)\in C^{\infty}([-b,b])$ satisfies
the following boundary conditions{\rm :}
\begin{align*}
u_{0}^{\v}(\pm b)=0,\qquad \big(P((\rho^{\v}_0(x)))-\v\mu(\rho^{\v}_0(x))(u^{\v}_0(x))_x\big)\big|_{x=\pm b}=0.
\end{align*}
A direct calculation following \eqref{A.37} implies
that
\begin{align}\label{A.1111}
&\int^b_{-b}\Big|\sqrt{\rho_{0}^{\v}(x)}u_{0}^{\v}(x)-\sqrt{\rho_{0}^{\v}(x)}\bar{u}_{0}^{\v}(x)\Big|^2 \,\dd x\nonumber\\
&\leq\frac{C}{\varepsilon^2}\Big(\int^{-b+\frac{1}{2}}_{-b}\rho^{\v}_0(x)\big|S(-4x-4b+2)\int^x_{-b}(\rho^{\v}_0(z))^{\gamma-\alpha}\,\dd z\big|^2\,\dd x\nonumber\\
&\qquad\quad +\int^b_{b-\frac{1}{2}}\rho^{\v}_0(x)\big|S(4x-4b+2)\int^b_x(\rho^{\v}_0(z))^{\gamma-\alpha}\,\dd z\big|^2 \,\dd x\Big)\nonumber\\
&\leq\frac{C}{\v^2}\Big(\int^b_{b-\frac{1}{2}}b^{-\frac{1}{^{\gamma}}}\big|\int^b_xb^{-\frac{\gamma-\alpha}{\gamma}}\,\dd z\big|^2\,\dd x+\int^{-b+\frac{1}{2}}_{-b}b^{-\frac{1}{^{\gamma}}}\big|\int^{x}_{-b}b^{-\frac{\gamma-\alpha}{\gamma}}\,\dd z\big|^2\,\dd x\Big)\nonumber\\
&\leq\frac{C}{\v^2}b^{-\frac{2(\gamma-\alpha)+1}{\gamma}}\leq C\v^{{\frac{2(\gamma-\alpha)+1}{\gamma}p}-2}\rightarrow 0 \qquad\text{ as } \varepsilon\rightarrow 0^+,
\end{align}
for $p>\frac{\gamma}{\gamma-\alpha},$ where we have used \eqref{A.19-1} and $b=\v^{-p}.$

Therefore, similarly to those as in \cite{Chen2021,He}, using Lemma \ref{lemA.5} and \eqref{A.1111},
we conclude
\begin{lemma}\label{lemA.7}
For fixed $\v>0$,
\begin{align*}
&\lim_{\varepsilon\rightarrow0^+}\int^b_{-b}\big|\sqrt{\rho_{0}^{\v}(x)}u_{0}^{\v}(x)\big|^2\,\dd x
=\int_{\mathbb{R}} \big|\frac{m_0(x)}{\sqrt{\rho_0(x)}}\big|^2\,\dd x,\\
&\lim_{\varepsilon\rightarrow0^+}\int^b_{-b}|\rho^{\v}_0(x)u^{\v}_0(x)-m_0(x)|\,\dd x=0.
\end{align*}	
\end{lemma}

\smallskip

In the previous argument, we have already defined the cut-off for the density. Now,
with $\rho^{\v}_0(x)$ and $u_0^{\v}(x)$ defined respectively in \eqref{A.23-1} and \eqref{A.37},
we also need to define the cut-off for the velocity function. For $b\gg1$,  define
\begin{align}\label{A.61}
\big(\rho^{\v}_0, u^{\v}_0\big)(x):=\big(\rho_{0}^{\v}(x), u^{\v}_0(x)\big) \mathbf{I}_{[-b,b]}(x).
\end{align}
Then, collecting all the above estimates, we have the following results.

\begin{lemma}\label{propA.1}
Let $\frac{2}{3}<\alpha\leq1,$ $b=\v^{-p},$ $p>\frac{\gamma}{\gamma-\alpha}$,
$q\in \{1,\gamma\}$, and $\v\in(0,1].$
Let $\big(\rho^{\v}_0, u^{\v}_0\big)(x)$ be the functions defined in \eqref{A.61}
so that $(\rho^{\v}_0, u^{\v}_0)\in C^\infty([-b,b]).$
Then
\begin{enumerate}
\item[\rm (i)]  For all $\v\in(0,1]$,
\begin{align*}
&\,\, \, b^{\frac{1}{\gamma}}\rho^{\v}_0(-b)=b^{\frac{1}{\gamma}}\rho^{\v}_0(b)\leq C_0,\quad u^{\v}_0(\pm b)=0,
\quad \big(P(\rho^{\v}_0(x)-\v\mu(\rho^{\v}_0(x)))(u^{\v}_0(x)))_x\big)\big|_{x=\pm b}=0,\\
&\,\,\,\int^b_{-b}\rho^{\v}_0(x)\,\dd x=M, \qquad \mathcal{E}^{\v}_0+\v^{-1}\mathcal{E}^{\v}_1\le C_0,
\end{align*}	
where $C_0>0$  is a constant independent of $\v\in (0,1]$ but may depend on $(\mathcal{E}_0, M, M_2, \gamma, \alpha)$.

\smallskip
\item[\rm (ii)]
As $\varepsilon\rightarrow 0^+$,
\begin{align*}
& (\mathcal{E}_0^{\v}, \mathcal{E}_1^{\v})\rightarrow (\mathcal{E}_0, 0), \\
& \int^{\infty}_{-\infty}x^2\rho^{\v}_0(x)\,\dd x\rightarrow \int^{\infty}_{-\infty}x^2\rho_0 (x) \,\dd x=M_2, \\
&(\rho_0^{\v}, m_0^{\v})\rightarrow (\rho_0, m_0)(x) \qquad\,
  \mbox{\rm in $L^q_{\rm loc}(\R)\times L^1_{\rm loc}(\R)$},
\end{align*}
where $\mathcal{E}_0^{\v}$ and $\mathcal{E}_1^{\v}$ are defined in \eqref{3.5} and \eqref{3.6}, respectively.
\end{enumerate}
\end{lemma}

\medskip
\noindent{\bf Acknowledgments.}
JAC, GQC, and DF acknowledge support by EPSRC grant EP/V051121/1. JAC was also partially supported by the Advanced Grant Nonlocal-CPD (Nonlocal PDEs for Complex Particle Dynamics: Phase Transitions, Patterns and Synchronization) of the European Research Council Executive Agency (ERC) under the European Union Horizon 2020 research and innovation programme (grant agreement No. 883363), and by EPSRC grant EP/T022132/1. GQC was also partially supported by EPSRC grants EP/L015811/1 and EP/V008854. DF was also partially supported by the National Natural Sciences Foundation of China Grant 12001045. The research of EZ was supported by the EPSRC Early Career Fellowship  EP/V000586/1.

\end{document}